\documentclass[review]{elsarticle}
\usepackage[OT1]{fontenc}  	
\usepackage[utf8]{inputenc}
\usepackage[spanish,french,english,USenglish]{babel} 		
\usepackage{amssymb}  			
\usepackage{amsmath}	
\usepackage{amsthm}	
\usepackage{esint}  				
\usepackage{bbm}   				
\usepackage{color,graphicx} 			
\usepackage[small]{caption}
\usepackage{cases}							
\usepackage{mathpazo}

\newcommand{\be}{\begin{equation}}
\newcommand{\ee}{\end{equation}}

\renewcommand{\Im}{\mathop{\rm Im}\nolimits}

\newcommand{\sh}{\mathop{\rm sh}\nolimits}
\newcommand{\ch}{\mathop{\rm ch}\nolimits}

\newcommand{\tg}{\mathop{\rm tg}\nolimits}
\newcommand{\ctg}{\mathop{\rm ctg}\nolimits}
\newcommand{\arctg}{\mathop{\rm arctg}\nolimits}

\makeatletter
\newcommand{\specialnumber}[1]{
\def\tagform@##1{\maketag@@@{(\ignorespaces##1\unskip\@@italiccorr#1)}}}
\makeatother

\makeatletter
\def\ps@pprintTitle{%
\let\@oddhead\@empty
\let\@evenhead\@empty
\let\@oddfoot\@empty
\let\@evenfoot\@oddfoot
}\makeatother

\newtheorem{teorema}{Theorem}
\newtheorem{lemma}{Lemma}
\newtheorem{pf}{Proof}
\newtheorem{corollary}{Corollary}

\journal{Journal of Number Theory (Elsevier)}

\begin{document}

\begin{frontmatter}

\title{A theorem for
the closed-form evaluation of the first generalized Stieltjes constant
at rational arguments and some related summations\footnote{\phantom{-} \\[6mm]
\texttt{\underline{Note to the readers of the 3rd arXiv version:}
this version is a copy of the journal version of the article, which has been published in the Journal of Number Theory (Elsevier), vol.~148, pp.~537-592, 2015.
DOI 10.1016/J.JNT.2014.08.009  http://www.sciencedirect.com/science/article/pii/S0022314X14002820 \\
Artcile history: submitted 14 January 2014, accepted 19 August 2014, published on-line 7 October 2014.\\ 
The layout of the present version and its page numbering differ from the journal version, but the content, the numbering of equations and the numbering 
of references are the same. This version also incorporates some minor corrections to the final journal version, which were published on-line
in the same journal on February 3, 2015 (DOI 10.1016/J.JNT.2015.01.001), and which are mainly due to the publisher's typesetters.
For any further reference to the material published here, please, use the journal version of the paper, 
which you can always get for free by writing a kind e-mail to the author.}}}

\author{Iaroslav V.~Blagouchine\corref{cor1}} 
\ead{iaroslav.blagouchine@univ-tln.fr}
\cortext[cor1]{Corresponding author. Phones: +33--970--46--28--33, +7--953--358--87--23.}
\address{University of Toulon, France.}

\begin{abstract}
Recently, it was conjectured that the first generalized Stieltjes
constant~at rational argument may be always expressed by means of
Euler's constant, the first Stieltjes constant, the $\Gamma$-function
at rational argument(s) and some relatively simple, perhaps even
elementary, function. This conjecture was based on the evaluation of
$\gamma_1(1/2)$, $\gamma_1(1/3)$, $\gamma_1(2/3)$, $\gamma_1(1/4)$,
$\gamma_1(3/4)$, $\gamma_1(1/6)$, $\gamma_1(5/6)$, which could be
expressed in this way. This article completes this previous study and provides 
an elegant theorem which allows to evaluate the first generalized Stieltjes constant at any rational argument.
Several related summation formul\ae~involving the first gener- alized
Stieltjes constant~and the Digamma function are also presented. In
passing, an interesting integral representation for the logarithm of
the $\Gamma$-function at rational argument is also obtained. Finally, it
is shown that similar theorems may be derived for higher Stieltjes
constants~as well; in particular, for the second Stieltjes constant~the
theorem is provided in an explicit form.
\end{abstract}


\end{frontmatter}

\section{Introduction and notations}\label{sec1}

\subsection{Introduction}\label{sec1.1}

The $\zeta$-functions are one of more important special functions in
modern analysis and theory of functions. The most known and frequently
encountered $\zeta$-functions are Riemann and Hurwitz
$\zeta$-functions. They are classically introduced as the following series
\begin{eqnarray*}
\zeta(s)=\sum_{n=1}^\infty \frac{1}{ n^{s}}
, \qquad \zeta(s,v)=\sum_{n=0}^\infty
\frac{1}{ (n+v)^{s}} , \quad v\neq0, -1, -2, \ldots
\end{eqnarray*}
convergent for $\operatorname{Re}{s}>1$, and may be extended to other
domains of $s$ by the principle of analytic continuation. It is well
known that $\zeta(s)$ and $\zeta(s,v)$ are meromorphic on the entire
complex $s$-plane and that their only pole is a simple pole at $ s=1 $
with residue 1. They can be, therefore, expanded in the Laurent series
in a neighborhood of $s=1$ in the following way
%
\begin{eqnarray}
\label{dhd73vj6s1} \zeta(s) = \frac{1}{ s-1 } + \sum_{n=0}^\infty
\frac{(-1)^n
(s-1)^n}{n!} \gamma_n , \quad s\neq1,
\end{eqnarray}
and
%
\begin{eqnarray}
\label{dhd73vj6s2} \zeta(s,v) = \frac{1}{ s-1 } + \sum
_{n=0}^\infty\frac{(-1)^n
(s-1)^n}{n!} \gamma_n(v)
, \quad s\neq1,
\end{eqnarray}
respectively.
Coefficients $ \gamma_n $ appearing in the regular part of expansion
{\eqref{dhd73vj6s1}} are called \emph{Stieltjes constants} or \emph{generalized Euler's constants}, while
those appearing in the regular part of {\eqref{dhd73vj6s2}}, $ \gamma
_n(v) $, are called \emph{generalized Stieltjes constants}.
It is obvious that $\gamma_n(1)=\gamma_n$ since $\zeta(s,1)=\zeta(s)$.

The study of these coefficients is an interesting subject
and may be traced back to the works of Thomas Stieltjes and Charles Hermite
\cite[vol.~I, letter 71 and following]{stieltjes_01}.
In 1885, first Stieltjes and then Hermite, proved that
%
\begin{eqnarray}
\label{k98y9g87fcfcf} \gamma_n = \lim_{m\to\infty} \left\{
\sum_{k=1}^m \frac{\ln^n k}{k} -
\frac{\ln^{n+1} m}{n+1} \right\} , \quad n=0, 1, 2,\ldots
\end{eqnarray}
Later, this formula was also obtained or simply stated in works of
Johan Jensen \cite{jensen_02,jensen_03}, J{\o}rgen Gram \cite{gram_01},
Godfrey Hardy \cite{hardy_03}, Srinivasa Ramanujan \cite{ramanujan_01}
and many others. From~{\eqref{k98y9g87fcfcf}}, it is visible
that $\gamma_0$ is Euler's constant $\gamma$. However, the study of
other Stieltjes constants~revealed to be more difficult and, at the
same time, interesting. In 1895 Franel \cite{franel_01}, by using
contour integration techniques, showed that\footnote{There was a
priority dispute between Jensen, Kluyver and Franel related to this
formula \cite{franel_01,jensen_03}. In fact, it can be
straightforwardly deduced from the first integral formula for the
$\zeta$-function {\eqref{kjd02jddnsa}} which was first
obtained by Jensen in 1893 \cite{jensen_04}. By the way, we corrected
the original Franel's formula which was not valid for $n=0$ [this
correction comes from {\eqref{lkjc982dhndgs}} and
{\eqref{lkjc982dhndgs2}} here after].}
%
\begin{eqnarray}
\label{kljc3094jcmfd} \gamma_n = \frac{1}{2}\delta_{n,0}+
\frac{1}{i} \! \int\limits_0^\infty \frac{dx}{e^{2\pi x}-1}
\left\{ \frac{\ln^n(1-ix)}{1-ix} - \frac{\ln^n(1+ix)}{1+ix} \right\} ,
\quad n=0, 1, 2,\ldots
\end{eqnarray}
Ninety years later this integral formula was discovered independently
by Ainsworth and Howell who also provided a very detailed proof of it
\cite{ainsworth_01}. Following Franel's line of reasoning, one can
also obtain these formul\ae\footnote{The proof is analogous to that
given for the formula {\eqref{lkjc982dhndgs}} here after,
except that the Hermite representation should be replaced by the third
and second Jensen's formul\ae~for $\zeta(s)$
{\eqref{kjd02jddnsa}} respectively.}
%
\be\label{i2jd298jx2onxjcnf}
\gamma_n \,=\,-\frac{\pi}{2(n+1)}\! \int\limits_{-\infty}^{+\infty} 
\frac{\ln^{n+1}\!\!\big(\frac{1}{2}\pm ix\big)}{\ch^2\!\pi x}\, dx 
\qquad\qquad\qquad\qquad\qquad\qquad n=0, 1, 2,\ldots
\ee
and
\be\label{h8ye2bd43h}
\begin{array}{l}
\displaystyle 
\gamma_1 =-\left[\gamma -\frac{\ln2}{2}\right]\ln2+\,i\!\int\limits_0^\infty \! \frac{dx}{e^{\pi x}+1} \left\{
\frac{\ln(1-ix)}{1-ix} - \frac{\ln(1+ix)}{1+ix} 
\right\}\, \\[6mm]
\displaystyle
\gamma_2 =-\left[2\gamma_1 +\gamma\ln2-\frac{\ln^2\!2}{3}\right]\ln2
+\,i\!\int\limits_0^\infty \! \frac{dx}{e^{\pi x}+1} \left\{
\frac{\ln^2(1-ix)}{1-ix} - \frac{\ln^2(1+ix)}{1+ix} \right\}\,   \\[6mm]
\displaystyle
\gamma_3 =-\left[3\gamma_2 +3\gamma_1\ln2+\gamma\ln^2\!2-\frac{\ln^3\!2}{4}\right]\ln2
+\,i\!\int\limits_0^\infty \! \frac{dx}{e^{\pi x}+1} \left\{
\frac{\ln^3(1-ix)}{1-ix} - \frac{\ln^3(1+ix)}{1+ix} \right\}\, \\[3mm]
\qquad\qquad\qquad\qquad\qquad\qquad\qquad\qquad\ldots\ldots\ldots
\end{array}
\ee
first of which is particularly simple.\footnote{Despite the surprising
simplicity of these formul\ae, we have not found them in the
literature devoted to Stieltjes constants. In contrast, formula
{\eqref{i2jd298jx2onxjcnf}} seems to be known; at least its
variant for the generalized Stieltjes constant~may be found in
\cite{coffey_04}.} Other important results concerning the Stieltjes
constants~lie in the field of rational expressions of natural numbers,
as well as in the closely related field of integer parts of functions.
In 1790 Lorenzo Mascheroni \cite[p.~23]{mascheroni_01}, by using some previous findings of Gregorio Fontana, 
showed that\footnote{The series itself was given by Fontana, who, however, failed to find a constant to which 
it converges (he only proved that it should be lesser than 1). Mascheroni identified this \emph{Fontana's constant} 
and showed that it equals Euler's constant  \cite[pp.~21--23]{mascheroni_01}. Taking into account that both 
Fontana and Mascheroni did practically the equal work, we propose to call {\eqref{nc3bknkvb8}} \emph{Fontana--Mascheroni's series} for Euler's constant.}
\begin{eqnarray}
\label{nc3bknkvb8}
\gamma= \sum_{k=1}^\infty
\frac{|a_k|}{k} , \quad\mbox{where }\;  \frac{z}{\ln(1+z)} = 1+\sum
_{k=1}^\infty a_k z^k ,\ |z|<1
\end{eqnarray}
i.e.~$a_k$ are coefficients in the Maclaurin expansion of $z/\ln(1+z)$ and are usually referred as to \emph{(reciprocal) logarithmic numbers} or \emph{Gregory's coefficients} 
(in particular $a_1=\frac{1}{2}$, \mbox{$a_2=-\frac{1}{12}$}, $a_3=\frac{1}{24}$, $a_4=-\frac{19}{720}$, $a_5=\frac{3}{160}$, $a_6=-\frac{863}{60\,480}$, \ldots).\footnote{These coefficients have a venerable history and were named after James Gregory 
who gave first six of them in November 1670 in a letter to John Collins
\cite[vol.~1, p.~46]{newton_01} (although in the fifth coefficient there is an error or misprint: $\frac{3}{164}$ should be replaced by $\frac{3}{160}$).
Coefficients $a_k$ are also closely related
to the \emph{Cauchy numbers} of the first kind $C_{1, k}$, to the
\emph{generalized Bernoulli numbers}, to the \emph{Stirling
polynomials} and to the \emph{signed Stirling numbers} of the first
kind $ S_1(k,l) $. In particular, $ a_k=\frac{C_{1,
k}}{k!}=\frac{1}{k!}\sum\frac {S_1(k,l)}{l+1} $, where the
summation extends over $ l=[1,k] $, see e.g.~\cite{davis_02},
\cite[pp.~293--294, \no13]{comtet_01}, \cite[vol.~III,
p.~257]{bateman_01}, \cite{skramer_01}, \cite{candelpergher_01}. } Fontana--Mascheroni's series {\eqref{nc3bknkvb8}} seems to be the first known series representation for Euler's
constant containing rational coefficients only and was 
subsequently rediscovered several times, in particular by Kluyver in 1924 \cite{kluyver_02}, by Kenter in
1999 \cite{kenter_01} and by Kowalenko in 2008 \cite{kowalenko_01} (this list is far from exhaustive, see e.g.~\cite{skramer_01}).
In 1897 Niels Nielsen \cite[Eq.~(6)]{nielsen_02} showed that
%
\begin{eqnarray}
\label{d34rf13d436f} \gamma=1-\sum_{k=1}^\infty\sum
_{l=2^{k-1}}^{2^{k}-1} \frac{k}{(2l+1)(2l+2)}
\end{eqnarray}
This formula was also the subject of several rediscoveries, e.g.~by Addison
in 1967 \cite{addison_01} and
by Gerst in 1969 \cite{gerst_01}.\footnote{The actual Addison's formula
\cite{addison_01} is slightly different, but it straightforwardly
reduces to {\eqref{d34rf13d436f}} by partial fraction decomposition.
In \cite{addison_01}, we also find a misprint: the upper bound in the
second sum on p.~823 should be the same as in {\eqref{d34rf13d436f}}.
As regards Gerst's formula \cite{gerst_01}, it is exactly the same
as {\eqref{d34rf13d436f}}.}
In 1906 Ernst Jacobsthal \cite[Eq.~(9)]{jacobsthal_01} and in 1910 Giovanni Vacca \cite{vacca_01}, apparently independently, derived a closely related series
%
\begin{eqnarray}
\label{lk2jd029jde} \gamma= \sum_{k=2}^\infty
\frac{(-1)^k}{k} \lfloor\log_2{k}\rfloor
\end{eqnarray}
which was also rediscovered in numerous occasions, in particular by
H.F.~Sandham in 1949 \cite{sandham_01},
by D.F.~Barrow, M.S.~Klamkin and N.~Miller in 1951 \cite{barrow_01}
or by Gerst in 1969 \cite{gerst_01}.\footnote{Series {\eqref{lk2jd029jde}}, thanks to the error of Glaisher, Hardy and Kluyver,  was long-time attributed to Giovanni Vacca and is widely known as \emph{Vacca's series}, see e.g.~\cite{glaisher_01,hardy_03,kluyver_02}. It was only in 1993 that Stefan Kr\"amer found that this series was first obtained by Jacobsthal in 1906. Besides, Kr\"amer also showed that Nielsen's series {\eqref{d34rf13d436f}} and Jacobsthal--Vacca's series {\eqref{lk2jd029jde}} are closely related and can be derived one from another via a simple geometrical progression $\frac{1}{2}=\frac{1}{4}+\frac{1}{8}+\frac{1}{16}+\ldots$ \cite{skramer_01}.} In 1910 James Glaisher \cite{glaisher_01} proposed yet another proof of the same result and derived a number of other series with rational terms for $\gamma$. In 1912 Hardy \cite{hardy_03} extended {\eqref{lk2jd029jde}} to the first Stieltjes constant
%
\begin{eqnarray}
\label{jm9c204dj} \gamma_1 = \frac{\ln2}{2}\sum
_{k=2}^\infty\frac{(-1)^k}{k} \lfloor
\log_2{k}\rfloor\cdot \left(2\log_2{k} - \lfloor
\log_2{2k}\rfloor \right)
\end{eqnarray}
However, this expression is not a full generalization of
{\eqref{lk2jd029jde}}
since it also contains irrational coefficients.
In 1924 Jan Kluyver \cite{kluyver_02} generalized Jacobsthal--Vacca's series {\eqref{lk2jd029jde}} in the another direction and showed that
%
\begin{eqnarray}
\label{f413f14cx45} \gamma= \sum_{k=m}^\infty
\frac{ \beta_k }{k} \lfloor\log _m{k}\rfloor ,\qquad
\beta_k = 
\begin{cases}
m-1 , & k = \mbox{multiple of } m \vspace{3pt}\cr
-1 , & k \neq\mbox{multiple of } m
\end{cases}
\end{eqnarray}
where $m$ is an arbitrary chosen positive integer.\footnote{For
example, if $ m=2 $ then $ \beta_k=(-1)^k $.}
In 1924--1927 Kluyver \cite{kluyver_03,kluyver_01} also tried to obtain series with rational
coefficients for higher Stieltjes constants, but these attempts were not
successful. Currently, apart from $\gamma_0$, no closed-form
expressions are known for $\gamma_n$. However, there are works devoted
to their estimations and to the asymptotic series representations for them
\cite{berndt_02,israilov_01,lammel_01,lavrik_01_eng,stankus_01_eng,zhang_01,matsuoka_01,matsuoka_02,knessl_01,adell_01,eddin_01,eddin_02}. Besides, there are also works devoted to the
behavior of their sign \cite{briggs_01,mitrovic_01}. In particular,
Briggs in 1955 \cite{briggs_01} demonstrated that there are infinitely
many changes of sign for them. Finally, aspects related to the
computation of Stieltjes constants were considered in \cite{ainsworth_01,gram_01,kreminski_01,todd_01}. 

As regards generalized Stieltjes constants, they are much less studied
than the usual Stieltjes constants.
In 1972 Berndt,
by employing the Euler--Maclaurin summation formula and by proceeding
analogously to Lammel \cite{lammel_01},
showed that $\gamma_n(v)$ can be given by an asymptotic representation
of the same kind as {\eqref{k98y9g87fcfcf}}
%
\begin{eqnarray}
\label{kjhd928hdnd} \gamma_n(v) = \lim_{m\to\infty}
\left\{ \sum_{k=0}^m
\frac{\ln^n (k+v)}{k+v} - \frac{\ln^{n+1} (m+v)}{n+1} \right\} , \quad \begin{array}[c]{l}
n=0, 1, 2,\ldots \\[6pt]
v\neq0, -1, -2, \ldots
\end{array}
\end{eqnarray}
see \cite{berndt_02}.\footnote{We, however, note that Wilton, by using Vall\'ee--Poussin's
expansion of the Hurwitz $\zeta$-function,
provided a similar formula already in 1927 \cite{wilton_01}.}
Similarly to Franel's method of the derivation of {\eqref
{kljc3094jcmfd}}, one may also derive
the following integral representation for the $n$th generalized
Stieltjes constant
%
\be\label{lkjc982dhndgs}
\gamma_n(v) \,=\,\left[\frac{1}{2v}-\frac{\ln{v}}{n+1} \right]\ln^n\!{v}
-i\!\int\limits_0^\infty \! \frac{dx}{e^{2\pi x}-1} \left\{
\frac{\ln^n(v-ix)}{v-ix} - \frac{\ln^n(v+ix)}{v+ix} 
\right\}
\ee	
$n=0, 1, 2,\ldots{}$,
$\operatorname{Re}{v}>0$.\footnote{This formula follows
straightforwardly from the well-known Hermite representation for
$\zeta(s,v)$, see e.g.~\cite[p.~66]{hermite_01}, \cite[p.~106]{lindelof_01},
\cite[vol.~I, p.~26, Eq.~1.10(7)]{bateman_01}. First,
recall that
$2(v^2+x^2)^{-s/2}\sin[s\operatorname
{arctg}(x/v)]=-i[(v-ix)^{-s}-(v+ix)^{-s}]$, and then,
expand $\frac{1}{2}v^{-s} + (s-1)^{-1} v^{1-s}$
into the Laurent series about $s=1$. Performing the term-by-term
comparison of the
derived expansion with the Laurent series {\eqref{dhd73vj6s2}} yields
{\eqref{lkjc982dhndgs}}.}
This formula was rediscovered several times, for example, by Mark
Coffey in 2009 \cite{coffey_04,connon_01}.
From both latter formul\ae, it follows that $\gamma_0(v)=-\Psi
(v)$. Consider, for instance, {\eqref{lkjc982dhndgs}} and put $n=0$.
Then, the latter equation takes the form
%
\be\label{lkjc982dhndgs2}
\gamma_0(v) \,=\,\frac{1}{2v} - \ln{v} +2 \!
\underbrace{\int\limits_0^{\,\infty} \!\frac{\,x\,dx\,}{(e^{2\pi x}-1)(v^2+x^2)}}_{
-\frac{1}{4v} + \frac{1}{\,2\,}\ln{v} - \frac{1}{\,2\,}\Psi(v)}  = -\Psi(v)
\ee	
where the last integral was first calculated by Legendre.\footnote{And
not by Binet as stated in
\cite[vol.~I, p.~18, Eq.~1.7.2(27)]{bateman_01}, see \cite[tome II,
p.~190]{legendre}
and \cite[p.~83, \no40, Eq.~(55)]{iaroslav_06}.}
The demonstration of the same result from formula {\eqref{kjhd928hdnd}}
may be found, for example, in \cite{zhang_01}.
For rational $v$, the $0$th Stieltjes constant~may be, therefore,
expressed by means of Euler's constant $\gamma$ and a finite combination
of elementary functions [thanks to the Gauss' Digamma theorem (\ref{dh38239djws}a,b)].
However, things are much more complicated for higher generalized
Stieltjes constants;
currently, no closed-form expressions are known for them
and little is known
as to their arithmetical properties.
Basic properties, such as the multiplication theorem
\begin{eqnarray*}
\sum_{l=0}^{n-1}
\gamma_p \left( v+\frac{l}{ n } \right) =
(-1)^p n \left[\frac{\ln n}{ p+1 } - \Psi(nv) \right]
\ln^p n + n\sum_{r=0}^{p-1}(-1)^r
C_p^r \gamma_{p-r}(nv) \cdot\ln^r
{n} ,
\end{eqnarray*}
$n=2, 3, 4,\ldots{}$, where $C_p^r$ denotes the binomial coefficient $
C_p^r = \frac{p!}{ r! (p-r)! } $,
and the recurrent relationship
%
\begin{eqnarray}
\label{qkw9123jssn} \gamma_p(v+1) = \gamma_p(v) -
\frac{ \ln^p v }{v} , \quad %
\begin{array}[c]{l}
p=1, 2, 3,\ldots\\[6pt]
v\neq0, -1, -2,\ldots
\end{array} %
\end{eqnarray}
may be both straightforwardly derived from those for the Hurwitz
$\zeta$-function, see
e.g.~\cite[pp.~101--102]{iaroslav_06}.\footnote{As regards the
multiplication theorem, see e.g.~\cite[Eq.~(6.6)]{connon_02} or
\cite[p.~101]{iaroslav_06}. We can also find its particular case for $
v=1/n $ in \cite[p.~1830, Eq.~(3.28)]{coffey_01}.} In attempt to obtain
other properties, several summation relations involving single and
double infinite series were quite recently obtained in
\cite{coffey_03,coffey_02}. Also, many important aspects regarding the
Stieltjes constants~were considered by Donal Connon
\cite{connon_01,connon_02,connon_03}.

Let now focus our attention on the first generalized Stieltjes constant.
The most strong and pertinent results in the field of its closed-form
evaluation
is the formula for the difference between the first generalized
Stieltjes constant~at rational argument and its reflected version
%
\be\label{hgxcw3b2}
\gamma_1 \biggl(\frac{m}{n} \vphantom{\frac{1}{2}} \biggr)
- \gamma_1 \biggl(1-\frac{m}{n} \vphantom{\frac{1}{2}}
\biggr) =\, 2\pi\!\sum_{l=1}^{n-1}
\sin\frac{2\pi m l}{n} \cdot\ln\Gamma \biggl(\frac{l}{n} \biggr) -\pi(
\gamma+\ln2\pi n)\ctg\frac{m\pi}{n}
\ee
In the literature devoted to Stieltjes constants this result is usually
attributed to Almkvist and Meurman who obtained it by deriving the
functional equation for $\zeta(s,v)$, Eq.~\eqref{kljc928c2ndbn}, with respect to $s$ at rational $v$, see
e.g.~\cite{adamchik_01}, \cite[p.~261, \S12.9]{apostol_01},
\cite[Eq.~(6)]{miller_01}. However, it was comparatively recently that
we discovered that this formula, albeit in a slightly different form,
was obtained by Carl Malmsten already in 1846. On pp.~20 and 38
\cite{malmsten_01}, we, \emph{inter alia}, find the following
expression
%
\be\label{odi239dn}
\begin{array}{l}
\displaystyle 
\sum_{l=0}^\infty\left\{\frac{\ln\!\big[(2l+1)n-m\big]}{(2l+1)n-m} - \frac{\ln\!\big[(2l+1)n+m\big]}{(2l+1)n+m}\right\}=\\[8mm]
\displaystyle \qquad
 \,=\,
\begin{cases}
\displaystyle -\frac{\pi(\gamma+\ln2\pi)}{2n}\tg\frac{\,\pi m\,}{2n} - \frac{\pi}{n}\cdot\!\sum_{l=1}^{n-1} (-1)^{l-1}  \sin\frac{\,\pi ml\,}{n}\cdot 
\ln\left\{\!\frac{\Gamma\!\left(\!\frac{n+l}{2n}\!\right) }{\Gamma\!\left(\!\frac{l}{2n}\!\right)}\right\} \\[6mm]
\qquad\qquad \qquad \qquad \qquad \qquad\qquad \qquad \qquad \qquad  \text{ if $\;m+n\;$ is odd,} \\[3mm]
\displaystyle -\frac{\pi(\gamma+\ln\pi)}{2n}\tg\frac{\,\pi m\,}{2n}  - \frac{\pi}{n}\cdot\!\!\!\!\sum_{l=1}^{\lfloor\!\frac{1}{2}(n-1)\!\rfloor} \!\!\!\! (-1)^{l-1} \sin\frac{\,\pi ml\,}{n}\cdot 
\ln\left\{\!\frac{\Gamma\!\left(\!\frac{n-l}{n}\!\right) }{\Gamma\!\left(\!\frac{l}{n}\!\right)}\right\} \\[6mm]
\qquad\qquad \qquad \qquad \qquad \qquad\qquad \qquad \qquad \qquad \text{ if $\;m+n\;$ is even,} \\
\end{cases}
\end{array}
\ee
where $m$ and $n$ are integers such that $m<n$.\footnote{Unfortunately,
this Malmsten's work contains a huge quantity of misprints in
formul\ae. We already corrected many of them
in our previous work \cite[Sections 2.1 \& 2.3]{iaroslav_06}. As regards the
above-referenced Malmsten's original
equation (55), case $ m+n $ even, note that $\Gamma (\frac
{n-i}{2n} )$ should be replaced
by $\Gamma (\frac{n-i}{n} )$. Formula (56) also has an
error: $\Gamma (\frac{n+i}{n} )$ should be replaced
by $\Gamma (\frac{n-i}{n} )$.} It is visible that the left
part of this equality contains
the difference of two first-order derivatives of $\zeta(s,v)$ at
$s\to1$ and $v=\frac{1}{2}\pm\frac{m}{2n}$.
Putting $2m-n$ instead of $m$ and using the Laurent series expansion
{\eqref{dhd73vj6s2}} yields, after some algebra,
formula {\eqref{hgxcw3b2}}. A somewhat
different way to get {\eqref{hgxcw3b2}} is to directly apply the
Mittag--Leffler theorem to
one of Malmsten's integrals at rational points; we developed such a method
in our preceding study \cite[pp.~97--98, \no63 and pp.~106--107, \no67]{iaroslav_06}.

Recently, Coffey \cite{coffey_01} derived several formul\ae~for the
linear combination of the first generalized Stieltjes constants~at some
rational arguments.
From these expressions, one may conjecture that in some cases (author
gave only two examples of such cases
\cite[p.~1821, Eqs.~(3.33)--(3.34)]{coffey_01}),
not only the $\Gamma$-function, but
also the second-order derivative of the Hurwitz $\zeta$-function
could be related, in some way, to the first
generalized Stieltjes constant. However, these preliminary findings do
not permit to precisely identify their roles
in the general problem of the closed-form evaluation of the first
Stieltjes constant~at any rational argument (the problem which we come
to solve here).

Very recently, it has been
conjectured in \cite[p.~103]{iaroslav_06}
that similarly to the Digamma theorem for $\gamma_0(v)$, the first
generalized Stieltjes constant~$\gamma_1(v)$ at rational $v$ may be expressed
by means of the Euler's constant~$\gamma$,
the first Stieltjes constant $ \gamma_1 $, the $\Gamma$-function
and some ``relatively simple'' function.
For seven rational values of $v$ in the range $(0,1)$, namely for
$\frac{1}{6}$, $\frac{1}{4}$,
$\frac{1}{3}$, $\frac{1}{2}$, $\frac{2}{3}$, $\frac{3}{4}$ and
$\frac{5}{6}$,
we showed in \cite[pp.~98--101, \no64]{iaroslav_06} that this ``relatively
simple'' function is
elementary.\footnote{Further to remarks we received after the
publication of \cite{iaroslav_06},
we note that similar closed-form expressions for $\gamma(1/4)$,
$\gamma(3/4)$ and $\gamma(1/3)$ were also obtained in \cite[pp.~17--18]{connon_03}.}
In this manuscript, we extend these preceding researches by providing a
theorem which allows to evaluate the
first generalized Stieltjes constant~at any rational argument
in a closed-form by precisely identifying this ``relatively simple''
function. The latter consists of elementary functions
containing the Euler's constant $\gamma$ and
of the reflected sum of two second-order derivatives of the Hurwitz
$\zeta$-function at zero
$\zeta''(0,p)+\zeta''(0,1-p)$, parameter $p$ being rational
in the range $(0,1)$. A close study
of this reflected sum reveals that it has several important integral
and series representations, one of which is
quite similar to an integral representation for the logarithm of the
$\Gamma$-function at rational argument
(see Section \ref{kjd12093ddmne3dd} and \ref{894yf3hedbe}).
Moreover, the derived theorem represents also the finite Fourier series
for the first generalized Stieltjes constant,
so that classic Fourier analysis tools may be used at their full
strength. With the help of the latter, we
derive several summation formul\ae~including summation with
trigonometric functions, summation with square, summation
with the Digamma function and summation giving the first-order moment (see
Section \ref{ij2093dm32ddf}).
Obviously, the same method can be applied to other discrete functions
allowing similar representations. In particular,
its application to a variant of Gauss' Digamma theorem
yields several beautiful
summation formul\ae~for the Digamma function which are derived in
 \ref{lo109sj1s2m2w}.
We also derive, in passing
[in \ref{894yf3hedbe}, Eq.~{\eqref{kjh298hdnd2}}],
an interesting integral representation
for the logarithm of the $\Gamma$-function at rational argument.
Finally, in Section~\ref{iou2ch20983j3e},\break we discuss extensions of the
derived theorem to the higher
Stieltjes constants~and provide closed-form expressions for the second
generalized Stieltjes constant~at rational arguments.

\subsection{Notations}\label{notations}

Throughout the manuscript, the following abbreviated notations are used:
$\gamma=0.5772156649\ldots$ for Euler's constant, $\gamma_n$
for the $n$th Stieltjes constant, $\gamma_n(p)$ for the $n$th
generalized Stieltjes constant at point $p$, $\lfloor x\rfloor$ for the
integer part of $x$, $\operatorname{tg}z$ for the tangent of $z$,
$\operatorname{ctg}z$ for the cotangent of $z$, $\operatorname{ch}z$
for the hyperbolic cosine of $z$, $\operatorname{sh}z$ for the
hyperbolic sine of $z$, ${\operatorname{th}}z$ for the hyperbolic
tangent of $z$.\footnote{Most of these notations come from Latin,
e.g.~``$\operatorname{ch}$'' stands for \emph{cosinus hyperbolicus},
``$\operatorname{sh}$'' stands for \emph{sinus hyperbolicus}, etc.} In
order to avoid any confusion between compositional inverse and
multiplicative inverse, inverse trigonometric and hyperbolic functions
are denoted as $\arccos, \arcsin, \operatorname{arctg}, \ldots$
and not as $\cos^{-1}, \sin^{-1}, \operatorname{tg}^{-1}, \ldots{}$.
Writings $\Gamma(z)$, $\Psi(z)$, $\zeta(s)$ and $\zeta(s,v)$
denote respectively the $\Gamma$-function, the $\Psi$-function
(or Digamma function), the Riemann $\zeta$-function and the Hurwitz
$\zeta$-function. When referring to the derivatives of the Hurwitz
$\zeta$-function, we always refer to the derivative with respect to
its first argument $s$ (unless otherwise specified).
$\operatorname{Re}{z}$ and $\operatorname{Im}{z}$ denote, respectively,
real and imaginary parts of $z$. Natural numbers are defined in a
traditional way as a set of positive integers, which is denoted by
$\mathbbm{N}$. Kronecker symbol of arguments $l$ and $k$ is denoted by
$\delta_{l,k}$. Letter $i$ is never used as index and is $\sqrt{-1}$.
The writing $\operatorname{res}_{z=a} f(z)$ stands for the residue of
the function $f(z)$ at the point $z=a$. By Malmsten's integral we mean
any integral of the form
\begin{eqnarray*}
\int\limits_{0}^\infty \frac{R(\operatorname{sh}{px}, \operatorname
{ch}{px}) \cdot\ln{x}}{R(\operatorname{sh}{x}, \operatorname
{ch}{x})} \,dx
\end{eqnarray*}
where $R$ denotes a rational function
and the parameter $p$ is such that the convergence is guaranteed.
Other notations are standard.

\section{Evaluation of the first generalized Stieltjes constant~at
rational argument}\label{sec2}
\subsection{Generalized Stieltjes constants~and their relationships to
Malmsten's integrals}\label{sec2.1}

Formula {\eqref{hgxcw3b2}} provides a closed-form
expression for the difference of two first Stieltjes constants~at
rational arguments. It should be therefore interesting to investigate
if there could be some expressions containing other combinations of
Stieltjes constants. In our previous work
\cite[pp.~97--107]{iaroslav_06}, we already demonstrated that some
Malmsten's integrals are connected with the first generalized Stieltjes
constants. This connection was quite fruitful and permitted not only to
prove by another method formula
{\eqref{hgxcw3b2}}, but also to evaluate the first generalized
Stieltjes constant $ \gamma_1(p) $ at $ p=\frac{1}{2}, \frac{1}{3},
\frac{1}{4}, \frac{1}{6}, \frac {2}{3}, \frac{3}{4}, \frac{5}{6} $ by
means of elementary functions, Euler's constant $\gamma$, the first
Stieltjes constant $\gamma_1$ and the $\Gamma$-function, see for
more details \cite[pp.~98--101, \no64]{iaroslav_06}. Taking into account
that aforementioned manuscript was quite long, many results and
theorems were given as exercises with hints and without rigorous
proofs. Below, we provide several useful proofs and unpublished results
(given as lemmas and corollaries) showing that Malmsten's integrals of
the first and second orders may be expressed by means of the first
generalized Stieltjes constants. This connection between Malmsten's
integrals and Stieltjes constants~is crucial and plays the central role
in the proof of the main theorem of this manuscript.

\begin{lemma}\label{bnhgcewucv}
For any $|\operatorname{Re}{p}|<1$ and $\operatorname{Re}{a}>-1$,
%
\begin{eqnarray}
\label{lkij289dnxnwx} \int\limits_0^\infty \frac{ x^{a-1} (\operatorname{ch}{px}-1) }{
\operatorname{sh}{x} } \,dx
= \frac{\Gamma(a)}{2^{a}} \left\{\zeta \left(a,\frac{1}{2}-
\frac{p}{
2 } \right) + \zeta \left(a,\frac{1}{2}+
\frac{p}{ 2 } \right) -2 \left(2^{a}-1 \right)
\zeta(a) \right\}
\end{eqnarray}
\end{lemma}

\begin{pf}
From elementary analysis it is well-known that $\operatorname{sh}^{-1}
{x}$, for $\operatorname{Re}{x}>0$,
may be represented by the following geometric series
\begin{eqnarray*}
\frac{1}{ \operatorname{sh}{x} } = 2 \sum_{n=0}^\infty
e^{-(2n+1)x} , \quad\operatorname{Re} {x}>0 .
\end{eqnarray*}
This series, being uniformly convergent, can be integrated
term-by-term. Hence
\be\notag
\begin{array}{ll}
\displaystyle
\int\limits_0^\infty \!\!\frac{\,x^{a-1} (\ch{px}-1)\,}{\,\sh{x}\,}\,dx\,&=\displaystyle
\sum\limits_{n=0}^\infty \int\limits_0^\infty \!\! x^{a-1} \left\{e^{-(2n+1-p)x} + e^{-(2n+1+p)x}  - 2e^{-(2n+1)x}  \right\} dx  \\[5mm]
\displaystyle &=\displaystyle
\Gamma(a) \! \sum\limits_{n=0}^\infty \!\left\{ \frac{1}{\,(2n+1-p)^{a}\,} 
+  \frac{1}{\,(2n+1+p)^{a}\,} - \frac{2}{\,(2n+1)^{a}\,}  \right\} \\[6mm]
\displaystyle &=\displaystyle 
\frac{\Gamma(a)}{2^{a}}\left\{\zeta\!\left(\!a,\frac{1}{2}-\frac{p}{\,2\,}\right) 
+ \zeta\!\left(\!a,\frac{1}{2}+\frac{p}{\,2\,}\right) 
-2 \zeta\!\left(\!a,\frac{1}{2}\right)  \! \right\}, 
\end{array}
\ee
where the integral on the left converges if $|\operatorname{Re}{p}|<1$
and $\operatorname{Re}{a}>-1$.
In order to obtain {\eqref{lkij289dnxnwx}}, it suffices to notice that
$ \zeta (a,\frac{1}{2} )=  (2^a-1 )\zeta
(a) $.
\qed
\end{pf}

\begin{corollary}\label{k9wj2qpofncns}
For any $p$ lying in the strip $|\operatorname{Re}{p}|<1$, we always have
%
\begin{eqnarray}
\label{k290djd2s} \int\limits_0^{\infty} \frac{(\operatorname{ch}{px}-1) \ln
{x}}{\operatorname{sh}{x}}
\,dx &=& (\gamma+\ln2 )\cdot \left\{\Psi \left(
\frac{1}{2}+\frac{p}{2} \right) + \ln2 - \frac{\pi}{2}
\operatorname{tg}\frac{\pi p}{2} \right\}
\nonumber
\\
&&{}+ \gamma^2 + \gamma_1 - \frac{1}{2}
\gamma_1 \!\left(\frac
{1}{2}+\frac{p}{2}
\right) - \frac{1}{2}\gamma_1 \!\left(
\frac{1}{2}-\frac{p}{2} \right).
\end{eqnarray}
This result is straightforwardly obtained from {Lemma~\ref{bnhgcewucv}}
by differentiating {\eqref{lkij289dnxnwx}} with respect to $a$,
and then by making $a\to1$. In order to evaluate the limit in the
right-hand side, we make
use of Laurent series {\eqref{dhd73vj6s1}} and {\eqref{dhd73vj6s2}}.
\end{corollary}

Another Malmsten's integral of the first order which also contains
Stieltjes constants~appear in the next lemma.

\begin{lemma}\label{lem2}
For any $|\operatorname{Re}{p}|<1$ and $\operatorname{Re}{a}>-1$,
\begin{eqnarray*}
\int\limits_0^\infty \frac{ x^{a-1} \operatorname{sh}{px} }{ \operatorname
{ch}{x} } \,dx &=&
\frac{\Gamma(a)}{2^{a}} \left\{\zeta \left( a,\frac{1}{2}+
\frac
{p}{ 2 } \right) - \zeta \left( a,\frac{1}{2}-
\frac{p}{ 2 } \right) \right.
\\
&&{}\left. -2^{1-a} \zeta \left( a,\frac{1}{4}+
\frac{p}{ 4 } \right) +2^{1-a}\zeta \left( a,
\frac{1}{4}-\frac{p}{ 4 } \right) \right\}
\end{eqnarray*}
\end{lemma}

\begin{pf}
Analogous to that of {Lemma \ref{bnhgcewucv}}.\qed
\end{pf}

\begin{corollary}\label{kjhxc82hd}
For any $|\operatorname{Re}{p}|<1$,
\begin{eqnarray*}
\int\limits_0^\infty \frac{ \operatorname{sh}{px}\cdot\ln{x} }{
\operatorname{ch}{x} } \,dx &= &
\frac{1}{2} \left\{ \pi(\gamma+\ln 2)\operatorname{tg}
\frac{ \pi p }{2} - (\gamma+2\ln2) \left[\Psi \left(
\frac{1}{4}+\frac{p}{4} \right) - \Psi \left(
\frac{1}{4}-\frac{p}{4} \right) \right]\right.
\\
&&{}\;\;\left. + \gamma_1 \left( \frac{1}{2}-\frac{p}{2}
\right) - \gamma _1 \left( \frac{1}{2}+
\frac{p}{2} \right) - \gamma_1 \left(
\frac{1}{4}-\frac{p}{4} \right) + \gamma_1
\left( \frac{1}{4}+\frac{p}{4} \right) \right
\}.
\end{eqnarray*}
This result can be shown in the same way as that in {Corollary \ref{k9wj2qpofncns}}.
\end{corollary}

By the same line of argument, one may also prove that following
logarithmic integrals may be expressed in terms
of first generalized Stieltjes constants.
\be\label{hd2928h32iodnb}
\begin{array}{l}
\displaystyle
\int\limits_0^\infty \!\frac{\,\sh{px}\cdot \ln{x}\,}{\,\sh{x}\,}\,dx   \,=\,
-\frac{1}{2}\left\{\! \pi (\gamma+\ln2)\tg\frac{\,\pi p\,}{2} \, 
+ \,\gamma_1\!\left(\!\frac{1}{2}-\frac{p}{2}\! \right)  - \, \gamma_1\!\left(\!\frac{1}{2}+\frac{p}{2}\! \right) \!\right\} \\[7mm]
\displaystyle
\int\limits_0^\infty \!\!\frac{\,\ch{px}\cdot\ln{x}\,}{\,\ch{x}\,}\,dx\,=\,
\frac{1}{2}\left\{\gamma_1\!\left(\frac{1}{2}+\frac{p}{\,2\,}\right) + \gamma_1\!\left(\frac{1}{2}-\frac{p}{\,2\,}\right)
- \gamma_1\!\left(\frac{1}{4}+\frac{p}{\,4\,}\right) - \gamma_1\!\left(\frac{1}{4}-\frac{p}{\,4\,}\right)  \! \right\}		\\[6mm]
\displaystyle \phantom{mmm}
-\frac{1}{2}\ln^2\!2 
+\ln2\cdot \Psi\!\left(\!\frac{1}{\,2\,} + \,\frac{p}{\,2\,}\! \right) +\, \frac{\pi}{2}(\gamma+\ln2)\tg\frac{\pi p}{2} 
- \,\frac{\pi}{2}(\gamma+2\ln2) \ctg\!\left(\frac{\pi}{4}-\frac{\pi p}{4} \right) \\[4mm]
\displaystyle
\int\limits_0^\infty \!\frac{\,\sh^2\!{px}\cdot\ln{x}\,}{\,\sh^2\!{x}\,}\,dx\,=\frac{1}{2}\Big\{\! 
\ln\pi -\ln\sin\pi p +p\big[\gamma_1(p) - \gamma_1(1-p) \big] - \, (\gamma+\ln2) (1 -\pi p\ctg \pi p ) \! \Big\}  \\[7mm]
\displaystyle
\int\limits_0^\infty \!\frac{\,\ch{px}\cdot\ln{x}\,}{\,\ch^2\!{x}\,}\,dx\, = \displaystyle\frac{p}{2}
\left\{\gamma_1\biggl(\frac{p}{2} \vphantom{\frac{1}{2}} \biggr) - 
\gamma_1\biggl(1-\frac{p}{2} \vphantom{\frac{1}{2}} \biggr) -
\gamma_1\biggl(\frac{p}{4} \vphantom{\frac{1}{2}} \biggr)
+\gamma_1\biggl(1-\frac{p}{4} \vphantom{\frac{1}{2}} \biggr) \!\right\}+\ln\tg\frac{\pi p}{4}\\[6mm]
 \displaystyle\displaystyle \phantom{mmmmmmmmmmmmmmmmmmmm}
- \frac{\pi p}{2}\left\{(\gamma+2\ln2)\csc\frac{\pi p}{2} 
+ \ln2 \cdot\ctg\frac{\pi p}{2}\right\}
\end{array}
\ee
where parameter $p$ should be such that $|\operatorname{Re}{p}|<1$ in
the first three integrals and \mbox{$|\operatorname{Re}{p}|<2$} in the
fourth one. Interestingly, higher Malmsten's integrals seem to not
contain higher Stieltjes constants, but rather other $\zeta$-function
related constants. For instance, the evaluation of the third-order
Malmsten's integral by the same method yields:
\be\label{kx02jixdmnfdmc}
\begin{array}{l}
\displaystyle
\int\limits_0^\infty \!\frac{\,\sh^3\!{px}\cdot\ln{x}\,}{\,\sh^3\!{x}\,}\,dx\, = \displaystyle\frac{1}{4}
\left\{3\zeta'\biggl(\!-1,\frac{1}{2}+\frac{p}{2} \biggr) - 
3\zeta'\biggl(\!-1,\frac{1}{2}-\frac{p}{2} \biggr)  - \zeta'\biggl(\!-1,\frac{1}{2}+\frac{3p}{2} \biggr) \right.\\[6mm]
\displaystyle\qquad \qquad \qquad \qquad \quad
\left. + \zeta'\biggl(\!-1,\frac{1}{2}-\frac{3p}{2} \biggr)  \!\right\}+ \frac{3(1-p^2)}{16}
\left\{\gamma_1\biggl(\frac{1}{2}+\frac{p}{2} \biggr) - 
\gamma_1\biggl(\frac{1}{2}-\frac{p}{2} \biggr) \!\right\}
\\[6mm]
 \displaystyle \qquad\qquad\qquad\qquad\quad
 - \frac{1-9p^2}{16}\left\{\gamma_1\biggl(\frac{1}{2}+\frac{3p}{2} \biggr) - 
\gamma_1\biggl(\frac{1}{2}-\frac{3p}{2} \biggr) \!\right\} 
-\frac{3p}{4}\ln\!\big(2\cos\pi p - 1\big) 
\\[6mm]
 \displaystyle \qquad\qquad\qquad\qquad\quad
+ \frac{\pi(\gamma+\ln2)}{16}\left\{
3(p^2-1)\tg\frac{\pi p}{2} - (9p^2-1)\tg\frac{3\pi p}{2}\right\}
\end{array}
\ee
in the strip $|\operatorname{Re}{p}|<1$.
In contrast, the evaluation of Malmsten's integrals containing higher
powers of the logarithm
in the numerator of the integrand\footnote{We propose to call such
integrals \emph{generalized Malmsten's integrals}.}
leads precisely to higher Stieltjes constants.
In fact, differentiating twice {\eqref{lkij289dnxnwx}} with respect to $a$,
and then making $a\to1$, yields
\be\label{kjc294djcwx}
\begin{array}{l}
\displaystyle
\int\limits_0^{\infty}\!\frac{(\ch{px}-1)\, \ln^2\!{x}}{\sh{x}}\,dx\, =
-\gamma_2
+\frac{1}{2}\left\{\gamma_2\biggl(\frac{1}{2}+\frac{p}{2} \biggr) + 
\gamma_2\biggl(\frac{1}{2}-\frac{p}{2} \biggr) \!\right\}
-2 \gamma_1(\gamma-\ln2) \\[6mm]
\displaystyle
\qquad\qquad\qquad\quad
+(\gamma+\ln2)\cdot
\left\{\gamma_1\biggl(\frac{1}{2}+\frac{p}{2} \biggr) + 
\gamma_1\biggl(\frac{1}{2}-\frac{p}{2} \biggr) \!\right\} - \gamma^3 
-\frac{\gamma}{6}\big(\pi^2+6\ln^2\!2\big)  \\[6mm]
\displaystyle
\qquad\qquad\qquad\quad
-\biggl[(\gamma+\ln2)^2+\frac{\pi^2}{6}\biggr]\cdot \biggl\{\Psi \biggl(\frac{1}{2}+\frac{p}{2} \biggr) 
- \frac{\pi}{2}\tg\frac{\pi p}{2} \biggr\} 
-\frac{\ln2}{3}\big(\pi^2+2\ln^2\!2\big) 
\end{array}
\ee
where $|\operatorname{Re}{p}|<1$. More generally, the same integral
containing $\ln^n {x}$ instead of $\ln^2 {x}$ will lead to the $n$th
Stieltjes constants.

\medskip
\noindent\textbf{Nota Bene.}\label{hd91283dn9323hnd} As showed in
\cite[pp.~51--60, Sect.~4, \no3, 6, 11, 13]{iaroslav_06}, integrals
\eqref{hd2928h32iodnb}--\eqref{kx02jixdmnfdmc} for
rational $p\in(0,1)$ may be reduced to the $\Gamma$-function and its
logarithmic derivatives. Besides, integral
{\eqref{kx02jixdmnfdmc}}, for any $|\operatorname{Re}{p}|<1$,
may be written in terms of antiderivatives of $\ln\Gamma(z)$ instead
of $\zeta'(-1,z)$. We, however, noticed that currently there is no agree
about the exact definition of $ \Psi
_{-2}(z)\equiv\int\limits\ln\Gamma(z)\, dz $. From the well-known identity $
\zeta'(0,z)=\ln\Gamma(z)-\frac {1}{2}\ln2\pi$ and the fact that $
\partial\zeta(s,z)/\partial z= -s \zeta(s+1,z) $, it clearly follows
that
%
\begin{eqnarray}
\label{j019jxpo1indnd} \Psi_{-2}(z) = \zeta'(-1,z) -
\frac{z^2}{2}+\frac{z}{2} (1+\ln2\pi )+ C
\end{eqnarray}
where $C$ is the constant of integration.\footnote{The Hurwitz $\zeta
$-function whose first
argument is a negative integer may be trivially expressed in terms of
Bernoulli polynomials. In particular
$ \zeta(-1,z)=- \frac{1}{2}z^2+\frac{1}{2}z -\frac{1}{12} $.}
Notwithstanding, we found that Maple 12 uses a different definition
\begin{eqnarray*}
\Psi_{-2}(z) = \zeta'(-1,z) -
\frac{z^2}{2}+\frac{z}{2} -\frac{1}{12}
\end{eqnarray*}
Yet, we remarked that Wolfram Alpha Pro employs another expression,
which numerically corresponds to\footnote{Wolfram Alpha Pro does not
explain how $ \Psi_{-2}(z)$ is evaluated.
Expression given below is derived by the author by
calculating the antiderivative of Binet's integral formula
for the logarithm of the $\Gamma$-function subject to the initial
condition $\Psi_{-2}(0)=0$
(for Binet's formula for $\ln\Gamma(z)$,
see e.g.~\cite[pp.~335--336]{binet_01}, \cite[pp.~250--251]{whittaker_01},
\cite[vol.~I, p.~22, Eq.~1.9(9)]{bateman_01}, \cite[p.~83,
Eq.~(54)]{iaroslav_06}).}
%
\begin{eqnarray}
\nonumber
\Psi_{-2}(z)&=&z\ln\Gamma(z)-
\frac{z^2}{2} \ln{z}+\frac
{z^2}{4}+\frac{z}{2} +
\frac{\ln2\pi}{12} -\frac{1-\gamma}{12} -\frac{\zeta'(2)}{2\pi^2} +\!\int\limits
_0^\infty \frac{ x \ln (x^2+z^2 )}{e^{2\pi x}-1} \,dx
\end{eqnarray}
These three definitions are all different, but it may be easily seen
that first definition {\eqref{j019jxpo1indnd}}
differs from the last one only by a constant of integration, while
Maple's definition is really different.

\subsection{Malmsten's series and Hurwitz's reflection formula}
We now show that the integral from {Lemma \ref{bnhgcewucv}} may be also
evaluated via a trigonometric series.
%
\begin{lemma}\label{bkjcb723bsa}
In the vertical strip $|\operatorname{Re}{a}|<1 $, the following equality holds
%
\begin{eqnarray}
\label{hgu8g7gvijvfdd} \int\limits_0^\infty \frac{ x^{a-1} (\operatorname{ch}{px}-1) }{
\operatorname{sh}{x} } \,dx
= \pi^a \sec\frac{\pi a}{2} \sum_{n=1}^\infty(-1)^n
\frac{ \cos
p\pi n - 1 }{n^{1-a}}
\end{eqnarray}
for $-1<{p}<+1$.
\end{lemma}

\begin{pf}
The Mittag--Leffler theorem is a fundamental theorem in the theory of
functions of a complex variable and allows to expand meromorphic
functions into a series accordingly to its poles.\footnote{For more
details, see \cite{makrushevi_01_eng}, \cite[pp.~147--148, \no
994--1002]{volkovyskii_01_eng}, \cite[Chap.~V, \S27, \no
27.10-2]{evgrafov_01_eng_1st}, \cite[Chap.~VII, p.~175]{spie2},
\cite{lindelof_01}.} Application of this theorem to the meromorphic
function $(\operatorname{ch}{pz}-1)/\operatorname{sh}{z}$, $p\in
(-1,+1)$, having only first-order
poles at $z=\pi n i$, $n\in\mathbbm{Z}$, with residue $(-1)^n (\cos
\pi p
n-1)$, leads to the following expansion
\begin{eqnarray*}
\frac{ \operatorname{ch}{pz}-1 }{ \operatorname{sh}{z} } = 2z\sum_{n=1}^\infty(-1)^n
\frac{ \cos p\pi n -1 }{z^2+\pi^2n^2} , \quad z\in\mathbbm{C} , \ z\neq\pi n i , \ n
\in\mathbbm{Z} ,
\end{eqnarray*}
which is uniformly convergent on the entire complex $z$-plane except
discs $|z-\pi i n|<\varepsilon$, $n\in\mathbbm{Z}$, where the positive parameter $\varepsilon$ can be made as small as we please. Therefore
%
\begin{eqnarray}
\label{h789tv5} %
\int\limits_0^\infty
\frac{ x^{a-1}(\operatorname{ch}{px}-1) }{
\operatorname{sh}{x} } \,dx & = & 2\sum_{n=1}^\infty(-1)^n
(\cos p\pi n -1) \underbrace{\int\limits_0^\infty
\frac{x^a}{x^2+\pi^2n^2} \,dx}_{\frac{1}{2}\pi^a n^{a-1}
\sec\frac{1}{2}\pi a}
\nonumber
\\
& = & \pi^a \sec\frac{\pi a}{2} \sum
_{n=1}^\infty(-1)^n \frac{ \cos
p\pi n - 1 }{n^{1-a}}
\end{eqnarray}
which holds only for $-1<{p}<+1$ and $|\operatorname{Re}{a}|<1$ (the
elementary integral in the middle, whose evaluation is due to Euler,
is convergent only in the strip $|\operatorname{Re}{a}|<1$, see
e.g.~\cite[p.~126, \no880]{volkovyskii_01_eng},
\cite[p.~197, \no856.2]{dwigth_01_en}, \cite[p.~256, \no
6.1.17]{abramowitz_01},
\cite[p.~67, \no587]{gunter_03_eng}, \cite[p.~51]{lindelof_01}).
However, the above equality can be analytically continued for other
values of $a$:
the integral is the analytic continuation of the sum for $
\operatorname{Re}{a}\geqslant+1 $, while the sum analytically
continues the integral for
$ \operatorname{Re}{a}\leqslant-1 $.
We obviously have to expect trouble with the right-hand part at $a=\pm
1, \pm3, \pm5, \ldots$
because of the secant. Since when $ a=-1, -3, -5,\ldots$ the sum in
the right-hand side converges, these points
are poles of the first order for the analytic continuation of integral
{\eqref{h789tv5}}. In contrast, for $a=1, 3, 5,\ldots{}$,
the integral on the left remains bounded, and thus, these points are
removable singularities for the right-hand side of {\eqref{h789tv5}}.
In other words, formally \mbox{$\sum(-1)^n (\cos p\pi n -
1)n^{a-1}$}, $ n\geqslant1 $, must vanish identically
for any odd positive $a$ (exactly as $\eta(1-a)$, the result which has
been derived by Euler, see e.g.~\cite[p.~85]{euler_03}).
These matters are treated in detail in the next corollary.\qed
\end{pf}

\begin{corollary}\label{khciec83cjn}
For $0<{p}<1$
%
\begin{equation}
\everymath{\displaystyle} \label{higce672ecbx}
\left\{\begin{array}{l@{\quad}l}
\sum_{n=1}^\infty\frac{ \cos2\pi p n }{n^{1-a}} =
\Gamma(a) (2\pi)^{-a}\cos\frac{\pi a}{2}  \left\{ \zeta
(a,p )
+ \zeta (a,1-p )  \right\} &(\mathrm{a})\\[16pt]
\sum_{n=1}^\infty\frac{ \sin2\pi p n }{n^{1-a}} =
\Gamma(a) (2\pi)^{-a}\sin\frac{\pi a}{2} \left \{ \zeta
(a,p )
- \zeta (a,1-p )  \right\}&(\mathrm{b})\end{array}\right.
\end{equation}
where both series on the left-hand side are uniformly convergent in
$\operatorname{Re}{a}<1$ and are absolutely convergent
in the half-plane $\operatorname{Re}{a}<0$.
These important formul\ae~seem to be obtained for the first time by
Malmsten in 1846.
\end{corollary}

\begin{pf}
In view of the fact the alternating $\zeta$-function $\eta(s)$ may
be reduced to the ordinary $\zeta$-function
and by making use of Euler--Riemann's reflection formula
for the $\zeta$-function
$ \zeta(1-s)=2\zeta(s)\Gamma(s) (2\pi)^{-s}\cos\frac{1}{2}\pi s$,
we may continue {\eqref{h789tv5}} as follows
\begin{eqnarray*}
\pi^a \sec\frac{\pi a}{2} \sum_{n=1}^\infty(-1)^n
\frac{ \cos
p\pi n - 1 }{n^{1-a}} & =& \pi^a \sec\frac{\pi a}{2} \left\{
\sum_{n=1}^\infty(-1)^n
\frac{
\cos p\pi n }{n^{1-a}} -\left(2^a-1\right)\zeta(1-a)
\right\}
\\[3pt]
&= & \pi^a \sec\frac{\pi a}{2} \sum
_{n=1}^\infty(-1)^n \frac{ \cos p\pi n }{n^{1-a}} - 2
\left(1-2^{-a}\right)\Gamma(a)\zeta(a)
\end{eqnarray*}
Comparing the latter expression to the result of
{Lemma~\ref{bnhgcewucv}} gives
\begin{eqnarray*}
\sum_{n=1}^\infty(-1)^n
\frac{ \cos p\pi n }{n^{1-a}} = \Gamma(a) (2\pi)^{-a}\cos\frac{\pi a}{2}
\left\{ \zeta \left(a,\frac{1}{2}+\frac{p}{ 2 }
\right) + \zeta \left(a,\frac{1}{2}-\frac{p}{ 2 }
\right) \right\}
\end{eqnarray*}
Writing in this expression $2p-1$ instead of $p$ yields immediately
(\ref{higce672ecbx}a).
Now, by partially differentiating
(\ref{higce672ecbx}a)
with respect to $p$ and by remarking that
$a\Gamma(a)=\Gamma(a+1)$, and then, by writing $a$ instead of
$a+1$, we arrive at
(\ref{higce672ecbx}b).
Note also that both sums
(\ref{higce672ecbx}a,b)
may be analytically continued to other domains of $a$ by means of
expressions in corresponding right parts.\qed
\end{pf}

Interestingly, nowadays,
formul\ae~(\ref{higce672ecbx}a,b) seem to be not particularly
well-known (for instance, advanced calculators such as Wolfram Alpha
Pro expresses both series in terms of polylogarithms). Notwithstanding,
Eq.~(\ref{higce672ecbx}b)
can be found in an old Malmsten's work published as early as 1849
\cite[p.~17, Eq.~(48)]{malmsten_01}, and
(\ref{higce672ecbx}a)
is a straightforward consequence of
(\ref{higce672ecbx}b).

\begin{corollary}
If we notice that
\begin{eqnarray*}
\Gamma(a)=\frac{\pi}{\sin\pi a \cdot\Gamma(1-a)}=\frac{\pi}{2
\sin\frac{1}{2}\pi a \cdot\cos\frac{1}{2}\pi a
\cdot\Gamma(1-a)}
\end{eqnarray*}
then, the sum of
(\ref{higce672ecbx}a)
with
(\ref{higce672ecbx}b)
leads to an important formula
%
\begin{eqnarray}
\label{khjx98qh2xbn21} \zeta(a,p) = \frac{2 \Gamma(1-a)}{(2\pi)^{1-a}} \left[ \sin
\frac
{\pi a}{2} \sum_{n=1}^\infty
\frac{ \cos2\pi p n }{n^{1-a}} + \cos\frac{\pi a}{2} \sum_{n=1}^\infty
\frac{ \sin2\pi p n }{n^{1-a}} \right] ,
\end{eqnarray}
with $0<p\leqslant1$ and $\operatorname{Re}{a}<1$, which is usually
attributed to Adolf Hurwitz who derived it in 1881, see
\cite[p.~93]{hurwitz_02},\footnote{\label{lkd024kdfm}
Hurwitz derived all his results
for the function $f(s,a)$ which is related to the modern Hurwitz $\zeta
$-function as
$f(s,a)\equiv f_m(s,a)=m^{-s}\zeta(s,a/m)$, see \cite[p.~89]{hurwitz_02}. By the way, this famous Hurwitz's paper begins
with several factual errors.
The reflection formula for the $L$-function, which he attributed to
Oscar Schl\"omilch \cite[p.~86, first two formul\ae~for $f(s)$]{hurwitz_02},
was first deduced by mathematical induction by Leonhard Euler in 1749
\cite[p.~105]{euler_03}.
Then, it was rigorously proved by two different methods by Malmsten in
1842 \cite{malmsten_00} and in 1846 \cite{malmsten_01}. As regards
Schl\"omilch's contribution, he gave the same formula only in 1849
\cite{schlomilch_01}, and this, without the proof (the proof \cite{schlomilch_02}
was published 9 years later). Similarly, Hurwitz erroneously attributed
the reflection formula for the $\zeta$-function
to Bernhard Riemann, although it was first given also by
Euler \cite[p.~94]{euler_03}, albeit in a slightly form, and
Riemann's contribution consists mainly in the more rigorous proof
of it \cite{riemann_01}.
Further information about the history of these two important formul\ae~may be found in
\cite{wieleitner_01}, \cite[p.~23]{hardy_02}, \cite[p.~861]{davis_01}, \cite[pp.~35--37]{iaroslav_06}.} \cite[p.~269]{whittaker_01},
\cite[p.~107]{lindelof_01}, \cite[p.~37]{titchmarsh_02}, \cite[p.~156]{berndt_02}, \cite[vol.~I, p.~26,
Eq.~1.10(6)]{bateman_01}.\footnote{There is
a slight error in this formula in the latter reference: it remains
valid not
only for $\operatorname{Re}{a}<0$, but also for $\operatorname
{Re}{a}<1$.} Sometimes, it is written in a complex form
\begin{eqnarray*}
\zeta(a,p) = \frac{i \Gamma(1-a)}{(2\pi)^{1-a}} \left[ e^{-\frac
{1}{2}\pi ia} \sum
_{n=1}^\infty\frac{ e^{-2\pi i p n }}{n^{1-a}} -
e^{+\frac
{1}{2}\pi ia} \sum_{n=1}^\infty
\frac{ e^{+2\pi i p n }}{n^{1-a}} \right] ,
\end{eqnarray*}
$0<p\leqslant1$, $\operatorname{Re}{a}<1$,
see e.g.~\cite[p.~87]{chudakov_01_eng}, which is completely
equivalent to {\eqref{khjx98qh2xbn21}}.
\end{corollary}
\textbf{Nota Bene.}
It is quite rarely emphasized that latter representations coincide with
the trigonometric Fourier series for $\zeta(a,p)$. Remarking this
permits to immediately derive several
integral formul\ae, whose demonstration by other means is more difficult
%
\begin{eqnarray}
\everymath{\displaystyle}
\label{k9230920jhdsn} 
\begin{cases}
\int\limits_0^1 \zeta(a,p)\, dp = 0 \vspace{2pt}\cr
\int\limits_0^1 \zeta(a,p) \cos2\pi p n\, dp = \Gamma(1-a) (2\pi
n)^{a-1}\sin\frac{\pi a}{2} \vspace{2pt}\cr
\int\limits_0^1 \zeta(a,p) \sin2\pi p n\, dp = \Gamma(1-a) (2\pi
n)^{a-1}\cos\frac{\pi a}{2}
\end{cases}
\quad \operatorname{Re} {a}<1
\end{eqnarray}
for $n=1, 2, 3,\ldots{}$.
Furthermore, in virtue of Parseval's theorem, we also have
%
\begin{eqnarray}
\label{lkd92jr4n} \int\limits_0^1 \zeta^2(a,p)
dp = 2\Gamma^2(1-a) (2\pi)^{2a-2}\zeta(2-2a) , \quad
\operatorname{Re} {a}<1 , \ a\neq\frac{1}{2}
\end{eqnarray}
Differentiating this formula with respect to $a$ and then setting
$a=0$, yields:
\begin{eqnarray*}
2 \!\int\limits_0^1 \!\underbrace{
\left(\frac{1}{2}-p \right)}_{\zeta(0,p)} \cdot
\underbrace{ \left(\ln\Gamma(p)-\frac{1}{2}\ln2\pi
\right)}_{\zeta'(0,p)} dp = \frac{\gamma+\ln2\pi}{6} - \frac{\zeta'(2) }{\pi^2}
\end{eqnarray*}
Whence, accounting for the well-known result\footnote{The value of
integral {\eqref{jn893dhwsa}}, as well as that of
{\eqref{jn893dhjhx8783}}, may be both straightforwardly deduced from a similar
Fourier series expansion for the logarithm of the $\Gamma$-function,
see e.g.~\cite[vol.~I, pp.~23--24, \S1.9.1]{bateman_01} or
\cite[p.~17, Eq.~(36)]{srivastava_03}.
This expansion, attributed erroneously to Ernst Kummer,
was first derived by Malmsten and colleagues from the Uppsala
University in 1842.
This interesting historical question is discussed in details in \cite[Sect.~2.2, Fig.~2 and exercise \no20 on pp.~66--68]{iaroslav_06}.
By the way, the evaluation of integral {\eqref{jn893dhjhx8783}} may be
also found in several modern works,
see e.g.~\cite[p.~177, Eq.~(7.3)]{espinosa_01}, \cite[p.~14,
Eq.~(3.19)]{connon_04}.}
%
\begin{eqnarray}
\label{jn893dhwsa} \int\limits_0^1 \ln\Gamma(p)\, dp =
\frac{1}{2}\ln2\pi
\end{eqnarray}
we obtain
\begin{eqnarray*}
\int\limits_0^1 \!p \ln\Gamma(p)\, dp =
\frac{\zeta'(2)}{2\pi^2} - \frac
{\gamma-2\ln2\pi}{12}
\end{eqnarray*}
Integration by parts of the latter expression leads to the
antiderivatives of $\ln\Gamma(x)$
which are currently not well-studied yet (see the \emph{Nota Bene} on
p.~\pageref{hd91283dn9323hnd}).
Similarly, differentiating twice {\eqref{lkd92jr4n}} with respect to $a$
at $a=0$, and accounting for\setcounter{footnote}{20}\footnotemark
%
\begin{eqnarray}
\label{jn893dhjhx8783} %
\int\limits_0^1
\ln^2 \Gamma(p)\, dp &= & \frac{\gamma^2}{12} + \frac{\pi^2}{48} +
\frac{\gamma\ln2\pi}{6} + \frac{\ln^2 2\pi}{3} - \frac
{(\gamma+\ln2\pi) \zeta'(2)}{\pi^2} +
\frac{\zeta''(2)}{2\pi^2}
\nonumber
\\
& = & \frac{1}{6}+ \frac{\pi^2}{36} +\frac{\ln^2 2\pi}{4} -2
\zeta'(-1)-\zeta''(-1)
\end{eqnarray}
yields another integral\setcounter{footnote}{21}%
\begin{eqnarray*}
\int\limits_0^1 p \zeta''(0,p)
\,dp &= & \frac{\pi^2}{144} - \frac{\gamma^2}{12} - \frac{\gamma\ln2\pi}{6} -
\frac{\ln^2 2\pi}{12} + \frac{(\gamma+\ln2\pi) \zeta'(2)}{\pi
^2} - \frac{\zeta''(2)}{2\pi^2}
\\
& = & -\frac{1}{6}+2\zeta'(-1)+\zeta''(-1)
\end{eqnarray*}
Some further results related to the Fourier series expansion of the
Hurwitz $\zeta$-function
are provided in \cite{espinosa_01}.\footnote{However, in many
formul\ae~domains of validity remain unspecified, and sometimes, are
incorrect
(e.g.~compare \cite[Eqs.~(1.26) and (3.5)]{espinosa_01}
with {\eqref{k9230920jhdsn}} and {\eqref{lkd92jr4n}} respectively).}

\begin{corollary} \label{kluwch92hd}
In {\eqref{khjx98qh2xbn21}}, the index $n$ may be represented as
$n=mk+l$, where for
each $k=0, 1, 2, \ldots, \infty$, the index $l$ runs over $[1, 2,
\ldots, m]$ and where $m$ is some positive integer. Then,
{\eqref{khjx98qh2xbn21}} may be written in the form:
%
\begin{eqnarray}
\;\,\zeta(a,p) = \frac{2 \Gamma(1-a)}{(2\pi)^{1-a}} \left[ \sin
\frac
{\pi a}{2} \sum_{l=1}^m \sum
_{k=0}^\infty\frac{ \cos2\pi p (mk+l)
}{(mk+l)^{1-a}} + \cos
\frac{\pi a}{2} \sum_{l=1}^m \sum
_{k=0}^\infty\frac{ \sin2\pi p (mk+l)
}{(mk+l)^{1-a}} \right]\nonumber
\end{eqnarray}
Now, let $p$ be a rational part of $m$, i.e.~$p=r/m$, where $r$ and $m$
are positive integers such that $r\leqslant m$. Then $\cos [ 2\pi p
(mk+l) ]=\cos(2\pi rl/m)$, and similarly for the sine. Hence, for
positive rational $p$ not greater than 1, the previous formula
takes the form
\be\label{khjhbcubdbvxhbs}
\begin{array}{ll}
\displaystyle
 \zeta\!\left(a,\frac{r}{m}\right)&\displaystyle = \,\frac{2 \Gamma(1-a)}{(2\pi)^{1-a}}\Biggl[ \sin\frac{\pi a}{2}
\sum_{l=1}^m \cos\frac{2\pi rl}{m} \underbrace{\sum_{k=0}^\infty\frac{\, 1\,}{(mk+l)^{1-a}}}_{m^{a-1}\zeta(1-a,\,l/m)} + \\[10mm]
&\displaystyle\qquad\qquad\qquad\qquad\qquad\qquad\qquad\qquad
+ \cos\frac{\pi a}{2}
\sum_{l=1}^m \sin \frac{2\pi rl}{m} \sum_{k=0}^\infty \frac{\,1 \,}{(mk+l)^{1-a}}\Biggr]  \\[7mm]
& \displaystyle 
= \frac{2 \Gamma(1-a)}{(2\pi m)^{1-a}}
\sum_{l=1}^m \sin\!\left(\frac{2\pi rl}{m} + \frac{\pi a}{2} \!\right)\!\cdot\zeta\!\left(\!1-a,\,\frac{l}{m}\!\right) ,
\qquad\qquad r=1, 2,\ldots, m.
\end{array}
\ee
This equality holds in the entire complex $a$-plane for any positive
integer $m\geqslant2$. Furthermore, by putting in the latter formula
$1-a$ instead of $a$, it may be rewritten as
%
\begin{eqnarray}
\label{kljc928c2ndbn} \zeta \left(1-a,\frac{r}{m} \right) =
\frac{2 \Gamma(a)}{(2\pi m)^{a}} \sum_{l=1}^m \cos
\left(\frac{2\pi rl}{m} - \frac{\pi a}{2} \right) \cdot\zeta
\left( a, \frac{l}{m} \right) , \quad r=1, 2,\ldots, m.
\end{eqnarray}
In the case $r=m$, the above formul\ae~reduce
to Euler--Riemann's reflection formul\ae~for the $\zeta$-function
(simply use the multiplication theorem for
the Hurwitz $\zeta$-function, see e.g.~\cite[p.~101]{iaroslav_06}).
Formulae~{\eqref{khjhbcubdbvxhbs}} and {\eqref{kljc928c2ndbn}} are
known as functional equations for the Hurwitz $\zeta$-function and
were both
obtained by Hurwitz in the same article \cite[p.~93]{hurwitz_02} in 1881.
By the way, the above demonstration also shows that they can be
elementary derived from Malmsten's results (\ref{higce672ecbx}a,b)
obtained as early as 1840s.
\end{corollary}
\textbf{Nota Bene.} Malmsten's series
(\ref{higce672ecbx}a,b)
are actually particular cases of a more general series
%
\begin{eqnarray}
\label{pio2ux092mdd} f(s) \equiv\sum_{n=1}^\infty
\frac{ a_n}{n^s} ,\quad a_n\in\mathbbm{C},\; |a_n|\leqslant1 ,
\end{eqnarray}
which is uniformly and absolutely convergent in the region $\operatorname{Re}{s}>1$
(it may also converge, albeit non-absolutely, in the half-plane
$\operatorname{Re}{s}>0$).\footnote{Indeed $\,\sum |a_n n^{-s}|\leqslant\sum |n^{-s}|=\zeta\left(\operatorname{Re}{s}\right)$, the latter being uniformly and absolutely convergent in \mbox{$\operatorname{Re}{s}>1$}.} Such a series is known as the
\emph{Dirichlet series}. Let now focus our attention on a particular
case of this series in which coefficients $a_n$ are $m$-periodic,
i.e.~$a_{n}=a_{n+m}=a_{n+2m}=\ldots $ (period $m$ being
natural).\footnote{If $a_n$ is a character, the above series may be, in
turn, an example of the Dirichlet $L$-function.} The first important consequence
of such a particular case is that $f(s)$ may be reduced to a linear
combination of Hurwitz $\zeta$-functions at rational argument.
Representing again the
summation's index $n=mk+l$ yields
%
\begin{eqnarray}
\label{jhcwiknx83} f(s) = \sum_{l=1}^m \sum
_{k=0}^\infty\frac{ a_{mk+l} }{(mk+l)^{s}} = \sum
_{l=1}^m a_l \sum
_{k=0}^\infty\frac{ 1 }{(mk+l)^{s}} = \frac{1}{m^s}
\sum_{l=1}^m a_l \zeta
\left( s,\frac{l}{m} \right)
\end{eqnarray}
The right-hand side continues $f(s)$ to the entire complex $s$-plane,
except possibly the point $s=1$.\footnote{This is the unique point
where the Hurwitz $\zeta$-function is not regular.} In order to
identify the character of the point $s=1$, we evaluate the
corresponding residue
\begin{eqnarray*}
\mathop{\mathrm{res}}\limits_{s=1} f(s)= \lim_{s\to1}
\left[(s-1)f(s) \right] = \frac{1}{m}\sum
_{l=1}^m a_l \lim_{s\to1}
\left[(s-1)\zeta \left( s,\frac{l}{m} \right)
\right] = \frac{1}{m} \sum_{l=1}^m
a_l
\end{eqnarray*}
Therefore, if the coefficients $a_1, a_2,\ldots, a_m$ are chosen so
that the latter sum vanishes, then
$f(s)$ is holomorphic; otherwise $f(s)$ is a meromorphic function with
a unique pole at $s=1$. Typical examples
of cases when $f(s)$ is regular everywhere are Malmsten's series
(\ref{higce672ecbx}a,b) at rational $p$ because
\begin{eqnarray*}
\sum_{l=1}^m a_l = \sum
_{l=1}^m \sin\frac{2\pi rl}{m} = 0 \quad
\mbox{and} \quad \sum_{l=1}^m \cos
\frac{2\pi rl}{m} = 0\,,\quad r=1, 2,\ldots, m-1.
\end{eqnarray*}

As to the reflection formula for the Dirichlet series $f(s)$, it may be
easily deduced with the help of {\eqref{kljc928c2ndbn}}.
Writing in {\eqref{jhcwiknx83}} $1-s$ for $s$, and then, using
{\eqref{kljc928c2ndbn}}, yields
%
\begin{eqnarray}
\label{io2d029jdmd} %
f(1-s) & = & \frac{1}{m^{1-s}} \sum
_{l=1}^m a_l \zeta \left( 1-s,
\frac{l}{m} \right) =\frac{2 \Gamma(s)}{m(2\pi)^{s}} \sum
_{l=1}^m a_l \sum
_{k=1}^m \cos \left(\frac{2\pi lk}{m} -
\frac
{\pi s}{2} \right) \cdot\zeta \left( s, \frac{k}{m}
\right)
\nonumber
\\[6pt]
& = & \frac{2 \Gamma(s)}{m(2\pi)^{s}} \left[ \sin\frac{\pi s}{2} \sum
_{k=1}^m \alpha_k \zeta \left( s,
\frac{k}{m} \right) + \cos \frac{\pi s}{2} \sum
_{k=1}^m \beta_k \zeta \left( s,
\frac{k}{m} \right) \right]
\end{eqnarray}
where
\begin{eqnarray*}
\alpha_k = \sum_{l=1}^m
a_l \sin\frac{2\pi lk}{m} \quad\mbox{and} \quad
\beta_k = \sum_{l=1}^m
a_l \cos\frac{2\pi lk}{m}
\end{eqnarray*}
holding in the entire complex $s$-plane except at points $s=1, 0, -1,
-2, \ldots{}$. This formula is also very useful in that the expression on
the right represents the analytic continuation for $f(1-s)$ to the
domains where the series {\eqref{pio2ux092mdd}} does not
converge. Finally, remark that the latter formula may be also written
in a complex form
\begin{eqnarray*}
f(1-s) = \frac{\Gamma(s)}{m(2\pi)^{s}} \left[ e^{+\frac{1}{2}\pi i s} \sum
_{k=1}^m \widetilde{\alpha}_k \zeta
\left( s, \frac{k}{m} \right) +e^{-\frac{1}{2}\pi i s} \sum
_{k=1}^m \widetilde{\beta}_k
\zeta \left( s, \frac{k}{m} \right) \right]
\end{eqnarray*}
$s\neq1, 0, -1, -2, \ldots{}$,
where
\begin{eqnarray*}
\widetilde{\alpha}_k = \sum_{l=1}^m
a_l e^{-2\pi i lk/m} \quad\mbox{and} \quad \widetilde{
\beta}_k = \sum_{l=1}^m
a_l e^{+2\pi i lk/m}
\end{eqnarray*}
and some authors precisely prefer this form, see
e.g.~\cite[pp.~88--91]{chudakov_01_eng}. This form is more appropriated
if one wishes to emphasize the Fourier series aspect (coefficients
$\widetilde{\alpha}_k$ and $\widetilde{\beta}_k$ may be regarded as
$m$-points Fourier transforms of coefficients $a_l$).

\subsection{Closed-form evaluation of the first generalized Stieltjes
constant~at rational argument}

We now state the main result of this manuscript allowing to evaluate in
a closed-form the first generalized Stieltjes constant~at any rational
argument.

\begin{teorema}
The first generalized Stieltjes constant~of any rational argument in
the range $(0,1)$ may be expressed in a closed form via a finite
combination of
logarithms of the $\Gamma$-function, of second-order derivatives of
the Hurwitz $\zeta$-function
at zero, of Euler's constant $\gamma$, of the first Stieltjes
constant~$\gamma_1$ and of elementary functions:
\be\label{hg784bdn}
\begin{array}{l}
\displaystyle
\gamma_1 \biggl(\frac{r}{m} \vphantom{\frac{1}{2}} \biggr)
 = \displaystyle\,
\gamma_1 - \gamma\ln2 m -\ln^2\!{2} - \ln2\cdot\ln{\pi m} -\frac{1}{2}\ln^2\!{m}
- \frac{\pi}{2}(\gamma+\ln2\pi m)\ctg\frac{\pi r}{m}  \\[6mm]
\displaystyle \quad
-\frac{(-1)^r}{4}\big[1-(-1)^{m+1}\big]\!\cdot  (3\ln2+2\ln\pi)\ln2   
- \pi\ln\pi\cdot
\csc\frac{\pi r}{m}\cdot
\sin\!\left(\!\frac{\pi r}{m}\left\lfloor\!\frac{m+1}{2}\!\right\rfloor \right)\times\\[6mm]
\displaystyle \quad
\!\times\sin\!\left(\!\frac{\pi r}{m}\left\lfloor\!\frac{m-1}{2}\!\right\rfloor \right) \!
+ 2(\gamma+\ln2\pi m)\cdot\!\!\!\!\!\!\!\sum_{l=1}^{\lfloor\!\frac{1}{2}(m-1)\!\rfloor} \!\!\!\!\!
\cos\frac{2\pi rl}{m}\cdot\ln\sin\frac{\pi l}{m}
+ \pi \!\!\!\!\!\!\sum_{l=1}^{\lfloor\!\frac{1}{2}(m-1)\!\rfloor} \!\!\!\!\!\!
\sin\frac{2\pi rl}{m}  \cdot\ln\sin\frac{\pi l}{m}     \\[6mm]
\displaystyle  \quad
+ 2\pi\!\!\!\!\!\!\sum_{l=1}^{\lfloor\!\frac{1}{2}(m-1)\!\rfloor} \!\!\!\!\!
\sin\frac{2\pi rl}{m}  \cdot\ln\Gamma\!\left(\!\frac{l}{m}\!\right) 
+ \!\!\!\!\!\!\sum_{l=1}^{\lfloor\!\frac{1}{2}(m-1)\!\rfloor} \!\!\!\!\!
\cos\frac{2\pi rl}{m}  \cdot\left\{\zeta''\!\left(\! 0,\,\frac{l}{m}\!\right) +\,
\zeta''\!\left(\! 0,\,1-\frac{l}{m}\!\right) \!\right\} 
\end{array}
\ee
This elegant formula holds for any $r=1, 2, 3,\ldots, m-1$, where $m$
is positive integer greater than 1.
The Stieltjes constants~for other ``periods'' may be obtained from the
recurrent relationship:
%
\begin{eqnarray}
\label{asjd83a3} \gamma_1(v+1) = \gamma_1(v) -
\frac{ \ln v }{v} , \quad v\neq0 ,
\end{eqnarray}
see, e.g.~\cite[p.~102, Eq.~(64)]{iaroslav_06}.
The above theorem is an equivalent of Gauss' Digamma theorem for
the $0th$ Stieltjes constant~$\gamma_0(r/m)=-\Psi(r/m)$.
Three alternative forms of the same theorem are given in Eqs.
{\eqref{kx98eujms}}, {\eqref{khjwe9jdns}}
and {\eqref{kjkhnc823h}} respectively.
\end{teorema}

\begin{pf}
Consider the integral {\eqref{lkij289dnxnwx}}. Put $2p-1$ instead of $p$
and denote the resulting integral via $J_a(p)$:
\be\label{lkijihgcqwiwx}
J_a(p)\equiv\int\limits_0^\infty \!\!\frac{\,x^{a-1} (\ch{[(2p-1)x]}-1)\,}{\,\sh{x}\,}\,dx\,=\,
\frac{\Gamma(a)}{2^{a}}\biggl\{\zeta(a,p)+ \zeta(a,1-p) -2\big(2^{a}-1\big)\zeta(a) \biggr\}
\ee
converging in the strip $0<\operatorname{Re}{p}<1$.
Let now $p$ be rational $p=r/m$, where $r$ and $m$ are positive
integers such that $r<m$.
Then, the preceding equation becomes
%
\begin{eqnarray}
\label{lkijihgcqwiwx2} J_a \left(\frac{r}{m} \right)=
\frac{\Gamma(a)}{2^{a}} \left\{\zeta \left( a, \frac{r}{m}
\right) + \zeta \left( a, 1- \frac{r}{m} \right) - 2
\left(2^{a}-1 \right)\zeta(a) \right\}
\end{eqnarray}
The sum of first two terms in curly brackets may be evaluated via
Hurwitz's reflection formula {\eqref{khjhbcubdbvxhbs}}:
\be\notag
\begin{array}{l}
\displaystyle
\zeta\biggl(\!a, \frac{r}{m} \vphantom{\frac{1}{2}} \biggr) +
\zeta\biggl(\!a, 1-\frac{r}{m} \vphantom{\frac{1}{2}} \biggr) =\,
\frac{2 \Gamma(1-a)}{(2\pi m)^{1-a}} \sum_{l=1}^m 
\left[\sin\!\left(\frac{2\pi rl}{m} + \frac{\pi a}{2} \!\right)\!
+\sin\!\left(\frac{2\pi (m-r)l}{m} + \frac{\pi a}{2} \!\right) \!\right]\\[8mm]
\displaystyle \qquad\qquad\qquad\qquad\qquad
\times \zeta\!\left(\!1-a,\,\frac{l}{m}\!\right) = \,\frac{4 \Gamma(1-a)}{(2\pi m)^{1-a}}\,\sin\frac{\pi a}{2} \cdot
\sum_{l=1}^m \cos\frac{2\pi rl}{m}\cdot\zeta\!\left(\!1-a,\,\frac{l}{m}\!\right) 
\end{array}
\ee
Thus, by noticing that $ \Gamma(a)\Gamma(1-a)=\frac{1}{2}\pi\csc
\frac{1}{2}\pi a \cdot\sec\frac{1}{2}\pi a $,
the integral $J_a(r/m)$ takes the form:
%
\begin{eqnarray}
\label{khd287chejbcx} J_a \left(\frac{r}{m} \right) &= &
\underbrace{\frac{\pi}{(\pi m)^{1-a}} \sec\frac{\pi
a}{2}}_{f_1} \cdot\!
\underbrace{\sum_{l=1}^m \cos
\frac{2\pi rl}{m}\cdot\zeta \left( 1-a, \frac{l}{m}
\right) }_{f_2} - \underbrace{\frac{\Gamma(a) (2^{a}-1 )\zeta(a) }{2^{a-1}}}_{f_3}\qquad\;
\end{eqnarray}
which is third expression for the integral $J_a$, other two expressions
being given by
{\eqref{lkij289dnxnwx}} and {\eqref{hgu8g7gvijvfdd}}.
Let now study each term of the right part, denoted for brevity $f_1$,
$f_2$ and $f_3$ respectively, in a neighborhood of $a=1$.
The first and the third terms have poles of the first order at this
point, while the second term $f_2$
is analytic at $a=1$. Thus, in a neighborhood of $a=1$, terms $f_1$
and $f_3$ may be expanded in the Laurent series
as follows
%
\begin{eqnarray}
\label{jhxc28dn1} f_1 = -\frac{2}{a-1}-2\ln\pi m - \left(
\frac{\pi^2}{12}+\ln^2\pi m \right) \cdot(a-1) +
O(a-1)^2
\end{eqnarray}
and
%
\begin{eqnarray}
\label{jhxc28dn3} f_3 = \frac{1}{a-1} + \ln2 + \left(
\frac{\pi^2}{12}-\frac{\ln^2 2}{2}-\frac{\gamma
^2}{2}-\gamma_1
\right) \cdot(a-1) + O(a-1)^2
\end{eqnarray}
while $f_2$ may be represented by the following Taylor series
%
\begin{eqnarray}
\label{jhxc28dn2} && \!\!\!f_2 =  \sum_{l=1}^m
\cos\frac{2\pi rl}{m}\cdot\underbrace{\zeta \left( 0,
\frac{l}{m} \right)}_{\frac{1}{2}-l/m} - (a-1) \sum
_{l=1}^m \cos\frac{2\pi rl}{m}\cdot \!\!\!\!\!\!\underbrace{
\zeta' \left( 0, \frac{l}{m} \right)}_{\ln\Gamma(l/m)-\frac{1}{2}\ln2\pi}\!\!\!\!\!\!+
\nonumber
\\
&&\qquad\qquad\qquad\qquad\qquad\qquad + \frac{(a-1)^2}{2} \sum_{l=1}^m
\cos\frac{2\pi rl}{m}\cdot\zeta'' \left( 0,
\frac{l}{m} \right) + O(a-1)^3
\\[2pt]
&&\!= -\frac{1}{2} - (a-1) \sum_{l=1}^m
\cos\frac{2\pi rl}{m}\cdot\ln \Gamma \left( \frac{l}{m}
\right)
+ \frac{(a-1)^2}{2} \sum_{l=1}^m
\cos\frac{2\pi rl}{m}\cdot\zeta '' \left( 0,
\frac{l}{m} \right) + O(a-1)^3
\nonumber
\end{eqnarray}
because
%
\begin{eqnarray}
\everymath{\displaystyle}
\label{hjg2873gdb} 
\begin{cases}
\sum_{l=1}^{m} \cos\frac{2\pi rl}{m} = 0 & r=1, 2, 3,\ldots,
m-1\vspace{6pt}\cr
\sum_{l=1}^{m} l\cdot\cos\frac{2\pi rl}{m} = \frac{m}{2} ,
 & r=1, 2, 3,\ldots, m-1
\end{cases}
\end{eqnarray}
In the final analysis, the substitution of {\eqref{jhxc28dn1}}, {\eqref
{jhxc28dn3}} and {\eqref{jhxc28dn2}} into {\eqref{khd287chejbcx}},
yields the following representation for the integral $J_a(r/m)$ in a
neighborhood of $a=1$:
\be\label{kjhg87fgf}
\begin{array}{l}
\displaystyle
J_a\biggl(\frac{r}{m} \vphantom{\frac{1}{2}} \biggr) =\, 
\ln\frac{\pi m}{2} + 2A_m(r) + (a-1)\cdot\!\left[-B_m(r) + 2A_m(r)\ln\pi m - \frac{\pi^2}{24}+\frac{\ln^2\!\pi m}{2} \right.\\[6mm]
\displaystyle \qquad\qquad\qquad\qquad\qquad\qquad\qquad\qquad\quad\;\;
\left. + 
\frac{\gamma^2}{2}+\frac{\ln^2\!{2}}{2} +\gamma_1\right]  + O(a-1)^2
\end{array}
\ee
where\vspace{-4pt}
\begin{eqnarray*}
\everymath{\displaystyle}
\begin{cases}
A_m(r) \equiv\sum_{l=1}^m \cos\frac{2\pi rl}{m}\cdot\ln\Gamma
 \left( \frac{l}{m}  \right) \vspace{6pt}\cr
B_m(r) \equiv\sum_{l=1}^m \cos\frac{2\pi rl}{m}\cdot\zeta''
\left( 0, \frac{l}{m}  \right)
\end{cases}
\end{eqnarray*}
Now, if we look at the integral $J_a(r/m)$ defined in {\eqref{lkijihgcqwiwx}},
we see that it is uniformly convergent and regular near $a=1$,
and hence, may be expanded in the following Taylor series
\be\label{kj3fh78fhh}
J_a(r/m) \, =\, J_1(r/m) + (a-1)\left.\frac{\partial J_a(r/m)}{\partial a}\right|_{a=1} + O(a-1)^2
\ee
Equating right-hand sides of {\eqref{kjhg87fgf}} and
{\eqref{kj3fh78fhh}}, and then,
searching for terms with same powers of $(a-1)$, gives
\begin{eqnarray*}
\everymath{\displaystyle}
\begin{cases}
\int\limits_0^\infty \frac{ \operatorname{ch}{[(2p-1)x]}-1 }{
\operatorname{sh}{x} } \,dx = \ln\frac{\pi m}{2} + 2A_m(r) \vspace
{2pt}\cr
\int\limits_0^\infty \frac{ (\operatorname{ch}{[(2p-1)x]}-1)\ln{x} }{
\operatorname{sh}{x} } \,dx =
\gamma_1-B_m(r) + 2A_m(r)\ln\pi m - \frac{\pi^2}{24}+\frac{\ln^2
\pi m}{2}\vspace{-2pt}\cr
\qquad\qquad\qquad\qquad\qquad\qquad\qquad
+\frac{\ln^2 {2}}{2} + \frac{\gamma^2}{2}
\end{cases}
\end{eqnarray*}
where $p\equiv r/m$. Remarking that the reflection formula for the
logarithm of the $\Gamma$-function
reduces the sum $A_m(r)$ to elementary functions\footnote{By using
Malmsten's representation for the Digamma function,
see (\ref{dh38239djws}c), the sum $A_m(r)$ may be also written in
terms of the $\Psi$-function and Euler's constant $\gamma$.
}
\begin{eqnarray}
\nonumber
A_m(r) \equiv\sum_{l=1}^m
\cos\frac{2\pi rl}{m}\cdot\ln\Gamma \left( \frac{l}{m}
\right) = -\frac{1}{2} \left\{\ln\pi+ \sum
_{l=1}^{m-1} \cos\frac{2\pi
rl}{m} \cdot\ln\sin
\frac{\pi l}{m} \right\}
\end{eqnarray}
yields for the first integral
\begin{eqnarray*}
\int\limits_0^\infty \frac{ \operatorname{ch}{[(2p-1)x]}-1 }{
\operatorname{sh}{x} } \,dx =
\ln\frac{ m}{2} - \sum_{l=1}^{m-1}
\cos\frac{2\pi rl}{m} \cdot\ln \sin\frac{\pi l}{m} , \quad p\equiv
\frac{r}{m}
\end{eqnarray*}
while the second one reads
%
\begin{eqnarray}
\label{lqjhc2908hx} %
 \int\limits_0^\infty
\frac{ (\operatorname{ch}{[(2p-1)x]}-1)\ln{x} }{
\operatorname{sh}{x} } \,dx &= & \ln^2 2 + \ln2\cdot\ln\pi+
\frac{1}{2} \ln^2 {m} -\sum_{l=1}^{m-1}
\cos\frac{2\pi rl}{m}\cdot\zeta'' \left( 0,
\frac{l}{m} \right)
\nonumber
\\[-1pt]
&&{} - \ln\pi m \cdot \sum_{l=1}^{m-1} \cos
\frac{2\pi rl}{m} \cdot\ln \sin\frac{\pi l}{m} , \quad p\equiv
\frac{r}{m}
\end{eqnarray}
where, at the final stage, we separate the last term in the sum
$B_m(r)$ whose value is known 
$
\,\zeta''(0,1)=\zeta''(0)=
\gamma_1+\frac{1}{2}\gamma^2 - \frac{1}{24}\pi^2 -
\frac{1}{2}\ln^2 2\pi
$.
But the integral {\eqref{lqjhc2908hx}} was also evaluated in
{\eqref{k290djd2s}} by means of first generalized Stieltjes constants. Hence,
the comparison of {\eqref{k290djd2s}} to
{\eqref{lqjhc2908hx}} yields
\be\label{d33f0rfifm}
\begin{array}{l}
\displaystyle
\gamma_1 \biggl(\frac{r}{m} \vphantom{\frac{1}{2}} \biggr)
+ \gamma_1 \biggl(1-\frac{r}{m} \vphantom{\frac{1}{2}} \biggr)  =\, 2\gamma_1 
-2\gamma\ln2m +
2\!\sum_{l=1}^{m-1} \cos\frac{2\pi rl}{m}\cdot\zeta''\!\left(\! 0,\,\frac{l}{m}\!\right) 
\\[6mm]
\displaystyle  \qquad\quad
+ 2(\gamma+\ln2\pi m)\! \sum_{l=1}^{m-1} \cos\frac{2\pi rl}{m} \cdot\ln\sin\frac{\pi l}{m}
 - 2\ln^2\!2 - 2\ln2\cdot\ln{\pi m} -  \ln^2\!m 
\end{array}
\ee
for each $r=1, 2,\ldots, m-1$.
Adding this to Malmsten's reflection formula for the first generalized Stieltjes constant {\eqref{hgxcw3b2}} finally gives
\be\label{kx98eujms}
\begin{array}{lll}
&&\displaystyle\gamma_1   \biggl(\frac{r}{m} \vphantom{\frac{1}{2}} \biggr)= 
\,\gamma_1 - \gamma\ln2m+ \sum_{l=1}^{m-1} \cos\frac{2\pi rl}{m}\cdot\zeta''\!\left(\! 0,\,\frac{l}{m}\!\right) 
+\pi\sum_{l=1}^{m-1} \sin\frac{2\pi r l}{m} \cdot\ln\Gamma \biggl(\frac{l}{m} \biggr) 
   - \ln^2\!2 
\\[6mm]
&&\displaystyle\quad
- \frac{\pi}{2}(\gamma+\ln2\pi m)\ctg\frac{\pi r}{m}
+ (\gamma+\ln2\pi m)\!\sum_{l=1}^{m-1} \cos\frac{2\pi rl}{m} \cdot\ln\sin\frac{\pi l}{m} 
- \ln2\cdot\ln{\pi m} - \frac{1}{2}\ln^2\!{m}\qquad\qquad
\end{array}
\ee
This is the most simple form of the theorem which we are stating here
and can be used as is.
It can be also written in several other forms. For instance, one may
notice that
sums over $l\in[1,m-1]$ may be further simplified. Since each pair of
terms which occupy symmetrical positions
relatively to the center (except for $l=m/2$ when $m$ is even) may be
grouped together, the first sum may be reduced to
\be\label{oi2d298dehd1}
\begin{array}{ll}
\displaystyle
\sum_{l=1}^{m-1} \cos\frac{2\pi rl}{m}\cdot\zeta''\!\left(\! 0,\,\frac{l}{m}\!\right) = \,
\frac{1}{2}\!\sum_{l=1}^{m-1} \!\cos\frac{2\pi rl}{m}\cdot\left\{\zeta''\!\left(\! 0,\,\frac{l}{m}\!\right) +\,
\zeta''\!\left(\! 0,\,1-\frac{l}{m}\!\right) \!\right\}   \\[6mm]
\quad =  \displaystyle 
\begin{cases}
\displaystyle \sum_{l=1}^{\frac{1}{2}(m-1)} \!\!\! \cos\frac{2\pi rl}{m}\cdot\left\{\zeta''\!\left(\! 0,\,\frac{l}{m}\!\right) +\,
\zeta''\!\left(\! 0,\,1-\frac{l}{m}\!\right) \!\right\}\,, \quad &\text{if $m$ is odd} \\[6mm]
\displaystyle \sum_{l=1}^{\frac{1}{2}m-1} \!\!\! \cos\frac{2\pi rl}{m}\cdot\left\{\zeta''\!\left(\! 0,\,\frac{l}{m}\!\right) +\,
\zeta''\!\left(\! 0,\,1-\frac{l}{m}\!\right) \!\right\} + (-1)^r\zeta''\!\left(\! 0,\,\frac{1}{2}\right)\,, 
\quad & \text{if $m$ is even} 
\end{cases} \\[16mm]
\quad =  \displaystyle  \!\!\!\!\!\!\sum_{l=1}^{\lfloor\!\frac{1}{2}(m-1)\!\rfloor} \!\!\!\!
\cos\frac{2\pi rl}{m}  \cdot\left\{\zeta''\!\left(\! 0,\,\frac{l}{m}\!\right) +\,
\zeta''\!\left(\! 0,\,1-\frac{l}{m}\!\right) \!\right\} - \\[6mm]
\qquad\qquad\qquad\qquad\qquad\qquad\qquad\qquad\qquad
-  \displaystyle 
\frac{(-1)^r}{4}\big[1-(-1)^{m+1}\big]\!\cdot  (3\ln2+2\ln\pi)\ln2
\end{array}
\ee
because $\zeta''\left(0,\frac{1}{2} \right)= -\frac{3}{2}\ln^2 2 -
\ln\pi\ln2$, see e.g.~\cite[p.~72, \no24]{iaroslav_06}.
Analogously, the second sum may be written as
\be\label{oi2d298dehd2}
\begin{array}{l}
\displaystyle
\sum_{l=1}^{m-1} \sin\frac{2\pi rl}{m}\cdot\ln\Gamma\!\left(\!\frac{l}{m}\!\right)
= \!\!\!\!\!\!\sum_{l=1}^{\lfloor\!\frac{1}{2}(m-1)\!\rfloor} \!\!\!\!
\sin\frac{2\pi rl}{m}  \cdot\underbrace{\left\{\ln\Gamma\!\left(\!\frac{l}{m}\!\right) -\,
\ln\Gamma\!\left(\!1-\frac{l}{m}\!\right) \!\right\}}_{2\ln\Gamma(l/m)+\ln\sin({\pi l/m})-\ln\pi} \\[9mm]
\displaystyle \qquad\qquad\qquad
= 2\!\!\!\!\!\!\sum_{l=1}^{\lfloor\!\frac{1}{2}(m-1)\!\rfloor} \!\!\!\!
\sin\frac{2\pi rl}{m}  \cdot\ln\Gamma\!\left(\!\frac{l}{m}\!\right) +
 \!\!\!\!\!\!\sum_{l=1}^{\lfloor\!\frac{1}{2}(m-1)\!\rfloor} \!\!\!\!
\sin\frac{2\pi rl}{m}  \cdot\ln\sin\frac{\pi l}{m}    \\[6mm]
\displaystyle \qquad\qquad\qquad\quad
- \ln\pi\cdot
\csc\frac{\pi r}{m}\cdot
\sin\!\left(\!\frac{\pi r}{m}\left\lfloor\!\frac{m+1}{2}\!\right\rfloor \right)\cdot
\sin\!\left(\!\frac{\pi r}{m}\left\lfloor\!\frac{m-1}{2}\!\right\rfloor \right)
\end{array}
\ee
because for natural $n$\vspace{-5pt}
\begin{eqnarray*}
\sum_{l=1}^{n} \sin(lx) = \csc
\frac{x}{2}\cdot\sin\frac{nx}{2}\cdot\sin \left[
\frac{x}{2}(n+1) \right]
\end{eqnarray*}
see e.g.~\cite[\no58, p.~12]{gunter_03_eng}. In like manner\vspace{-6pt}
\begin{eqnarray*}
\sum_{l=1}^{m-1} \cos \frac{2\pi rl}{m}
\cdot\ln\sin\frac{\pi
l}{m} = 2 \!\!\!\!\!\!\sum_{l=1}^{\lfloor\frac{1}{2}(m-1) \rfloor}\!\!\!\!\!
\cos\frac{2\pi rl}{m}\cdot\ln\sin \frac{\pi l}{m}
\end{eqnarray*}
By using the last three identities, Eq.~{\eqref{kx98eujms}} reduces to
{\eqref{hg784bdn}}.

The theorem may be also written by means of the Digamma function. In
fact, by recalling that
Gauss' Digamma theorem (\ref{dh38239djws}b) provides a connection
between the last
sum in {\eqref{kx98eujms}} and the $\Psi$-function,
formula {\eqref{kx98eujms}} may be also written in the following form:
\be\label{khjwe9jdns}
\begin{array}{lr}
\displaystyle
\gamma_1 \biggl(\frac{r}{m} \vphantom{\frac{1}{2}} \biggr)
 =& \displaystyle\,
\gamma_1 +\gamma^2 + \gamma\ln2\pi m + \ln2\pi\cdot\ln{m} +\frac{1}{2}\ln^2\!{m}
+ (\gamma+\ln2\pi m)\cdot\Psi\!\left(\vphantom{\frac{1}{2}}\!\frac{r}{m}\!\right) \\[6mm]
\displaystyle & \displaystyle\qquad
+\pi\sum_{l=1}^{m-1} \sin\frac{2\pi r l}{m} \cdot\ln\Gamma \biggl(\frac{l}{m} \biggr) 
+ \sum_{l=1}^{m-1} \cos\frac{2\pi rl}{m}\cdot\zeta''\!\left(\! 0,\,\frac{l}{m}\!\right) 
\end{array}
\ee

In some cases, it may be more advantageous to have the complete finite
Fourier series form. For this aim, it suffices
to take again {\eqref{kx98eujms}} and to recall that
\be\label{jk90d2dmekwe}
\sum_{l=1}^{m-1} l\cdot\sin\dfrac{2\pi rl}{m} \,=\,-\frac{m}{2}\ctg\dfrac{\pi r}{m}\,,
\qquad\qquad\qquad
r=1, 2, \ldots,m-1.
\ee
This yields the following expression
\be\label{kjkhnc823h}
\begin{array}{ll}
\displaystyle
\gamma_1 \biggl(\frac{r}{m} \vphantom{\frac{1}{2}} \biggr)
 =& \displaystyle\,
\gamma_1 - \gamma\ln2 m - \ln^2\!2 - \ln2\cdot\ln{\pi m} - \frac{1}{2}\ln^2\!{m} \\[6mm]
& \displaystyle\,
+\pi\sum_{l=1}^{m-1} \sin\frac{2\pi r l}{m} \cdot\left\{\ln\Gamma \biggl(\frac{l}{m} \biggr) 
+\frac{l(\gamma+\ln2\pi m)}{m} \right\} \\[6mm]
& \displaystyle\,
+ \sum_{l=1}^{m-1} \cos\frac{2\pi rl}{m}\cdot\left\{\zeta''\!\left(\! 0,\,\frac{l}{m}\!\right) 
+(\gamma+\ln2\pi m) \ln\sin\frac{\pi l}{m} \right\}
\end{array}
\ee
where $r=1, 2, 3,\ldots, m-1$, and $m$ is positive integer greater
than 1.
\qed
\end{pf}

Formulae~{\eqref{hg784bdn}} {\eqref{kx98eujms}},
{\eqref{khjwe9jdns}},
{\eqref{kjkhnc823h}} and {\eqref{asjd83a3}} permit
to readily obtain closed-form expressions for $\gamma_1(v)$ at any
rational $v$.
We, however, remark in passing that in some cases, these expressions
may be further simplified so that the resulting
formul\ae~may not contain at all $\zeta''(0,l/m)+\zeta
''(0,1-l/m)$, or contain only one combination (or transcendent) of
them. More detailed information related to these two special cases are
provided in  \ref{pi3ofj94cm3pimc3opij}.

\subsection{Summation formul\ae~with the first generalized
Stieltjes constant~at rational argument}\label{ij2093dm32ddf}

The derived theorem is very useful for many purposes, and in
particular, for the derivation of
summation formul\ae~involving the first generalized Stieltjes
constant~at rational argument.

\begin{teorema}
For the first generalized Stieltjes constant~at rational argument take
place following summation formul\ae%
\be\label{jw38endnslx03}
\specialnumber{a,b}
\begin{cases}
\displaystyle
\sum_{r=1}^{m-1} \gamma_1 \biggl(\!\frac{r}{m} \vphantom{\frac{1}{2}} \!\biggr)
\!\cdot\cos\dfrac{2\pi rk}{m} \,=\, -\gamma_1  + m(\gamma+\ln2\pi m)
\ln\!\left(\!2\sin\frac{\,k\pi\,}{m}\!\right)   \\[5mm]
\displaystyle\qquad\qquad\qquad\qquad\qquad\qquad\qquad\qquad\qquad\qquad
+\frac{m}{2}
\left\{\zeta''\!\left(\! 0,\,\frac{k}{m}\!\right)+ \, \zeta''\!\left(\! 0,\,1-\frac{k}{m}\!\right) \! \right\}  \\[7mm]
\displaystyle
\sum_{r=1}^{m-1} \gamma_1\biggl(\!\frac{r}{m} \vphantom{\frac{1}{2}} \!\biggr)
\!\cdot\sin\dfrac{2\pi rk}{m} \,=\,\frac{\pi}{2} (\gamma+\ln2\pi m)(2k-m) 
- \frac{\pi m}{2} \left\{\ln\pi -\ln\sin\frac{k\pi}{m} \right\}   \\[5mm]
\displaystyle\qquad\qquad\qquad\qquad\qquad\qquad\qquad\qquad\qquad\qquad\qquad\qquad\qquad
+ m\pi\ln\Gamma \biggl(\frac{k}{m} \biggr)  
\end{cases}
\ee
for $k=1, 2, 3, \ldots, m-1$, where $m$ is natural greater than
1.\footnote{One of these formul\ae~also appears in an unpublished
work sent to the author by Donal Connon.}
\end{teorema}

\begin{pf}
Formula {\eqref{kjkhnc823h}} represents the finite Fourier
series of the type {\eqref{lkwcjh3098j}}. Comparing
{\eqref{kjkhnc823h}} to {\eqref{lkwcjh3098j}}, we
immediately identify
%
\begin{eqnarray}
\everymath{\displaystyle}
\label{lkjx209mxechgcbn} 
\begin{cases}
a_m(0) = \gamma_1 - \gamma\ln2 m - \ln^2 2 - \ln2\cdot\ln{\pi m}
- \frac{1}{2}\ln^2 {m} ,  \vspace{10pt}\cr
a_m(l) = \zeta''  \left( 0, \frac{l}{m}  \right)
+(\gamma+\ln2\pi m) \ln\sin\frac{\pi l}{m} ,
\quad  l=1, 2, 3,\ldots, m-1 \vspace{10pt} \cr
b_m(l) = \pi \left\{\ln\Gamma \left(\frac{l}{m}  \right) +\frac
{l(\gamma+\ln2\pi m)}{m}  \right\} ,
\quad  l=1, 2, 3,\ldots, m-1
\end{cases}
\end{eqnarray}
Thus, in virtue of {\eqref{lkjcx20djndmednsd}}, for any $k=1, 2,
3,\ldots, m-1$,
\be\notag
\begin{array}{ll}
\displaystyle
\sum_{r=1}^{m-1} \gamma_1 \biggl(\!\frac{r}{m} \vphantom{\frac{1}{2}} \!\biggr)
\!\cdot\cos\dfrac{2\pi rk}{m} \,=\,
-\gamma_1 + \gamma\ln2 m + \ln^2\!2 + \ln2\cdot\ln{\pi m} + \frac{1}{2}\ln^2\!{m}
- \!\!\!\underbrace{\sum_{l=1}^{m-1}\zeta''\!\left(\! 0,\,\frac{l}{m}\!\right)}_
{-\frac{1}{2}\ln^2\!{m}-\ln{m}\cdot\ln2\pi} \\[12mm]
\displaystyle\qquad
-(\gamma+\ln2\pi m) \!\underbrace{\sum_{l=1}^{m-1} 
\ln\sin\frac{\pi l}{m}}_{(1-m)\ln2+\ln{m}}
+ \, \frac{m(\gamma+\ln2\pi m)}{2}\Biggl[\ln\sin\frac{\pi k}{m} + \ln\sin\frac{\pi (m-k)}{m}\Biggr] \\[12mm]
\displaystyle\qquad
+ \frac{m}{2}
\left\{\zeta''\!\left(\! 0,\,\frac{k}{m}\!\right) + \, \zeta''\!\left(\! 0,\,1-\frac{k}{m}\!\right) \! \right\} =
-\gamma_1  + m(\gamma+\ln2\pi m)
\cdot\ln\!\left(\!2\sin\frac{\,k\pi\,}{m}\!\right) \\[5mm]
\displaystyle \qquad\qquad\qquad\qquad\qquad\qquad\qquad\qquad\qquad\qquad
+\frac{m}{2}
\left\{\zeta''\!\left(\! 0,\,\frac{k}{m}\!\right) + \, \zeta''\!\left(\! 0,\,1-\frac{k}{m}\!\right) \! \right\} 
\end{array}
\ee
where we respectively used the multiplication theorem for the Hurwitz
$\zeta$-function
\be\label{2piod29rjddfg}
\sum_{l=1}^{m-1}\zeta''\!\left(\! 0,\,\frac{l}{m}\!\right) =
\left.\frac{d^2}{ds^2}[(n^s-1)\zeta(s)]\right|_{s=0}\!\!=
-\frac{1}{2}\ln^2\!{m}-\ln{m}\cdot\ln2\pi
\ee
see e.g.~\cite[p.~101]{iaroslav_06},
and the well-known formula from elementary mathematical analysis
%
\begin{eqnarray*}
\label{jhc23984nd} \prod_{l=1}^{m-1} \sin
\frac{\pi l}{m} = \frac{m}{2^{m-1}}
\end{eqnarray*}
which is, by the way, due to Euler \cite[tomus I, art.~240,
p.~204]{euler_04}, \cite[tome II, art.~99, p.~445]{legendre_02}
or \cite[vol.~I, p.~752, \no6.1.2-2]{prudnikov_en}.
Analogously, by {\eqref{hiuylckja2kdd}}, we deduce
\be
\begin{array}{ll}
\displaystyle
\sum_{r=1}^{m-1} \gamma_1 \biggl(\!\frac{r}{m} \vphantom{\frac{1}{2}} \!\biggr)
\!\cdot\sin\dfrac{2\pi rk}{m} \,= \, \frac{\pi m}{2}\Biggl\{
\ln\Gamma\!\left(\!\dfrac{k}{m}\!\right) - \ln\Gamma\!\left(\!1-\dfrac{k}{m}\!\right) 
+\frac{\gamma+\ln2\pi m}{m}\big[k-(m-k) \big] \Biggr\} \\[8mm]
\displaystyle\qquad\qquad\qquad\quad
=\,\frac{\pi}{2} (\gamma+\ln2\pi m)(2k-m) 
- \frac{\pi m}{2} \left\{\ln\pi -\ln\sin\frac{\pi k}{m} \right\}  + m\pi\ln\Gamma \biggl(\frac{k}{m} \biggr) 
\end{array}
\ee
which holds for $k=1, 2, 3,\ldots, m-1$.\qed
\end{pf}

\begin{teorema}
Parseval's theorem for the first generalized Stieltjes constant~at
rational argument has the following form
\be\label{jhc983jc3ncgbf}
\begin{array}{l}
\displaystyle
\sum_{r=1}^{m-1} \gamma^2_1 \biggl(\!\frac{r}{m} \vphantom{\frac{1}{2}} \!\biggr)\!=\, 
(m-1)\gamma^2_1 - m\gamma_1(2\gamma+\ln{m})\ln m  +
\frac{m}{4}\sum_{l=1}^{m-1} \!\left\{\zeta''\!\left(\! 0,\,\frac{l}{m}\!\right) 
+ \, \zeta''\!\left(\! 0,\,1-\frac{l}{m}\!\right) \! \right\}^{\!2}   \\[8mm]
\displaystyle\qquad
+ m(\gamma+\ln2\pi m) \!\sum_{l=1}^{m-1} \! \left\{\zeta''\!\left(\! 0,\,\frac{l}{m}\!\right) 
+ \, \zeta''\!\left(\! 0,\,1-\frac{l}{m}\!\right) \! \right\} \!\cdot\ln\sin\frac{\pi l}{m}
+m\pi^2\! \sum_{l=1}^{m-1} \ln^2\!\Gamma\!\left(\! \frac{l}{m}\!\right)\\[8mm]
\displaystyle \qquad
+2\pi^2 (\gamma+\ln2\pi m)\! \sum_{l=1}^{m-1} l\!\cdot\!\ln\Gamma\!\left(\! \frac{l}{m}\!\right)
+ \frac{m}{4}\big[4(\gamma+\ln2\pi m)^2-\pi^2\big]\! 
\sum_{l=1}^{m-1} \ln^2\!\sin\frac{\pi l}{m} + C_m
\end{array}
\ee
where, for convenience in writing, by $C_m$ we designated an elementary
function depending on $m$ and containing Euler's constant $\gamma$
\be\notag
\begin{array}{ll}
C_m \equiv - m(m-1)\ln^4\!2 - m(m-1)(2\ln{m}+2\gamma+3\ln\pi)\ln^3\!2 - m(m-2)\ln^2\! m \cdot\ln^2\! 2 \\[3mm]
- 2m\big[2(m-1)\ln\pi+\gamma (m-2)\big] \ln{m}\cdot\ln^2\! 2 
-m(m-1)\bigl[3\ln^2\!\pi+4\gamma \ln\pi + \gamma^2+\frac{5}{12}\pi^2+\frac{1}{6m}\pi^2\bigr]\times \\[3mm]
\times\ln^2\!2 
- m\left[(m-{\frac{5}{2}})\ln{\pi}-3\gamma\right]\ln{m}^2\cdot\ln2
+2m\big[(1-m)\ln^2\!\pi - (m-{\frac{5}{2}})\gamma\ln\pi\big]\ln{m}\cdot\ln2   \\[3mm]
+ \frac{1}{12}\Bigl[\left((6 \pi^2+24 \gamma^2) m+ 4 \pi^2 (1-m^2) \right)\ln{m} - 
4(m-1)\Bigl(3 m \ln^3\!\pi+6 m \gamma\ln^2\!\pi  +\gamma \pi^2 (m+1)  \\[3mm]
+(({\textstyle\frac{13}{4}}\pi^2 +3 \gamma^2) m+\pi^2) \ln\pi\Bigr)\Bigr] \ln2
+ \frac{1}{4} m\ln^4\!{m}  + m(\gamma+{\frac{1}{2}}\ln\pi)\ln^3\!{m} 
+ \frac{1}{12} \Bigl[6 m \ln^2\pi  \\[3mm]
+18 \gamma m \ln\pi+ \pi^2 m^2+(12 \gamma^2+3 \pi^2) m+2 \pi^2\Bigr]\ln^2\!{m} +m\ln^3\!m\cdot\ln2 +
\frac{1}{12}\Bigl[12 m \gamma \ln^2\!\pi \\[3mm]
+\left((12 \gamma^2+9 \pi^2)m + 4 \pi^2 (1-m^2)\right) \ln{\pi}+2 \pi^2 (2+m^2) \gamma\Bigr] \ln{m}\\[3mm]
-\frac{1}{12} (m-1) \bigl[2\pi^2 (4m+1) \ln^2\!{\pi}+4 \gamma \pi^2 (m+1) \ln{\pi} - \pi^2 \gamma^2 (m-2)\bigr]\\[3mm]
- \frac{1}{4}m\big[4(\gamma+\ln2\pi m)^2-\pi^2\big]\!\cdot\!\big[(1-m)\ln2+\ln m \big]\ln\pi
+ m(\gamma+\ln2\pi m)(\frac{1}{2}\ln m +\ln2\pi)\ln\pi\cdot\ln{m}
\end{array}
\ee
and where $m$ is natural greater than 1.
\end{teorema}

\begin{pf} Inserting Fourier series coefficients
{\eqref{lkjx209mxechgcbn}} into (\ref{hhxqwuybs}b)
and proceeding analogously to
\eqref{kjch328hdcnls}--\eqref{7826gdxhkwx39e4},
yields, after several pages of careful calculations and
simplifications, the above result.
The unique formula that should be used in addition to those employed in
derivations \eqref{kjch328hdcnls}--\eqref{7826gdxhkwx39e4}
 is
%
\begin{eqnarray}
\label{h203hdl2sdmn} \sum_{l=1}^{m-1} l\cdot\ln
\sin\frac{\pi l}{m} = \frac{m}{2}\sum_{l=1}^{m-1}
\ln\sin\frac{\pi l}{m} = \frac{m[(1-m)\ln2 + \ln m]}{2}
\end{eqnarray}
Also, the fact that the reflected sum
$\zeta''(0,l/m)+\zeta''(0,1-l/m)$, as well as the function $\ln\sin
(\pi
l/m)$, are both invariant with respect to a change of summation's index
$l\to m-l$ greatly helps when simplifying formula
{\eqref{jhc983jc3ncgbf}}. \qed
\end{pf}

Analogously, {\eqref{kjkhnc823h}} allows us to obtain a number of other
interesting summation formul\ae~for the first generalized
Stieltjes constant~at rational argument. For instance, with the help of \eqref{olkefj938hjcwdc}, 
we easily deduce these results
\be\label{jkcn3928echd}
\specialnumber{a,b}
\begin{cases}
\displaystyle
\sum_{r=0}^{m-1} \cos\dfrac{(2r+1)\pi k}{m}\cdot\gamma_1 \biggl(\!\frac{2r+1}{2m} \vphantom{\frac{1}{2}} \!\biggr)
\,=\,  m(\gamma+\ln4\pi m)\ln\tg\frac{\,\pi k\,}{2m}  +  \\[6mm]
\displaystyle \qquad\;\;
+\frac{m}{2}
\left\{\zeta''\!\left(\! 0,\,\frac{k}{2m}\!\right) + \, \zeta''\!\left(\! 0,\,1-\frac{k}{2m}\!\right) \! \right\}
-\frac{m}{2} \left\{\zeta''\!\left(\! 0,\,\frac{1}{2}+\frac{k}{2m}\!\right) + \, \zeta''\!\left(\! 0,\,\frac{1}{2}-\frac{k}{2m}\!\right) \! \right\} \\[8mm]
\displaystyle
\sum_{r=0}^{m-1} \sin\dfrac{(2r+1)\pi k}{m}\cdot\gamma_1 \biggl(\!\frac{2r+1}{2m} \vphantom{\frac{1}{2}} \!\biggr)
\,=\, m\pi\left\{\ln\Gamma \biggl(\frac{k}{2m} \biggr)  
+ \ln\Gamma \biggl(\frac{1}{2}-\frac{k}{2m} \biggr)  +\frac{1}{2}\ln\sin\frac{\,\pi k\,}{m} \right\}\\[6mm]
\displaystyle \qquad\qquad\qquad\qquad\qquad\qquad\qquad\qquad\qquad\qquad\qquad\qquad\;\;
-\frac{\pi m}{2}\big(3\ln2\pi+\ln m+\gamma \big)
\end{cases}
\ee
for $k=1, 2, 3, \ldots, m-1$, where $m$ is natural greater than 1. By a
similar line of argument, we also deduce \label{o2108wx0j0wxm}
%
\begin{eqnarray*}
\everymath{\displaystyle}
\left\{
\begin{array}{l}
\!\sum_{r=1}^{m-1} \cos\frac{(2k+1)\pi r}{m}\cdot
\gamma_1 \left( \frac{r}{m}   \right) =\\[12pt]
\qquad\qquad\quad
=-\pi\sum_{r=1}^{m-1} \frac{\sin\frac{2\pi r}{m}}{ \cos\frac
{2\pi r}{m} -\cos\frac{(2k+1)\pi}{m} }
\left\{\ln\Gamma \left(\frac{r}{m}  \right) +\frac
{r(\gamma+\ln2\pi m)}{m}  \right\}\\[18pt]
\!\sum_{r=1}^{m-1} \sin\frac{(2k+1)\pi r}{m}\cdot
\gamma_1  \left( \frac{r}{m}   \right)
 =
 \left[\gamma_1 - \gamma\ln2 m - \ln^2 2 - \ln2\cdot\ln{\pi m} -
\frac{1}{2}\ln^2 {m} \right]{\times} \\[13pt]
\quad{}{\times} \operatorname{ctg}\frac{(2k+1)\pi
}{2m}+(\gamma+\ln2\pi m) \sin\frac{(2k+1)\pi}{m}\cdot\!\sum_{r=1}^{m-1} \frac{1}
{ \cos\frac{2\pi r}{m} -\cos\frac{(2k+1)\pi}{m} }\cdot\ln\sin
\frac{\pi r}{m}
\\[13pt]
\quad {}
+\frac{1}{2} \sin\frac{(2k+1)\pi}{m}\cdot\sum_{r=1}^{m-1} \frac{1}
{ \cos\frac{2\pi r}{m} -\cos\frac{(2k+1)\pi}{m} } \cdot \left\{
\zeta''  \left( 0, \frac{r}{m}  \right)
+ \zeta''  \left( 0, 1-\frac{r}{m}  \right)
\right\}\end{array}\right.
\end{eqnarray*}
which are valid for any $k\in\mathbbm{Z}$. By the way, two
particular cases of (\ref{jw38endnslx03}a)
and (\ref{jkcn3928echd}b) may represent some special
interest. Thus putting in the former $k=m/2$ when $m$ is even yields
%
\begin{eqnarray}
\label{idh20djd} \sum_{r=1}^{2m-1}
(-1)^r \cdot\gamma_1 \left( \frac{r}{2m}
\right) = -\gamma_1+m(2\gamma+\ln2+2\ln m)\ln2
\end{eqnarray}
However, the same relationship may be also derived from the
multiplication theorem for the first Stieltjes constant\footnote{This
is a particular case of
the multiplication theorem for the first generalized Stieltjes
constant. More general case of this theorem
and equivalent theorems for higher Stieltjes constants~were derived in
exercise \no64 \cite[p.~101, Eqs.~(62)--(63)]{iaroslav_06}.
Some particular cases of these theorems
appear also in \cite[Eqs.~(3.28), (3.54)]{coffey_01}; Eq.~(3.54)
contains, unfortunately, an error (see
footnote 42 \cite[p.~101]{iaroslav_06}).}
%
\begin{eqnarray}
\label{kjhc29hdn} \sum_{r=1}^{m-1}
\gamma_1 \left( \frac{r}{ m } \right) = (m-1)
\gamma_1 - m\gamma\ln{m} - \frac{m}{2} \ln^2 m
\end{eqnarray}
Putting $2m$ instead of $m$, and then, treating separately odd and even terms,
we have
%
\begin{eqnarray}
\label{kjhc29hdn2} \sum_{r=0}^{m-1}
\gamma_1 \left( \frac{2r+1}{ 2m } \right) = m
\left\{\gamma_1 -\gamma\ln4m - \frac{1}{2}
\ln^2 m - \ln^2 2-2\ln2\cdot\ln{m} \right\}
\end{eqnarray}
Subtracting from the above sum even terms $\gamma_1(2r/2m)$ for $r=1,
2, \ldots, m-1$, immediately yields {\eqref{idh20djd}}.
In other words, {\eqref{idh20djd}} may be also regarded as a direct
consequence of the multiplication theorem for the first Stieltjes constant.
In contrast, the particular case of Eq.~ (\ref{jkcn3928echd}b)
corresponding to $k=m/2$ when $m$ is even
%
\begin{eqnarray}
\label{kjd943k3cke} \sum_{r=0}^{2m-1}
(-1)^r \cdot\gamma_1 \left( \frac{2r+1}{4m}
\right) = m \left\{4\pi\ln\Gamma \left( \frac{1}{4}
\right) - \pi (4\ln2+3\ln\pi+\ln m+\gamma ) \right\}
\end{eqnarray}
appears to be more interesting and cannot be derived solely from {\eqref
{kjhc29hdn}}.
Moreover, we can also combine
{\eqref{kjd943k3cke}} with
{\eqref{kjhc29hdn2}} by putting in the later
$2m$ instead of $m$. Adding and subtracting them respectively yields:
\be\label{984ch378hbd}
\specialnumber{a,b}
\begin{array}{ll}
\displaystyle\sum_{r=0}^{m-1} \gamma_1\biggl(\!\frac{4r+1}{4m} \vphantom{\frac{1}{2}} \!\biggr)= & \displaystyle 
\frac{m}{2}\left\{2\gamma_1 - \gamma\big(\pi+6\ln2+2\ln{m}\big)+4\pi\ln\Gamma \biggl(\frac{1}{4} \biggr) 
-4\pi\ln2   \right.\\[4mm]
&\displaystyle\!\!\!\!
\left. \vphantom{\ln\Gamma \biggl(\frac{1}{4} \biggr) }
-3\pi\ln\pi-\pi\ln{m}-7\ln^2\!2-6\ln2\cdot\ln{m}-\ln^2\!{m}\right\}\\[5mm]
\displaystyle\sum_{r=0}^{m-1} \gamma_1\biggl(\!\frac{4r+3}{4m} \vphantom{\frac{1}{2}} \!\biggr)= & \displaystyle 
\frac{m}{2}\left\{2\gamma_1 + \gamma\big(\pi-6\ln2-2\ln{m}\big)-4\pi\ln\Gamma \biggl(\frac{1}{4} \biggr) 
+4\pi\ln2   \right.\\[4mm]
&\displaystyle\!\!\!\!
\left. \vphantom{\ln\Gamma \biggl(\frac{1}{4} \biggr) }
+3\pi\ln\pi+\pi\ln{m}-7\ln^2\!2-6\ln2\cdot\ln{m}-\ln^2\!{m}\right\}
\end{array}
\ee
 From these equations it follows, \emph{inter alia}, that
sums $\gamma_1(1/8)+\gamma_1(5/8)$ and $\gamma_1(1/12)+\gamma_1(5/12)$
may be expressed in terms of $\Gamma(1/4)$, $\gamma_1$, $\gamma$ and
elementary functions.\footnote{For the value
of $\gamma_1(3/4)$, see \cite[p.~100]{iaroslav_06}.}
Besides, the role of $\ln\Gamma(1/4)$ in three latter identities
seems quite intriguing because the logarithm of the $\Gamma$-function
possesses
very similar properties
\be\notag
\begin{array}{ll}
\displaystyle 
\sum_{r=0}^{2m-1} \! (-1)^r\!\cdot \ln\Gamma\!\left(\!\frac{\,2r+1\,}{4m}\!\right) 
=\, 2\ln\Gamma\!\left(\!\frac{\,1\,}{4}\!\right) - \frac{1}{2}(\ln2+2\ln\pi-\ln{m}) \\[6mm]
\displaystyle 
\sum_{r=0}^{m-1} \ln\Gamma\!\left(\!\frac{\,4r+1\,}{4m}\!\right) \,=\,\ln\Gamma\!\left(\!\frac{\,1\,}{4}\!\right)
+\frac{1}{2}(m-1)\ln{2\pi} +\,\frac{1}{4}\ln{m}\,\\[6mm]
\displaystyle 
\sum_{r=0}^{m-1} \ln\Gamma\!\left(\!\frac{\,4r+3\,}{4m}\!\right) \,=\,-\ln\Gamma\!\left(\!\frac{\,1\,}{4}\!\right)
+\frac{m}{2}\ln{2\pi} +\,\frac{1}{4}\ln{\frac{\pi^2}{m}}\,
\end{array}
\ee
Particular cases of (\ref{jw38endnslx03}b) corresponding to $k=m/3$
and $k=m/6$ are also interesting.
Put in (\ref{jw38endnslx03}b) $3m$ instead of $m$, and then, set
$k=m$. This yields:
\be\label{9032udn2d}
\begin{array}{ll}
\displaystyle 
\gamma_1\biggl(\!\frac{1}{3m} \vphantom{\frac{1}{2}} \!\biggr) 
-\gamma_1\biggl(\!\frac{2}{3m} \vphantom{\frac{1}{2}} \!\biggr) 
+\gamma_1\biggl(\!\frac{4}{3m} \vphantom{\frac{1}{2}} \!\biggr) 
-\gamma_1\biggl(\!\frac{5}{3m} \vphantom{\frac{1}{2}} \!\biggr)  + \ldots 
+\gamma_1\biggl(\!\frac{3m-2}{3m} \vphantom{\frac{1}{2}} \!\biggr) 
-\gamma_1\biggl(\!\frac{3m-1}{3m} \vphantom{\frac{1}{2}} \!\biggr)\\[6mm]
\displaystyle \qquad\qquad\qquad\qquad\qquad
\qquad\,=\,\frac{\pi m}{\sqrt{3}}\left\{6\ln\Gamma \biggl(\frac{1}{3} \biggr)
-\gamma-4\ln2\pi +\frac{1}{2}\ln3 -\ln m \!\right\} 
\end{array}
\ee
But the multiplication theorem \eqref{kjhc29hdn} rewritten for $3m$ in
place of $m$ gives
\be\label{j093eiocmdn}
\begin{array}{ll}
\displaystyle 
\gamma_1\biggl(\!\frac{1}{3m} \vphantom{\frac{1}{2}} \!\biggr) 
+\gamma_1\biggl(\!\frac{2}{3m} \vphantom{\frac{1}{2}} \!\biggr) 
+\gamma_1\biggl(\!\frac{4}{3m} \vphantom{\frac{1}{2}} \!\biggr) 
+\gamma_1\biggl(\!\frac{5}{3m} \vphantom{\frac{1}{2}} \!\biggr)  + \ldots 
+\gamma_1\biggl(\!\frac{3m-2}{3m} \vphantom{\frac{1}{2}} \!\biggr) 
+\gamma_1\biggl(\!\frac{3m-1}{3m} \vphantom{\frac{1}{2}} \!\biggr)\\[6mm]
\displaystyle \qquad\quad
\qquad\,=\,2m\gamma_1 - m\gamma(2\ln{m}+3\ln3) - \frac{m}{2}\left(3\ln^2\!3+6\ln3\cdot\ln{m}+2\ln^2\!{m}\right)
\end{array}
\ee
and hence
\addtocounter{equation}{-1}
\be\label{kjxc23jr}
\specialnumber{b}
\begin{array}{ll}
\displaystyle\sum_{r=0}^{m-1} \gamma_1\biggl(\!\frac{3r+1}{3m} \vphantom{\frac{1}{2}} \!\biggr)=
\,m\left\{\gamma_1 -
\gamma \left(\frac{\pi}{2\sqrt{3}}+\ln{m}+\frac{3}{2}\ln3 \right)
+\pi\sqrt{3}\ln\Gamma \biggl(\frac{1}{3} \biggr)  \right.\\[6mm]
\displaystyle \qquad\quad
\left.
-\frac{\pi}{2\sqrt{3}}  \left(4\ln2\pi -\frac{1}{2}\ln3 + \ln m \right)
- \frac{1}{4}\left(3\ln^2\!3+6\ln3\cdot\ln{m}+2\ln^2\!{m}\right) \!\right\}
\end{array}
\ee
Consider now the particular case of (\ref{jw38endnslx03}a) corresponding to $k=m/6$. Recalling that  
$\,\ln\Gamma(1/6)=\frac{1}{2}\ln3 -\frac{1}{3}\ln2-\frac{1}{2}\ln\pi +2\ln\Gamma(1/3)\,$, we have
\be
\begin{array}{ll}
\displaystyle 
\gamma_1\biggl(\!\frac{1}{6m} \vphantom{\frac{1}{2}} \!\biggr) 
+\gamma_1\biggl(\!\frac{2}{6m} \vphantom{\frac{1}{2}} \!\biggr) 
-\gamma_1\biggl(\!\frac{4}{6m} \vphantom{\frac{1}{2}} \!\biggr) 
-\gamma_1\biggl(\!\frac{5}{6m} \vphantom{\frac{1}{2}} \!\biggr)  + \ldots 
-\gamma_1\biggl(\!\frac{3m-2}{6m} \vphantom{\frac{1}{2}} \!\biggr) 
-\gamma_1\biggl(\!\frac{3m-1}{6m} \vphantom{\frac{1}{2}} \!\biggr)\\[6mm]
\displaystyle \qquad\qquad\qquad
\qquad\,=\,\frac{2\pi m}{\sqrt{3}}\left\{12\ln\Gamma \biggl(\frac{1}{3} \biggr)
-2\gamma-9\ln2 +\ln3 -8\ln\pi -2\ln m \!\right\} 
\end{array}
\ee
By adding this to {\eqref{j093eiocmdn}} rewritten for $2m$
instead of $m$, and then, by subtracting (\ref{kjxc23jr}b) results in another summation relation
\be\label{kjce934jf2pkdmemf}
\begin{array}{l}
\displaystyle\sum_{r=0}^{m-1} \gamma_1\biggl(\!\frac{6r+1}{6m} \vphantom{\frac{1}{2}} \!\biggr)= 
\,m\left\{\gamma_1 -
\gamma \left(\!\frac{\sqrt{3}\pi}{2}+2\ln2 +\frac{3}{2}\ln3 +\ln{m}\!\right)
+3\pi\sqrt{3}\ln\Gamma \biggl(\frac{1}{3} \biggr)  \right.\\[7mm]
\displaystyle\qquad\qquad\qquad\qquad
-\frac{\pi}{2\sqrt{3}} \left(\!14\ln2 -\frac{3}{2}\ln3 + 12\ln\pi +3\ln m \right)
-\ln^2\!2 - \frac{3}{4}\ln^2\!3 \\[5mm]
\displaystyle \qquad\qquad\qquad\qquad \left.
-3\ln2\cdot\ln3 -2\ln2\cdot\ln m - \frac{3}{2}\ln3\cdot\ln{m} - \frac{1}{2}\ln^2\!m 
\vphantom{\frac{\sqrt{3}\pi}{2}}
\right\}
\end{array}
\ee
Previous relationships permit to derive several summation formul\ae~for $\gamma_1(\ldots/12m)$.
Put in {(\ref{984ch378hbd}a)} $3m$ instead of $m$ and then represent the
summation index $r$ as $3l+k$, where the new
summation index $l$ runs through $0$ to $m-1$ for each $k=0, 1, 2$.
Then {(\ref{984ch378hbd}a)} may be written as a sum of three terms last
of which equals {(\ref{984ch378hbd}b)}. Hence
\be\label{oxi2103ikdmdf}
\begin{array}{l}
\displaystyle\sum_{l=0}^{m-1} \gamma_1\biggl(\!\frac{12l+1}{12m} \vphantom{\frac{1}{2}} \!\biggr)
+ \sum_{l=0}^{m-1} \gamma_1\biggl(\!\frac{12l+5}{12m} \vphantom{\frac{1}{2}} \!\biggr)= \,
\,\frac{m}{2}\left\{4\gamma_1 - \gamma (4\pi+12\ln2 +6\ln3 +4\ln{m})
 \right.\\[6mm]
\displaystyle \qquad\qquad\qquad\qquad \;
+16\pi\ln\Gamma \biggl(\frac{1}{4} \biggr) 
-\pi(16\ln2+12\ln\pi+3\ln3+4\ln{m}) -14\ln^2\!2 \\[3mm]
\displaystyle \qquad\qquad\qquad\qquad  \left. \vphantom{\Gamma \biggl(\frac{1}{4} \biggr)  }
-3\ln^2\!3 -18\ln2\cdot\ln3-12\ln2\cdot\ln{m}-6\ln3\cdot\ln{m}-2\ln^2\!{m}
\right\}
\end{array}
\ee
Similarly, by separately treating odd and even terms in \eqref{kjce934jf2pkdmemf} written for $2m$ instead of $m$, we have
\be\notag
\begin{array}{l}
\displaystyle\sum_{l=0}^{m-1} \gamma_1\biggl(\!\frac{12l+1}{12m} \vphantom{\frac{1}{2}} \!\biggr)
+ \sum_{l=0}^{m-1} \gamma_1\biggl(\!\frac{12l+7}{12m} \vphantom{\frac{1}{2}} \!\biggr) = 
\,2m\left\{\gamma_1 -
\gamma \left(\!\frac{\sqrt{3}\pi}{2}+3\ln2 +\frac{3}{2}\ln3 +\ln{m}\!\right)  \right.\\[7mm]
\displaystyle \qquad\qquad\qquad\qquad 
+3\pi\sqrt{3}\ln\Gamma \biggl(\frac{1}{3} \biggr) -\frac{\pi}{2\sqrt{3}} \left(\!17\ln2 -\frac{3}{2}\ln3 + 12\ln\pi +3\ln m \right)\\[5mm]
\displaystyle\qquad\qquad\qquad\qquad  \left. 
-\frac{7}{2}\ln^2\!2 - \frac{3}{4}\ln^2\!3 
-\frac{9}{2}\ln2\cdot\ln3 -3\ln2\cdot\ln m - \frac{3}{2}\ln3\cdot\ln{m} - \frac{1}{2}\ln^2\!m 
\vphantom{\frac{\sqrt{3}\pi}{2}}
\right\}
\end{array}
\ee
From these relationships, it appears that the sum $\gamma
_1(1/12)+\gamma_1(5/12)$
may be expressed in terms of $\Gamma(1/4)$, $\gamma_1$, $\gamma$ and
elementary functions, while
$\gamma_1(1/12)+\gamma_1(7/12)$ contains $\Gamma(1/3)$ instead of
$\Gamma(1/4)$.\footnote{At the same time,
the difference $\gamma_1(1/12)- \gamma_1(7/12)$ may be written
as function of $\Gamma(1/4)$ and $ \zeta''(0,1/12)+\zeta''(0,11/12)
$. This follows from the argument developed here later.}
This is certainly correlated with the fact that $\Gamma(1/12)$ may be
written in terms of product $\Gamma(1/3)\cdot\Gamma(1/4)$,
see e.g.~\cite[p.~31]{campbell_01}. Many particular cases of
equations from pp.~{\pageref{jw38endnslx03}}--{\pageref{o2108wx0j0wxm}} will also imply $ \zeta''(0,p)+\zeta
''(0,1-p) $ at different rational $p$.
For instance, setting in (\ref{jw38endnslx03}a) $k=m/5$ and recalling
that $ \cos\frac{2}{5}\pi=\frac{1}{4} (\sqrt{5}-1 ) $,
$ \cos\frac{4}{5}\pi=-\frac{1}{4} (\sqrt{5}+1 ) $
and $ \sin\frac{1}{5}\pi=\frac{1}{4}\sqrt{10-2\sqrt{5}} $, as
well as using several times
the multiplication theorem {\eqref{kjhc29hdn}}, yields
\be\notag
\begin{array}{ll}
\displaystyle 
\sum_{l=0}^{m-1}\gamma_1\biggl(\frac{5l+1}{5m} \vphantom{\frac{1}{2}} \!\biggr)
+\sum_{l=0}^{m-1} \gamma_1\biggl(\!\frac{5l+4}{5m} \vphantom{\frac{1}{2}} \!\biggr)
 =\,\frac{m}{2\sqrt{5}}\left\{\! 4\gamma_1\sqrt{5}+ 10\!\left[\zeta''\!\left(\! 0,\,\frac{1}{5}\vphantom{\dfrac{7}{2}}\!\right) 
+ \zeta''\!\left(\! 0,\,\frac{4}{5}\vphantom{\dfrac{7}{2}}\!\right)\!\right] \right.\\[6mm]
\displaystyle \qquad\qquad\qquad\qquad\qquad\qquad\qquad\quad\;
-\gamma\!\left(4\sqrt{5} \ln{m} +10\ln\!\big(1+\sqrt{5}\big)-10\ln2+5\sqrt{5}\ln5\right) \\[3mm]
\displaystyle  \qquad\quad
- 10\big(\!\ln2+\ln5+\ln\pi+\ln m\big)\!\cdot\ln\!\big(1+\sqrt{5})    
+10\ln^2\!2 -\frac{10}{1+\sqrt{5}}\ln^2\!5 - 2\sqrt5\ln^2\!m\\[3mm]
\displaystyle\qquad\qquad\quad\left. 
 + 15\ln2\cdot\ln5
+10\ln2\cdot\ln\pi+5\ln5\cdot\ln\pi+10\ln2\cdot\ln{m}-5\sqrt5\ln5\cdot\ln m
\vphantom{\dfrac{7}{2}}\!\right\} 
\end{array}
\ee
Interestingly, the golden ratio $ \phi$ seems to play a certain role
in the above formula.

Let now consider the case $k=m/8$, where $k$ should be positive
integer.~Eq.~(\ref{jw38endnslx03}b),
employed together with both Eqs.~(\ref{984ch378hbd}a) and (\ref{984ch378hbd}b), provides
\be\notag
\begin{array}{ll}
\displaystyle 
\gamma_1\biggl(\!\frac{1}{8m} \vphantom{\frac{1}{2}} \!\biggr) 
+\gamma_1\biggl(\!\frac{3}{8m} \vphantom{\frac{1}{2}} \!\biggr) 
-\gamma_1\biggl(\!\frac{5}{8m} \vphantom{\frac{1}{2}} \!\biggr) 
-\gamma_1\biggl(\!\frac{7}{8m} \vphantom{\frac{1}{2}} \!\biggr)  + \ldots 
+\gamma_1\biggl(\!\frac{8m-7}{8m} \vphantom{\frac{1}{2}} \!\biggr) 
+\gamma_1\biggl(\!\frac{8m-5}{8m} \vphantom{\frac{1}{2}} \!\biggr)
\\[6mm]
\displaystyle \quad
-\gamma_1\biggl(\!\frac{8m-3}{8m} \vphantom{\frac{1}{2}} \!\biggr)
-\gamma_1\biggl(\!\frac{8m-1}{8m} \vphantom{\frac{1}{2}} \!\biggr)
=\,\pi m \sqrt{2}\biggl\{8\ln\Gamma \biggl(\frac{1}{8} \biggr) 
- 4\ln\Gamma \biggl(\frac{1}{4} \biggr) - 2\gamma -11\ln2 \\[5mm]
\displaystyle \qquad\qquad\qquad\qquad\qquad\qquad\qquad\qquad\qquad\qquad
- 4\ln\pi -2\ln m
- 2\ln\!\big(1+\sqrt{2}\big)\!\biggr\} 
\end{array}
\ee
At the same time, Eq.~(\ref{jw38endnslx03}a) for $k=m/8$, used
together with \eqref{idh20djd}, leads to
\be\notag
\begin{array}{ll}
\displaystyle
\gamma_1\biggl(\!\frac{1}{8m} \vphantom{\frac{1}{2}} \!\biggr) 
-\gamma_1\biggl(\!\frac{3}{8m} \vphantom{\frac{1}{2}} \!\biggr) 
-\gamma_1\biggl(\!\frac{5}{8m} \vphantom{\frac{1}{2}} \!\biggr) 
+\gamma_1\biggl(\!\frac{7}{8m} \vphantom{\frac{1}{2}} \!\biggr)  + \ldots 
+\gamma_1\biggl(\!\frac{8m-7}{8m} \vphantom{\frac{1}{2}} \!\biggr) 
-\gamma_1\biggl(\!\frac{8m-5}{8m} \vphantom{\frac{1}{2}} \!\biggr)\\[6mm]
\displaystyle \quad
-\gamma_1\biggl(\!\frac{8m-3}{8m} \vphantom{\frac{1}{2}} \!\biggr)
+\gamma_1\biggl(\!\frac{8m-1}{8m} \vphantom{\frac{1}{2}} \!\biggr)
=\,m \sqrt{2}\left\{4\!\left[\zeta''\!\left(\! 0,\,\frac{1}{8}\vphantom{\dfrac{7}{2}}\!\right) 
+ \zeta''\!\left(\! 0,\,\frac{7}{8}\vphantom{\dfrac{7}{2}}\right)\!\right] \right.\\[4mm]
\displaystyle \qquad\qquad\qquad\qquad
-4 \big(\gamma+4\ln2+\ln\pi+\ln{m}\big)\!\cdot\ln\!\big(1+\sqrt{2}\big)
\left.   + 7\ln^2\!2 +2\ln2\cdot\ln\pi 
\vphantom{\dfrac{7}{2}} \right\} 
\end{array}
\ee
Adding both equations together results in another summation relation
\be\notag
\begin{array}{l}
\displaystyle \notag
\sum_{r=0}^{m-1}\gamma_1\biggl(\!\frac{8r+1}{8m} \vphantom{\frac{1}{2}} \!\biggr) 
- \sum_{r=0}^{m-1}\gamma_1\biggl(\!\frac{8r+5}{8m} \vphantom{\frac{1}{2}} \!\biggr)   
=\,\frac{m}{\sqrt{2}}  \left\{4\!\left[\zeta''\!\left(\! 0,\,\frac{1}{8}\vphantom{\dfrac{7}{2}}\!\right) 
+ \zeta''\!\left(\! 0,\,\frac{7}{8}\vphantom{\dfrac{7}{2}}\right)\!\right] 
+8\pi\ln\Gamma \biggl(\frac{1}{8} \biggr) \right.\\[6mm]
\displaystyle \quad
- 4\pi\ln\Gamma \biggl(\frac{1}{4} \biggr)
-2\gamma\Big[\pi+2\ln\!\big(1+\sqrt{2}\big)\!\Big]
-2\big(\pi+8\ln2+2\ln\pi+2\ln{m}\big)\!\cdot\ln\!\big(1+\sqrt{2}\big)   +
\end{array}
\ee

\be\notag
\begin{array}{l}
\displaystyle \notag
\displaystyle \qquad\qquad\qquad\qquad\qquad\qquad\qquad
\left. 
+ 7\ln^2\!2 +2\ln2\cdot\ln\pi -\pi(11\ln2 +4\ln\pi +2\ln m)
\vphantom{\dfrac{7}{2}} \right\} 
\end{array}
\ee
Analogous relation with ``$+$'' instead of ``$-$'' in the left part has
much more simple form and follows directly from (\ref{984ch378hbd}a) rewritten for $2m$ in place of $m$
\be\notag
\begin{array}{l}
\displaystyle \notag
\sum_{r=0}^{m-1}\gamma_1\biggl(\!\frac{8r+1}{8m} \vphantom{\frac{1}{2}} \!\biggr) 
+ \sum_{r=0}^{m-1}\gamma_1 \biggl(\!\frac{8r+5}{8m} \vphantom{\frac{1}{2}} \!\biggr)  
= \, m\left\{2\gamma_1 - \gamma\big(\pi+8\ln2+2\ln{m}\big)+4\pi\ln\Gamma \biggl(\frac{1}{4} \biggr)  \right.\\[4mm]
\displaystyle\qquad\qquad\qquad\qquad\qquad\;\;
\left. \vphantom{\ln\Gamma \biggl(\frac{1}{4} \biggr) }
-5\pi\ln2  
-3\pi\ln\pi-\pi\ln{m}-14\ln^2\!2-8\ln2\cdot\ln{m}-\ln^2\!{m}\right\}
\end{array}
\ee
Similarly, one can obtain equations for $\sum [\gamma_1
(\frac{8r+3}{8m} ) \pm\gamma_1  (\frac{8r+7}{8m} ) ]$.

The above summation formul\ae~are not only interesting in themselves, but
also may be useful for the closed-form determination of certain first
Stieltjes constants~(expressions in \ref{pi3ofj94cm3pimc3opij}
are obtained precisely by means of such formul\ae). Besides,
summation formul\ae~akin to {\eqref{kjhc29hdn2}},
(\ref{984ch378hbd}a), (\ref{kjxc23jr}b),
{\eqref{kjce934jf2pkdmemf}} may be often more easily obtained
by the direct summation of {\eqref{kx98eujms}}. For the
derivation of such a formula, we, first, write in
{\eqref{kx98eujms}} $mn$ for $m$ and $rn+k$ for $r$, where
$n\in\mathbbm{N}$ and $k=1, 2, \ldots, n-1$. Then, we remark that for
$l=1, 2, 3, \ldots, mn-1$, we have
\be\notag
\begin{array}{ll}
\displaystyle 
\sum_{r=0}^{m-1} \cos\frac{2\pi l(nr+k)}{nm}
=\, m\cos\frac{2\pi lk}{nm}\cdot\Big\{\delta_{l,m}+\delta_{l,2m}+\delta_{l,3m}+\ldots+\delta_{l,(n-1)m} \Big\}\\[6mm]
\displaystyle 
\sum_{r=0}^{m-1} \sin\frac{2\pi l(nr+k)}{nm}
=\, m\sin\frac{2\pi lk}{nm}\cdot\Big\{\delta_{l,m}+\delta_{l,2m}+\delta_{l,3m}+\ldots+\delta_{l,(n-1)m} \Big\}\\[6mm]
\displaystyle 
\sum_{r=0}^{m-1} \ctg\frac{\pi(nr+k)}{nm}
=\, m\ctg\frac{\pi k}{n}
\end{array}
\ee
see e.g.~\cite[p.~8, \no33]{gunter_02_eng}, whence
\be\notag
\begin{array}{ll}
\displaystyle 
\sum_{r=0}^{m-1}\!\gamma_1   \biggl(\frac{nr+k}{nm} \vphantom{\frac{1}{2}} \biggr) =
m\!\left(\!\gamma_1-\gamma\ln2mn  - \ln^2\!2 - \ln2\cdot\ln{\pi mn} - \frac{1}{2}\ln^2\!{mn}\!\right) \\[5mm]
\displaystyle\qquad\qquad\qquad
+m\!\sum_{\lambda=1}^{n-1} \cos\frac{2\pi\lambda k}{n}\cdot\zeta''\!\left(\! 0,\,\frac{\lambda}{n}\!\right) 
+m\pi\!\sum_{\lambda=1}^{n-1} \sin\frac{2\pi\lambda k}{n} \cdot\ln\Gamma \biggl(\frac{\lambda}{n} \biggr) \\[5mm]
\displaystyle\qquad\qquad\qquad
- \frac{m\pi}{2}(\gamma+\ln2\pi mn)\ctg\frac{\pi k}{n}
+ m(\gamma+\ln2\pi mn)\!\sum_{\lambda=1}^{n-1} \cos\frac{2\pi\lambda k}{n} \cdot\ln\sin\frac{\pi \lambda}{n} 
\end{array}
\ee
Comparing the right-hand side of this equation with the parent
equation {\eqref{kx98eujms}} finally yields
\be
\frac{1}{m}\!\sum_{r=0}^{m-1}\!\gamma_1   \biggl(\frac{nr+k}{nm} \vphantom{\frac{1}{2}} \biggr) 
= \,\gamma_1   \biggl(\frac{k}{n} \vphantom{\frac{1}{2}} \biggr) 
-\left\{\!\gamma+\ln{2n}+\frac{1}{2}\ln{m}  + \frac{\pi}{2}\ctg\frac{\pi k}{n} 
- \sum_{\lambda=1}^{n-1} \cos\frac{2\pi\lambda k}{n} \cdot\ln\sin\frac{\pi \lambda}{n} \!\right\}\ln{m}
\ee
This relationship represents a special variant of the generalized
multiplication theorem for the first generalized Stieltjes
constant.\footnote{This variant may be also obtained from
\cite[Eq.~(6.6)]{connon_02} or \cite[p.~101, Eq.~(63)]{iaroslav_06} by
making use of Gauss' Digamma theorem
{\eqref{dh38239djws}}.}

Another summation formula with the first generalized Stieltjes constants
may be obtained by using respectively {\eqref{jk90d2dmekwe}},
(\ref{jw38endnslx03}b), \eqref{jk0909jjds2eja} and \eqref{h203hdl2sdmn}
\be
\begin{array}{l}
\displaystyle
\sum_{r=1}^{m-1} \ctg\frac{\pi r}{m} \cdot\gamma_1 \biggl(\!\frac{r}{m} \vphantom{\frac{1}{2}} \!\biggr) =\, \displaystyle
\frac{\pi }{6} \Big\{\!(1-m)(m-2)\gamma + 2(m^2-1)\ln2\pi - (m^2+2)\ln{m}\Big\} \\[4mm]
\displaystyle \qquad\qquad\qquad\qquad\qquad
-2\pi\!\sum_{l=1}^{m-1} l\!\cdot\!\ln\Gamma\!\left(\! \frac{l}{m}\!\right)
\end{array}
\ee

The normalized first-order moment of the first generalized Stieltjes
constant~may be derived from {\eqref{kjkhnc823h}}
by making use of {\eqref{jk0909jjds2eja}}, {\eqref{jk90d2dmekwe}}, \eqref
{hjg2873gdb},
{\eqref{jhc23984nd}}, as well as {\eqref{2piod29rjddfg}}.
This yields
\be
\begin{array}{ll}
\displaystyle
\sum_{r=1}^{m-1} \frac{r}{m} \cdot\gamma_1 \biggl(\!\frac{r}{m} \vphantom{\frac{1}{2}} \!\biggr) =\, &  \displaystyle
\frac{1}{2}\left\{\!(m-1)\gamma_1 - m\gamma\ln{m} - \frac{m}{2}\ln^2\!{m}\! \vphantom{\frac{73}{2}}  \right\}
-\frac{\pi}{2m}(\gamma+\ln2\pi m) \!\sum_{l=1}^{m-1} l\!\cdot\! \ctg\frac{\pi l}{m} \\[5mm]
&\displaystyle
-\frac{\pi}{2} \!\sum_{l=1}^{m-1} \ctg\frac{\pi l}{m} \cdot\ln\Gamma\biggl(\!\frac{l}{m} \vphantom{\frac{1}{2}} \!\biggr) 
\end{array}
\ee

More complicated summation relations may be obtained if considering
other functions.
For example, the summation formula with the Digamma function reads
\be
\begin{array}{ll}
&\displaystyle
\sum_{r=1}^{m-1} \displaystyle
\Psi \biggl(\!\frac{r}{m} \vphantom{\frac{1}{2}} \!\biggr)\!\cdot  
\gamma_1 \biggl(\!\frac{r}{m} \vphantom{\frac{1}{2}} \!\biggr) = 
\,\big[\gamma (1-m)-m\ln{m}\big]\gamma_1 + m\gamma^2\ln m 
+ \left\{\!\frac{(m-1)(m-2)\pi^2}{12}  \right.   \\[5mm]
&\displaystyle\qquad
\left. -m(m-1)\ln^2\!2  
+2m\ln2\cdot\ln{m}+ \frac{3m}{2} \ln^2\!{m}\right\}\gamma
-m(m-1)\ln^3\!2+ \frac{m}{2}\ln^3\!{m}   \\[5mm]
&\displaystyle\qquad
- \big[m(m-2)\ln{m} + m(m-1)\ln\pi\big]\ln^2\!2 
+ \frac{3m}{2}\ln2\cdot\ln^2\!{m} +m\ln2\cdot\ln\pi\cdot\ln{m}\\[4mm]
&\displaystyle\qquad
-\frac{(m^2-1)\pi^2}{6}\ln2\pi +\frac{(m^2+2)\pi^2}{12}\ln{m}  
+m(\gamma+\ln2\pi m)\! \sum_{l=1}^{m-1} \ln^2\!\sin\frac{\pi l}{m} \\[4mm]
&\displaystyle\qquad 
+ \frac{m}{2} \!\sum_{l=1}^{m-1} \! \left\{\zeta''\!\left(\! 0,\,\frac{l}{m}\!\right) 
+ \, \zeta''\!\left(\! 0,\,1-\frac{l}{m}\!\right) \! \right\} \!\cdot \ln\sin\frac{\pi l}{m}  
+\pi^2\!\sum_{l=1}^{m-1} l\!\cdot\!\ln\Gamma\!\left(\! \frac{l}{m}\!\right)
\end{array}
\ee
In order to obtain this expression we start from {\eqref{kx98eujms}} and
we successively employ \eqref{cj932230fdncdbds}, \eqref{jk920u3jd209rnd},
\eqref{2piod29rjddfg}, \eqref{jhc23984nd} as well as multiplication
theorems for the logarithm of the $\Gamma$-function and for the $\Psi
$-function
\be\label{kjd203jxndhe8}
\sum_{r=1}^{m-1} \ln\Gamma\biggl(\!\frac{r}{m} \vphantom{\frac{1}{2}} \!\biggr)=\,
\frac{1}{2}(m-1)\ln2\pi-\frac{1}{2}\ln{m}\,, \qquad\qquad
\sum_{r=1}^{m-1} \Psi \biggl(\!\frac{r}{m} \vphantom{\frac{1}{2}} \!\biggr)=\,
\gamma (1-m)-m\ln{m}
\ee 
Note that, generally, when summing the first generalized Stieltjes
constants with an odd function, one arrives at the logarithm of the
$\Gamma$-function, while summing with an even function leads to a
reflected sum of two second-order derivatives of the Hurwitz
$\zeta$-function. The latter sum is the subject of a more detailed
study presented in the next section.

\subsection{Several remarks related to the sum $\zeta''(0,p)+\zeta
''(0,1-p)$}\label{kjd12093ddmne3dd}

From the above formul\ae~it appears that the sum of $\zeta''(0,p)$
with its reflected version $\zeta''(0,1-p)$, at positive rational $p$
less than $1$, plays the fundamental role for the first generalized
Stieltjes constant~at rational argument. We do not know which is the
transcendence of such a sum, but it is not unreasonable to expect that
it is lower than that of solely $\zeta''(0,p)$. Furthermore, in our
previous work \cite[pp.~66--71]{iaroslav_06}, we demonstrated that this
sum has several comparatively simple integral and series
representations; below, we briefly present some of them. In exercises
\no20--21, we dealt with integral $\varPhi(\varphi)$, which we,
unfortunately, could not reduce to elementary functions (despite of its
simple and naive appearance). Written in terms of this integral, the
above sum reads\footnote{Put in \cite[p.~69, Eq.~49]{iaroslav_06}
$\varphi=\pi(2p-1)$.}
\be\label{kjc2wcbfae3}
\begin{array}{rl}
\displaystyle 
\zeta''(0,p)+\zeta''(0,1-p) =&\displaystyle \pi\ctg2\pi p\cdot\Big\{2\ln\Gamma(p) 
+ \ln\sin\pi p + (2p-1)\ln2\pi -\, \ln\pi \Big\}  \\[4mm] 
\displaystyle \phantom{aaaaaaaaa} 
&\displaystyle -2\ln{2\pi}\cdot\ln\!\big(2\sin\pi p\big) \,+
\int\limits_0^\infty \!\!\frac{\,e^{-x}\ln{x}\,}{\,\ch{x}-\cos{2\pi p}\,}\,dx
\end{array}
\ee
where parameter $p$ should lie within the strip $0<\operatorname{Re}{p}<1$.
By a simple change of variable, the last integral may be rewritten in a variety
of other forms. For instance,
\begin{eqnarray}\label{iuh928ychbdbs}
&&\displaystyle \notag
\int\limits_0^\infty \!\!\frac{\,e^{-x}\ln{x}\,}{\,\ch{x}-\cos{2\pi p}\,}\,dx \,=
2\!\int\limits_0^1\!\!\frac{\,x\ln{\ln{\frac{1}{x}}}\,}{\,x^2-2x\cos{2\pi p}+1\,}\,dx \,=
2\!\int\limits_1^\infty \!\!\frac{\,\ln{\ln{x}}\,}{\,x(x^2-2x\cos{2\pi p}+1)\,}\,dx \qquad\\[1mm]
&&\displaystyle \notag \qquad\qquad\qquad\;
=\frac{\pm2}{\sin2\pi p}\cdot\Im\!\int\limits_0^\infty \!\!\frac{\,\ln{x}\,}{\,e^x-e^{\pm 2\pi i p}\,}\,dx\,
=\frac{\pm2}{\sin2\pi p}\cdot\Im\!\int\limits_0^1 \!\!\frac{\,x\ln{\ln{\frac{1}{x}}}\,}{\,x-e^{\pm 2\pi i p}\,}\,dx\, \\[1mm]
&&\displaystyle \qquad\qquad\qquad\qquad\qquad\qquad\qquad\qquad\qquad\quad
=\frac{\pm2}{\sin2\pi p}\cdot\Im\!\int\limits_1^\infty \!\!\frac{\,\ln\ln{x}\,}{\,x\big(x-e^{\pm 2\pi i p}\big)\,}\,dx\,
\end{eqnarray}
The latter forms are particularly simple and display the close
connection to the polylogarithms.
Let now focus our attention on the last integral from the first line.
By partial fraction decomposition it may be written in terms of three
other integrals
\be\label{kj093djl2mwsnf}
\int\limits_1^\infty \!\!\frac{\,\ln{\ln{x}}\,}{\,x^n(x^2-2x\cos{2\pi p}+1)\,}\,dx\,,\qquad
\int\limits_1^\infty \!\!\frac{\,\ln{\ln{x}}\,}{\,x^{k}\,}\,dx\,
\qquad\text{and}\qquad
\int\limits_1^\infty \!\!\frac{\,\ln{\ln{x}}\,}{\,x^2-2x\cos{2\pi p}+1\,}\,dx
\ee
where $n$ and $k$ are positive integers greater than 1.
The values of the last two integrals, thanks to Euler, Legendre and
Malmsten, are known,\footnote{See
\cite[p.~24]{malmsten_01}, \cite[Sect.~4, \no2, 29-h,
30]{iaroslav_06}) or {\eqref{j2039drjm2d}} in \ref{894yf3hedbe}.}
so that the problem of
the evaluation of {\eqref{kjc2wcbfae3}} may be reduced to the first
integral. We, however, note that
the success of this technique depends on the appropriate choice of $p$
and $n$.
Indeed, by expanding the integrand of the first integral in {\eqref
{kj093djl2mwsnf}} into partial fractions, we have
\be\label{908djdsn}
\frac{\,1\,}{\,x^n\!\left(x^2-2x\cos2\pi p+1\right)\,}\,=
\,\frac{\,a_0\,}{\,x\!\left(x^2-2x\cos2\pi p+1\right)\,}+\frac{\,a_1\,}{\,x^2-2x\cos2\pi p+1\,}
+\sum_{l=2}^n \frac{a_l}{x^l}
\ee
with coefficients $a_l$ given by
\be\notag
\begin{array}{lllll}
\displaystyle 
a_0=\frac{\,\sin 2\pi p n\,}{\sin2\pi p}\,, \quad &\displaystyle 
a_1=-\frac{\,\sin 2\pi p(n-1)\,}{\sin2\pi p}\,, \quad&\displaystyle 
a_2=+\frac{\,\sin 2\pi p(n-1)\,}{\sin2\pi p}\,, \quad&\displaystyle 
\ldots, &\\[5mm]
\displaystyle 
a_l=\frac{\,\sin 2\pi p(n-l+1)\,}{\sin2\pi p}\,,\quad
&\displaystyle \ldots,\quad\quad
&\displaystyle a_{n-1}=2\cos2\pi p\,,\quad
&\displaystyle a_n=1\,.
\end{array}
\ee
But if parameter $p$ is such that $a_0=0$, the wanted integral cannot
be collared. The most unpleasant is that this situation occurs
precisely when $p=k/n$, where $k$ is positive integer or demi-integer
not greater than $n$. We, in turn, are able to evaluate
\be\label{kj093dhcg2387bh}
\int\limits_1^\infty \!\!\frac{\,\ln{\ln{x}}\,}{\,x^n(x^2-2x\cos{2\pi p}+1)\,}\,dx\,=\,
\int\limits_0^1\!\!\frac{\,x^n\ln{\ln{\frac{1}{x}}}\,}{\,x^2-2x\cos{2\pi p}+1\,}\,dx\,=
\frac{1}{2}\!\int\limits_0^\infty \!\!\frac{\,e^{-nx}\ln{x}\,}{\,\ch{x}-\cos{2\pi p}\,}\,dx 
\ee
only for those $p$ which may be written as $k/n$, in which case it can
be expressed in terms of $\ln\Gamma(k/n)$ [see \ref{894yf3hedbe}].
Thus, the evaluation of the integral
%
\begin{eqnarray}
\label{kcjhw398} \int\limits_0^\infty \frac{ e^{-nx} \cdot\ln{x} }{ \operatorname
{ch}{x}-\cos\frac{2\pi k}{m} }
\,dx =  2\! \int\limits_0^1 \frac{ x^n\ln{\ln{\frac{1}{x}}} }{ x^2-2x\cos\frac
{2\pi k}{m}+1 } \,dx =  2\! \int\limits
_1^\infty \frac{ \ln{\ln{x}} }{ x^n  \! \left(x^2-2x\cos\frac
{2\pi k}{m}+1 \right) } \,dx\nonumber\\[-5pt]
\end{eqnarray}
with $n=2, 3, 4, \ldots{}$, number $m$ being positive integer such that
$m\neq2kn/l$ for $l=\pm1, \pm2, \pm3,\ldots{}$, could bring the
solution to our problem, but
currently we do not know if this integral can be evaluated in terms of
lower transcendental functions.
However, it should be noted that
integrals closely related to {\eqref{kcjhw398}} and {\eqref{kjh298hdnd}}
were a subject of several investigations appeared already in the XIXth century.
The most significant contribution seems to belong to Malmsten who
showed in 1842 that
%
\begin{eqnarray}
\label{kew9e3jd} \frac{\sin a}{\Gamma(s)} \int\limits_0^1
\frac{ x^y \cdot\ln^{s-1} \frac
{1}{x} }{ x^2+2x\cos{a}+1 } \,dx &= & \int\limits_0^\infty
\frac{\operatorname{sh}ax}{\operatorname{sh}\pi
x}\cdot \frac{\cos \big(s \operatorname{arctg}\frac{x}{y} \big)}{
(x^2+y^2 )^{s/2}} \,dx= \sum
_{l=1}^\infty(-1)^{l-1} \frac{\sin al}{(y+l)^s}\qquad\;\;
\end{eqnarray}
$y, s\in\mathbbm{C}$, $-\pi<a<+\pi$, see \cite[pp.~20--25]{malmsten_00}
and \cite[p.~12]{malmsten_01}. He studied these integrals for different
values of parameters $y, s$ and $a$, and evaluated some of them in a
closed form. The above equality permitted to Malmsten to derive
numerous fascinating results, such as, for example,
formul\ae~{\eqref{odi239dn}} and (\ref{higce672ecbx}b).
Furthermore, his investigations devoted to the cases $y=0$, $a=\pi/2$
and $y=0$, $a=\pi/3$ resulted in two important reflection formul\ae~for
the $L$- and $M$-functions respectively \cite[p.~23,
Eq.~(36)]{malmsten_00}, \cite[pp.~17--18,
Eqs.~(51)--(52)]{malmsten_01}, \cite[pp.~35--36, Eq.~(21),
Fig.~3]{iaroslav_06} (these formul\ae~are similar to Euler--Riemann's reflection formula for the $\zeta$-function, see also
footnote \ref{lkd024kdfm}). Notwithstanding, Malmsten failed to show
that more generally, when $a$ is a rational multiple of $\pi$, one has
%
\begin{eqnarray}
\label{jd2093jmf} \sum_{l=1}^\infty(-1)^{l-1}
\frac{\sin al}{(y+l)^s} = \frac{1}{ (2n)^s }\sum_{l=1}^{2n-1}
(-1)^{l-1}\sin\frac{\pi ml}{n} \cdot\zeta \left( s,
\frac{ y+l }{2n} \right) , \quad a\equiv\frac{ m\pi}{n}
\end{eqnarray}
$m=1, 2, 3,\ldots, n-1$, which may be obtained by applying
Hurwitz's method used in
\eqref{khjx98qh2xbn21}--\eqref{jhcwiknx83} 
to series {\eqref{kew9e3jd}}.\footnote{Actually, Malmsten also
studied the case $a=m\pi/n$, but quite superficially and mainly for
$y=0$.}$\stackrel{,}{\vphantom{.}}\,$\footnote{Note that for $s=1, 2, 3,\ldots$ the right
part of {\eqref{jd2093jmf}} reduces to polygamma functions, see
e.g. \cite[pp.~71--72, \no23]{iaroslav_06}.} Now, Malmsten's integrals
from {\eqref{kew9e3jd}} are related to ours from
{\eqref{iuh928ychbdbs}} as follows
%
\begin{eqnarray}
\int\limits_0^1 \frac{ x \ln{\ln\frac{1}{x}} }{ x^2-2x\cos{2\pi p}+1 } \,dx =\lim
_{s\to1} \left\{\frac{\partial}{\partial s} \int\limits
_0^1 \frac{ x\cdot\ln^{s-1} \frac{1}{x} }{ x^2-2x\cos{2\pi
p}+1 } \,dx \right\}
\end{eqnarray}
Therefore, by {\eqref{kew9e3jd}} we have
%
\begin{eqnarray}
\label{lkjwc983ncnds} %
\int\limits_0^1
\frac{ x\cdot\ln^{s-1} \frac{1}{x} }{ x^2-2x\cos{2\pi
p}+1 } \,dx & = & -\frac{\Gamma(s)}{\sin2\pi p} \int\limits_0^\infty
\frac{\operatorname
{sh}  [\pi(2p-1)x  ]}{\operatorname{sh}\pi x}\cdot \frac{\cos(s \operatorname{arctg}x)}{ (x^2+1 )^{s/2}} \,dx
\nonumber
\\
& =& -\frac{\Gamma(s)}{2\sin2\pi p} \int\limits_{-\infty}^{+\infty}
\frac{\operatorname{sh}  [\pi(2p-1)x
 ]}{\operatorname{sh}\pi x}\cdot \frac{dx}{ (1\pm ix )^{s}} , \quad\operatorname{Re} {s}>0 .\nonumber\\[-5pt]
\end{eqnarray}
Integrals appearing on the right-hand side are also quite similar to
Jensen's formul\ae~for $\zeta(s)$ derived between
1893 and 1895 by contour integration methods, see \cite{jensen_04} and
\cite{jensen_03}.
Taking into account that these references are hard to find and that the same
formul\ae~were later reprinted with misprints,\footnote{In the
well-known monograph
\cite[vol.~I]{bateman_01}, in formula (13) on p.~33,
``$(e^{2\pi t}+1)^{-t}$'' should be replaced by ``$(e^{\pi t}+1)^{-t}$''.}
we find it useful to reproduce them here as well
\be\label{kjd02jddnsa}
\begin{array}{rl}
&\displaystyle 
\zeta(s) = \frac{1}{s-1} + \frac{1}{2} + 2\!\int\limits_0^{\pi/2} \! 
\frac{(\cos\theta)^{s-2}\sin s\theta}{e^{2\pi\tg\theta}-1} d\theta  \,=\,
\frac{1}{s-1} + \frac{1}{2} + 2\!
\int\limits_0^\infty \! 
\frac{\sin(s \arctg x)\, dx}{\left(e^{2\pi x}-1\right) \left(x^2+1\right)^{s/2}}
\\[8mm]
&\displaystyle 
\zeta(s) = \frac{2^{s-1}}{s-1} \,+\, i\, 2^{s-1}\!\!\int\limits_0^\infty \! 
\frac{\left(1+ix\right)^{s}-\left(1-ix\right)^{s}}{\left(e^{\pi x}+1\right)\left(x^2+1\right)^{s}} \,dx=
\frac{2^{s-1}}{s-1} \,-\, \, 2^s\!\!\int\limits_0^\infty \! 
\frac{\sin(s \arctg x)\, dx}{\left(e^{\pi x}+1\right) \left(x^2+1\right)^{s/2}}
\\[8mm]
&\displaystyle 
\zeta(s) = \frac{\pi}{2(s-1)}\!\!\int\limits_{-\infty}^{+\infty} \! 
\frac{1}{\ch^2\!\pi x} \cdot \frac{dx}{\big(\frac{1}{2}+ix\big)^{s-1}} \,=
\,\frac{\pi\,2^{s-2}}{s-1}\!\!\int\limits_{0}^{\infty} \! 
\frac{\cos\big[(s-1) \arctg x\big]}{\left(x^2+1\right)^{(s-1)/2}\ch^2\!\frac{1}{2}\pi x}\, dx
\end{array}
\ee
$s\in\mathbbm{C}$, $s\neq1$, where final simplifications were done later
by Lindel\"of \cite[p.~103]{lindelof_01} who also gave details of their
derivation.\footnote{Jensen did not provide proofs for these
formul\ae; he only stated that he had found them in his
notes,\footnotemark and added that they can be
easily derived by Cauchy's residue theorem. By the way, the first
of these three formul\ae~was also obtained by Franel
\cite{jensen_03,franel_01,jensen_04}.} 
\footnotetext{\emph{Je trouve encore, dans mes notes, entre autres, les formules\ldots} \cite{jensen_03}.}
Application of contour
integration methods to integrals {\eqref{lkjwc983ncnds}} seems
quite attractive as well (especially if $p$ is rational), but the
branch point at $\pm i$ is really annoying.

Other representations for $\zeta''(0,p)+\zeta''(0,1-p)$ may also
involve integrals
\begin{eqnarray}
\nonumber
\int\limits_0^\infty\frac{\ln (x^2+p^2 )\cdot\operatorname
{arctg} (x/p )}{e^{2\pi x}-1} \,dx
\quad\mbox{or}\quad \int\limits_0^\infty
\frac{\ln^2(ip+x)-\ln^2(ip-x)}{e^{2\pi x}-1} \,dx
\end{eqnarray}
which directly follow from the well-known Hermite representation for
the Hurwitz $\zeta$-function
\cite[p.~66]{hermite_01}, \cite[p.~106]{lindelof_01}, \cite[vol.~I,
p.~26, Eq.~1.10(7)]{bateman_01}.

The sum $\zeta''(0,p)+\zeta''(0,1-p)$ may be also reduced to an
important logarithmic--trigonometric series
\begin{eqnarray}
\nonumber
\zeta''(0,p)+\zeta''(0,1-p)
= -2(\gamma+\ln2\pi)\ln (2\sin \pi p ) +2 \sum_{n=1}^\infty\!
\frac{ \cos2\pi p n \cdot\ln{n} }{n}
\end{eqnarray}
see \cite[p.~69, \no22]{iaroslav_06}.
This series, unlike the similar sine-series, is not known to be
reducible to any elementary or classical function of analysis;
however, it was remarked by Landau \cite[pp.~180--182]{landau_02} that
it has some common properties with the logarithm
of the $\Gamma$-function. Besides, it also appeared in works of Lerch
\cite{lerch_02} and Gut \cite{gut_01}.

Another way to treat the problem could be to use the antiderivatives of
the first generalized Stieltjes constant~$\Gamma_1(p)$. In
\cite[p.~69, \no22]{iaroslav_06}, we showed that the sum $\zeta
''(0,p)+\zeta''(0,1-p)$
may be also written in terms of such functions
\begin{eqnarray}
\nonumber
\zeta''(0,p)+\zeta''(0,1-p)
= -(3\ln2+2\ln\pi)\ln2 - 4\Gamma _1(1/2)+2\Gamma_1(p)+2
\Gamma_1(1-p)
\end{eqnarray}
The latter
formula, inserted into {\eqref{kx98eujms}}, gives an equation which is
in some way analogous Malmsten's
representation for the Digamma function (\ref{dh38239djws}c)
[in the sense that for rational arguments it provides a connection
between the function and its derivative].

Finally, note that almost all above expressions remain valid everywhere
in the strip \mbox{$0<\operatorname{Re}{p}<1$},
so it is not impossible that for rational $p$ they could be further
simplified or reduced to less transcendental forms.
Thus, the question of the possibility to express any first generalized
Stieltjes constant~of a rational argument
not only via the $\Gamma$-function, $\gamma_1$, $\gamma$ and some
``relatively simple''
function, but solely via the $\Gamma$-function, $\gamma_1$, $\gamma
$ and elementary functions
remains open and is directly connected to the transcendence of the
reflected sum $\zeta''(0,p)+\zeta''(0,1-p)$ at rational $p$,
which is currently not sufficiently well studied.

\section{Extensions of the theorem to the second and higher Stieltjes
constants}\label{iou2ch20983j3e}
It can be reasonably expected that similar theorems could be derived
for the higher Stieltjes constants.
Such a demonstration could be carried out again with the help of
$J_a(p)$ and integral {\eqref{k290djd2s}} where $\ln x$
is replaced with $\ln^n x$ [see below how integral {\eqref
{kjc294djcwx}} is used for the determination of the second Stieltjes constant].
As regards the equation for the difference between generalized
Stieltjes constants, which is also necessary,
it is simply sufficient to note that from {\eqref{dhd73vj6s2}} and
{\eqref{khjhbcubdbvxhbs}} it follows that
\be\notag
\begin{array}{l}
\displaystyle
\gamma_n \biggl(\frac{r}{m} \vphantom{\frac{1}{2}} \biggr)
- \gamma_n \biggl(1-\frac{r}{m} \vphantom{\frac{1}{2}}
\biggr) = (-1)^n \lim_{a\to1}\!
\left\{\zeta^{(n)}\!\biggl(\!a,\frac{r}{m} \vphantom{\frac{1}{2}} \biggr)
- \zeta^{(n)}\!\biggl(\!a,1-\frac{r}{m} \vphantom{\frac{1}{2}}
\biggr)\!\right\}=\\[6mm]
\displaystyle\qquad\qquad\qquad
=4(-1)^n \lim_{a\to1}\frac{\partial^n}{\partial a^n}\!
\left\{ \!\frac{\Gamma(1-a)}{(2\pi m)^{1-a}}\,\cos\frac{\pi a}{2} \cdot\!
\sum_{l=1}^{m-1} \sin\frac{2\pi rl}{m}\cdot\zeta\!\left(\!1-a,\,\frac{l}{m}\!\right)  \!\right\}
\end{array}
\ee
$n=1, 2, 3, \ldots$ and
where $r$ and $m$ are positive integers such that $r<m$. In particular,
for the second generalized Stieltjes constant, the latter
formula takes the form\footnote{This formula also appears in an
unpublished work sent to the author by Donal Connon.}
\be\label{dwec8992wj}
\begin{array}{rl}
\displaystyle
\gamma_2 \biggl(\frac{r}{m} \vphantom{\frac{1}{2}} \biggr)
- \gamma_2 \biggl(1-\frac{r}{m} \vphantom{\frac{1}{2}}
\biggr) =&\,\, \displaystyle2\pi\!\sum_{l=1}^{m-1}
\sin\frac{2\pi r l}{m} \cdot\zeta''\!\left(\!0,\,\frac{l}{m}\!\right) +\pi\left[\frac{\pi^2}{12}+(\gamma+\ln2\pi m)^2\right]
\ctg\frac{\pi r}{m}
 \\[6mm]
\displaystyle
&\displaystyle
 - 4\pi(\gamma+\ln2\pi m)\!
\sum_{l=1}^{m-1}
\sin\frac{2\pi r l}{m} \cdot\ln\Gamma \biggl(\frac{l}{m} \biggr) 
\end{array}
\ee
In order to obtain a formula for $\gamma_2(r/m)$, we take again
expansion {\eqref{kjhg87fgf}} and write down its terms up to $O(a-1)^3$. Hence
\be\notag
\begin{array}{ll}
\displaystyle
\int\limits_0^\infty \!\!\frac{\, (\ch{[(2p-1)x]}-1)\ln^2{x}\,}{\,\sh{x}\,}\,dx\,=  \,
\frac{2 }{3}C_m(r)
-2B_m(r)\ln\pi m + \left\{2\ln^2\pi m + \frac{\pi^2}{6}\right\}A_m(r) - 
\end{array}
\ee

\be\notag
\begin{array}{ll}
\displaystyle\qquad\qquad\qquad\qquad\qquad\quad
-2\gamma_1(\gamma-\ln2) 
+\frac{2}{3}\zeta(3) - \frac{2}{3}\gamma^3 -\gamma_2 + \left(\!\gamma^2-\frac{\pi^2}{6}\right)\!\ln2 \\[5mm]
\displaystyle\qquad\qquad\qquad\qquad\qquad\quad
+\frac{\pi^2}{12}\ln\pi m + \ln\pi\cdot\ln m\cdot\ln\pi m +\frac{1}{3}\big(\ln^3\pi + \ln^3 m - \ln^32\big)
\end{array}
\ee
where $p\equiv r/m$ and 
\be\notag
\begin{array}{rl}
\displaystyle
C_m(r) \equiv \sum_{l=1}^m \cos\frac{2\pi rl}{m}\cdot\zeta'''\!\left(\!0,\,\frac{l}{m}\!\right) =& \displaystyle \,
\sum_{l=1}^{m-1} \cos\frac{2\pi rl}{m}\cdot\zeta'''\!\left(\!0,\,\frac{l}{m}\!\right)
+\frac{3}{2}\gamma_2  + \gamma^3 - \zeta(3) \\[6mm]
\displaystyle
& \displaystyle 
+3\gamma_1\gamma - \frac{1}{2}\ln^3\!{2\pi}
+\left\{3\gamma_1+ \frac{3}{2}\gamma^2 - \frac{\pi^2}{8}\right\}\ln2\pi
\end{array}
\ee
Comparing the latter integral to {\eqref{kjc294djcwx}} and then using
{\eqref{d33f0rfifm}}, we obtain
\be\notag
\begin{array}{l}
\displaystyle
\gamma_2 \biggl(\frac{r}{m} \vphantom{\frac{1}{2}} \biggr)
+ \gamma_2 \biggl(1-\frac{r}{m} \vphantom{\frac{1}{2}}
\biggr) = 2\gamma_2 - 4\gamma_1 \ln{m}  +2 \gamma^{2}\ln{2} + \frac{4}{3}\!\sum_{l=1}^{m-1}
\cos\frac{2\pi r l}{m} \cdot\zeta'''\!\left(\!0,\,\frac{l}{m}\!\right) \\[6mm]
\displaystyle\qquad
 - 4 (\gamma+\ln2\pi m)\! \sum_{l=1}^{m-1}
\cos\frac{2\pi r l}{m} \cdot\zeta''\!\left(\!0,\,\frac{l}{m}\!\right) 
+2\!\left[\frac{\pi^2}{12}-(\gamma+\ln2\pi m)^2\right]\!\times \\[6mm]
\displaystyle\qquad
\times\!\sum_{l=1}^{m-1} \cos\frac{\,2\pi r l\,}{m} \cdot\ln\sin\frac{\pi l}{m}
 -\frac{\pi^{2}}{6}\ln{2}
+ 2 \gamma\big(\ln^{2}\!{m}+2\ln^{2}\!{2}+2 \ln{2}\cdot\ln{\pi m}\big) \\[6mm]
\displaystyle\qquad
+2\big(\ln^{2}\!{2}+\ln^{2}\!{m}+\ln^{2}\!{\pi}+2 \ln{\pi}\ln{m}+2 \ln{2}\ln{\pi m}\big)\ln{2}
+ \frac{2}{3}\ln^{3}{m}
\end{array}
\ee
which, being added to \eqref{dwec8992wj}, finally yields
\be\label{nolabel}
\begin{array}{rl}
\displaystyle
\gamma_2 \biggl(\frac{r}{m} \vphantom{\frac{1}{2}} \biggr) = &\displaystyle
\gamma_2  + \frac{2}{3}\!\sum_{l=1}^{m-1}
\cos\frac{2\pi r l}{m} \cdot\zeta'''\!\left(\!0,\,\frac{l}{m}\!\right)  -
2 (\gamma+\ln2\pi m)\! \sum_{l=1}^{m-1}
\cos\frac{2\pi r l}{m} \cdot\zeta''\!\left(\!0,\,\frac{l}{m}\!\right) \\[6mm]
\displaystyle &\displaystyle
+ \pi\!\sum_{l=1}^{m-1}
\sin\frac{2\pi r l}{m} \cdot\zeta''\!\left(\!0,\,\frac{l}{m}\!\right) 
-2\pi(\gamma+\ln2\pi m)\!
\sum_{l=1}^{m-1}
\sin\frac{2\pi r l}{m} \cdot\ln\Gamma \biggl(\frac{l}{m} \biggr)    \\[6mm]
\displaystyle&\displaystyle
+\left[\frac{\pi^2}{12}-(\gamma+\ln2\pi m)^2\right]\!\cdot\!
\sum_{l=1}^{m-1} \cos\frac{\,2\pi r l\,}{m} \cdot\ln\sin\frac{\pi l}{m} + \gamma^{2}\ln{2}   - 2\gamma_1 \ln{m} 
  \\[6mm]
\displaystyle &\displaystyle
+\left[\frac{\pi^2}{12}+(\gamma+\ln2\pi m)^2\right]\cdot\frac{\pi}{2}\ctg\frac{\pi r}{m} 
+ \gamma\big(\ln^{2}\!{m}+2\ln^{2}\!{2}+2 \ln{2}\cdot\ln{\pi m}\big) -   \\[6mm]
\displaystyle &\displaystyle
- \frac{\pi^{2}}{12}\ln{2}
+\big(\ln^{2}\!{2}+\ln^{2}\!{m}+\ln^{2}\!{\pi}+2 \ln{\pi}\ln{m}+2 \ln{2}\ln{\pi m}\big)\ln{2} 
+ \frac{1}{3}\ln^{3}{m} 
\end{array}
\ee
This formula is an analog of {\eqref{kx98eujms}} for the second
generalized Stieltjes constant. It can be also reduced to other
forms if necessary. For instance, similarly to {\eqref{khjwe9jdns}}, we
may rewrite it in
the form containing the $\Psi$-function
\be\notag
\begin{array}{rl}
\displaystyle
\gamma_2 \biggl(\frac{r}{m} \vphantom{\frac{1}{2}} \biggr) = \,
\gamma_2 + \frac{2}{3}\!\sum_{l=1}^{m-1}
\cos\frac{2\pi r l}{m} \cdot\zeta'''\!\left(\!0,\,\frac{l}{m}\!\right)  -
2 (\gamma+\ln2\pi m)\! \sum_{l=1}^{m-1}
\cos\frac{2\pi r l}{m} \cdot\zeta''\!\left(\!0,\,\frac{l}{m}\!\right) \\[6mm]
\displaystyle \quad
+ \pi\!\sum_{l=1}^{m-1}
\sin\frac{2\pi r l}{m} \cdot\zeta''\!\left(\!0,\,\frac{l}{m}\!\right) 
-2\pi(\gamma+\ln2\pi m)\!
\sum_{l=1}^{m-1}
\sin\frac{2\pi r l}{m} \cdot\ln\Gamma \biggl(\frac{l}{m} \biggr) 
 - 2\gamma_1 \ln{m}    \\[6mm]
\displaystyle\quad
- \gamma^3 
-\left[(\gamma+\ln2\pi m)^2-\frac{\pi^2}{12}\right]\!\cdot\!
\Psi\!\biggl(\frac{r}{m} \vphantom{\frac{1}{2}} \biggr) + 
\frac{\pi^3}{12}\ctg\frac{\pi r}{m} 
  -\gamma^2\ln\big(4\pi^2 m^3\big) +\frac{\pi^2}{12}(\gamma+\ln{m}) \\[6mm]
\displaystyle\quad
 - \gamma\big(\ln^2\!{2\pi} +4\ln{m}\cdot\ln{2\pi}+2\ln^2\!{m}\big)
 -\left\{\!\ln^2\!{2\pi}+2\ln{2\pi}\cdot\ln{m}+\frac{2}{3}\ln^2\!{m}\!\right\}\!\ln{m}
\end{array}
\ee
Thus, corresponding expressions for higher generalized Stieltjes
constants~at rational points are expected to be quite long and to
contain higher derivatives of the
Hurwitz zeta-function at zero at rational points $\zeta^{(n)}(0,l/m)$
whose properties are currently little studied.

\section*{Acknowledgments}
The author would like to thank anonymous reviewers for their
suggestions and remarks which improved the quality of the manuscript.
The author is also grateful to Donal F.~Connon for sending alternative
demonstrations of some of the results derived in this
paper,\footnote{Private communication.} as well as to Mark W.~Coffey
for providing Refs.~\cite{wilton_01} and \cite{kluyver_01}. Finally, 
the author is very grateful to Stefan Kr\"amer for his precious remarks related to series {\eqref{lk2jd029jde}}, 
for drawing attention to Mascheroni's work \cite{mascheroni_01} and for providing Refs.~\cite{jacobsthal_01} 
and \cite{skramer_01}.

\vspace{19mm}


\begin{figure}[!h]   
\centering
\includegraphics[width=0.3\textwidth]{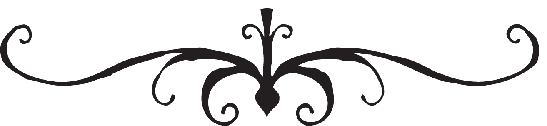}
\end{figure}

\newpage

\appendix
\section{Closed-form expressions for some Stieltjes constants}\label{pi3ofj94cm3pimc3opij}

In this first supplementary part of our work, we provide some
information about particular values of $\gamma_1(v)$ which are free
from $\zeta''(0,l/m)+\zeta''(0,1-l/m)$ or which contain only one
combination of it. The value of $\gamma_1(1/2)$ has been long-time
known and may be found in numerous works. The values of
$\gamma_1(1/4)$, $\gamma_1(3/4)$ and $\gamma_1(1/3)$ were independently
obtained by Donal Connon in \cite[pp.~17--18]{connon_03} and by the
author in \cite[p.~100]{iaroslav_06}. Closed-form expressions for
$\gamma_1(2/3)$, $\gamma_1(1/6)$ and $\gamma_1(5/6)$ were given by the
author in \cite[pp.~100--101]{iaroslav_06}. All these values do not
contain the Hurwitz $\zeta$-function. Below, we provide some further
values which may be of interest and which may be reduced to only one
transcendent $\zeta''(0,l/m)+\zeta''(0,1-l/m)$.
\be\notag
\begin{array}{ll}
\displaystyle 
 \gamma_1&\displaystyle\biggl(\!\frac{1}{5} \vphantom{\frac{1}{2}} \!\biggr)
 = \;\gamma_1  + \frac{\sqrt{5}}{2}\!\left\{\zeta''\!\left(\! 0,\,\frac{1}{5}\vphantom{\dfrac{7}{2}}\!\right) 
+ \zeta''\!\left(\! 0,\,\frac{4}{5}\vphantom{\dfrac{7}{2}}\!\right)\!\right\}
+ \frac{\pi\sqrt{10+2\sqrt5}}{2} \ln\Gamma \biggl(\!\frac{1}{5} \!\biggr)

\\[6mm]
& \displaystyle 
+ \frac{\pi\sqrt{10-2\sqrt5}}{2} \ln\Gamma \biggl(\!\frac{2}{5} \!\biggr)
+\left\{\!\frac{\sqrt{5}}{2} \ln{2} -\frac{\sqrt{5}}{2} \ln\!\big(1+\sqrt{5}\big) -\frac{5}{4}\ln5
-\frac{\pi\sqrt{25+10\sqrt5}}{10} \right\}\!\cdot\gamma \\[5mm]
& \displaystyle 
- \frac{\sqrt{5}}{2}\left\{\ln2+\ln5+\ln\pi+\frac{\pi\sqrt{25-10\sqrt5}}{10}\right\}\!\cdot\ln\!\big(1+\sqrt{5})    
+\frac{\sqrt{5}}{2}\ln^2\!2 + \frac{\sqrt{5}\big(1-\sqrt{5}\big)}{8}\ln^2\!5   \\[5mm]
& \displaystyle 
 +\frac{3\sqrt{5}}{4}\ln2\cdot\ln5 + \frac{\sqrt{5}}{2}\ln2\cdot\ln\pi+\frac{\sqrt{5}}{4}\ln5\cdot\ln\pi
- \frac{\pi\big(2\sqrt{25+10\sqrt5}+5\sqrt{25+2\sqrt5}  \big)}{20}\ln2\\[5mm]
& \displaystyle 
- \frac{\pi\big(4\sqrt{25+10\sqrt5}-5\sqrt{5+2\sqrt5}  \big)}{40}\ln5
- \frac{\pi\big(5\sqrt{5+2\sqrt5}+\sqrt{25+10\sqrt5}  \big)}{10}\ln\pi\\[5mm]
& \displaystyle 
= -8.030205511\ldots
\end{array}
\ee
Stieltjes constants $\gamma_1(2/5)$, $\gamma_1(3/5)$ and $\gamma_1(4/5)$ may be similarly expressed
in terms of $\zeta''(0,1/5)+\zeta''(0,4/5)$, $\Gamma(1/5)$, $\Gamma(2/5)$, $\gamma_1$, $\gamma$ and elementary functions, 
which, by the way, contain the golden ratio $\phi$.
\be\notag
\begin{array}{ll}
\displaystyle 
 \gamma_1\biggl(\!\frac{1}{8} \vphantom{\frac{1}{2}} \!\biggr)
 =& \displaystyle\;\gamma_1  + \sqrt{2}\left\{\zeta''\!\left(\! 0,\,\frac{1}{8}\vphantom{\dfrac{7}{2}}\!\right) 
+ \zeta''\!\left(\! 0,\,\frac{7}{8}\vphantom{\dfrac{7}{2}}\right)\!\right\}
+ 2\pi\sqrt{2}\ln\Gamma \biggl(\!\frac{1}{8} \!\biggr)
-\pi \sqrt{2}\big(1-\sqrt2\big)\ln\Gamma \biggl(\!\frac{1}{4} \!\biggr)
\\[5mm]
& \displaystyle 
-\left\{\!\frac{1+\sqrt2}{2}\pi+4\ln{2} +\sqrt{2}\ln\!\big(1+\sqrt{2}\big) \!\right\}\!\cdot\gamma 
- \frac{1}{\sqrt{2}}\big(\pi+8\ln2+2\ln\pi\big)\!\cdot\ln\!\big(1+\sqrt{2})  
\\[5mm]
& \displaystyle 
 - \frac{7\big(4-\sqrt2\big)}{4}\ln^2\!2   +  \frac{1}{\sqrt{2}}\ln2\cdot\ln\pi 
 -\frac{\pi\big(10+11\sqrt2\big)}{4}\ln2
  -\frac{\pi\big(3+2\sqrt2\big)}{2}\ln\pi\\[5mm]
& \displaystyle 
= -16.64171976\ldots \\[5mm]
\displaystyle 
 \gamma_1\biggl(\!\frac{3}{8} \vphantom{\frac{1}{2}} \!\biggr)
 =& \displaystyle\;\gamma_1  - \sqrt{2}\left\{\zeta''\!\left(\! 0,\,\frac{1}{8}\vphantom{\dfrac{7}{2}}\!\right) 
+ \zeta''\!\left(\! 0,\,\frac{7}{8}\vphantom{\dfrac{7}{2}}\right)\!\right\}
+  2\pi\sqrt{2}\ln\Gamma \biggl(\!\frac{1}{8} \!\biggr)
-\pi \sqrt{2}\big(1+\sqrt2\big)\ln\Gamma \biggl(\!\frac{1}{4} \!\biggr)
\\[5mm]
& \displaystyle 
+\left\{\!\frac{1-\sqrt2}{2}\pi -4\ln{2} +\sqrt{2}\ln\!\big(1+\sqrt{2}\big) \!\right\}\!\cdot\gamma 
- \frac{1}{\sqrt{2}}\big(\pi-8\ln2-2\ln\pi\big)\!\cdot\ln\!\big(1+\sqrt{2})  
\\[5mm]
& \displaystyle 
 - \frac{7\big(4+\sqrt2\big)}{4}\ln^2\!2   -  \frac{1}{\sqrt{2}}\ln2\cdot\ln\pi 
 +\frac{\pi\big(10-11\sqrt2\big)}{4}\ln2
  +\frac{\pi\big(3-2\sqrt2\big)}{2}\ln\pi\\[5mm]
& \displaystyle 
= -2.577714402\ldots\\[-5mm]
\end{array}
\ee

\be\notag
\begin{array}{ll}
\displaystyle 
 \gamma_1\biggl(\!\frac{5}{8} \vphantom{\frac{1}{2}} \!\biggr)
 =& \displaystyle\;\gamma_1  - \sqrt{2}\left\{\zeta''\!\left(\! 0,\,\frac{1}{8}\vphantom{\dfrac{7}{2}}\!\right) 
+ \zeta''\!\left(\! 0,\,\frac{7}{8}\vphantom{\dfrac{7}{2}}\right)\!\right\}
-  2\pi\sqrt{2}\ln\Gamma \biggl(\!\frac{1}{8} \!\biggr)
+\pi \sqrt{2}\big(1+\sqrt2\big)\ln\Gamma \biggl(\!\frac{1}{4} \!\biggr)
\\[5mm]
& \displaystyle 
-\left\{\!\frac{1-\sqrt2}{2}\pi+4\ln{2} -\sqrt{2}\ln\!\big(1+\sqrt{2}\big) \!\right\}\!\cdot\gamma 
+ \frac{1}{\sqrt{2}}\big(\pi+8\ln2+2\ln\pi\big)\!\cdot\ln\!\big(1+\sqrt{2})  
\\[5mm]
& \displaystyle 
 - \frac{7\big(4+\sqrt2\big)}{4}\ln^2\!2   -  \frac{1}{\sqrt{2}}\ln2\cdot\ln\pi 
 -\frac{\pi\big(10-11\sqrt2\big)}{4}\ln2
  -\frac{\pi\big(3-2\sqrt2\big)}{2}\ln\pi\\[5mm]
& \displaystyle 
= -0.7353809459\ldots\\[6mm]
\displaystyle 
 \gamma_1\biggl(\!\frac{7}{8} \vphantom{\frac{1}{2}} \!\biggr)
 =& \displaystyle\;\gamma_1  + \sqrt{2}\left\{\zeta''\!\left(\! 0,\,\frac{1}{8}\vphantom{\dfrac{7}{2}}\!\right) 
+ \zeta''\!\left(\! 0,\,\frac{7}{8}\vphantom{\dfrac{7}{2}}\right)\!\right\}
-  2\pi\sqrt{2}\ln\Gamma \biggl(\!\frac{1}{8} \!\biggr)
+\pi \sqrt{2}\big(1-\sqrt2\big)\ln\Gamma \biggl(\!\frac{1}{4} \!\biggr)
\\[5mm]
& \displaystyle 
+\left\{\!\frac{1+\sqrt2}{2}\pi -4\ln{2} -\sqrt{2}\ln\!\big(1+\sqrt{2}\big) \!\right\}\!\cdot\gamma 
+ \frac{1}{\sqrt{2}}\big(\pi-8\ln2-2\ln\pi\big)\!\cdot\ln\!\big(1+\sqrt{2})  
\\[5mm]
& \displaystyle 
 - \frac{7\big(4-\sqrt2\big)}{4}\ln^2\!2   + \frac{1}{\sqrt{2}}\ln2\cdot\ln\pi 
 +\frac{\pi\big(10+11\sqrt2\big)}{4}\ln2
  +\frac{\pi\big(3+2\sqrt2\big)}{2}\ln\pi\\[5mm]
& \displaystyle 
= -0.1906592305\ldots \\[6mm]
\displaystyle 
 \gamma_1\biggl(\!\frac{1}{12} \vphantom{\frac{1}{2}} \!\biggr)
 =& \displaystyle\;\gamma_1  + \sqrt{3}\left\{\zeta''\!\left(\! 0,\,\frac{1}{12}\vphantom{\dfrac{7}{2}}\!\right) 
+ \zeta''\!\left(\! 0,\,\frac{11}{12}\vphantom{\dfrac{7}{2}}\right)\!\right\}
+ 4\pi\ln\Gamma \biggl(\!\frac{1}{4} \!\biggr)
+3\pi \sqrt{3}\ln\Gamma \biggl(\!\frac{1}{3} \!\biggr)
\\[5mm]
& \displaystyle 
-\left\{\!\frac{2+\sqrt3}{2}\pi+\frac{3}{2}\ln3 -\sqrt3(1-\sqrt3)\ln{2} +2\sqrt{3}\ln\!\big(1+\sqrt{3}\big) \!\right\}\!\cdot\gamma 
\\[5mm]
& \displaystyle 
- 2\sqrt3\big(3\ln2+\ln3 +\ln\pi\big)\!\cdot\ln\!\big(1+\sqrt{3})  
 - \frac{7-6\sqrt3}{2}\ln^2\!2  - \frac{3}{4}\ln^2\!3  \\[5mm]
& \displaystyle 
+  \frac{3\sqrt3(1-\sqrt3)}{2}\ln3\cdot\ln2
 +  \sqrt3\ln2\cdot\ln\pi 
 -\frac{\pi\big(17+8\sqrt3\big)}{2\sqrt3}\ln2 \\[5mm]
& \displaystyle 
 +\frac{\pi\big(1-\sqrt3\big)\sqrt3}{4}\ln3
 -\pi\sqrt3(2+\sqrt3)\ln\pi
= -29.84287823\ldots \\[6mm]
\displaystyle 
 \gamma_1\biggl(\!\frac{7}{12} \vphantom{\frac{1}{2}} \!\biggr)
 =& \displaystyle\;\gamma_1  - \sqrt{3}\left\{\zeta''\!\left(\! 0,\,\frac{1}{12}\vphantom{\dfrac{7}{2}}\!\right) 
+ \zeta''\!\left(\! 0,\,\frac{11}{12}\vphantom{\dfrac{7}{2}}\right)\!\right\}
- 4\pi\ln\Gamma \biggl(\!\frac{1}{4} \!\biggr)
+3\pi \sqrt{3}\ln\Gamma \biggl(\!\frac{1}{3} \!\biggr)
\\[5mm]
& \displaystyle 
-\left\{\!\frac{-2+\sqrt3}{2}\pi+\frac{3}{2}\ln3 +\sqrt3(1+\sqrt3)\ln{2} -2\sqrt{3}\ln\!\big(1+\sqrt{3}\big) \!\right\}\!\cdot\gamma 
\\[5mm]
& \displaystyle 
+2\sqrt3\big(3\ln2+\ln3 +\ln\pi\big)\!\cdot\ln\!\big(1+\sqrt{3})  
 - \frac{7+6\sqrt3}{2}\ln^2\!2  - \frac{3}{4}\ln^2\!3  \\[5mm]
& \displaystyle 
-  \frac{3\sqrt3(1+\sqrt3)}{2}\ln3\cdot\ln2
 -  \sqrt3\ln2\cdot\ln\pi 
 -\frac{\pi\big(17-8\sqrt3\big)}{2\sqrt3}\ln2 \\[5mm]
& \displaystyle 
+\frac{\pi\big(1+\sqrt3\big)\sqrt3}{4}\ln3 
-\pi\sqrt3(2-\sqrt3)\ln\pi
= -0.900932495\ldots
\end{array}
\ee
Expressions for Stieltjes constants $\gamma_1(5/12)$ and $\gamma
_1(11/12)$ may be similarly written in terms of
$\zeta''(0,1/12)+\zeta''(0,11/12)$, $\Gamma(1/3)$, $\Gamma(1/4)$,
$\gamma_1$, $\gamma$ and elementary functions, see
e.g.~{\eqref{oxi2103ikdmdf}}.

\section{Some results from the theory of finite Fourier series.
Applications to certain summations involving
the $\Psi$-function and the Hurwitz $\zeta$-function}\label{lo109sj1s2m2w}

\subsection{Theoretical part}\label{ljhif938hlfws}

Finite Fourier series are well-known and widely used in discrete
mathematics, numerical analysis, engineering sciences (especially in
signal and image processing) and in a lot of related disciplines.
Unlike usual Fourier series, which are essentially variants or
particular cases of the same formula, finite Fourier series may take
quite different forms and expressions. For instance, in engineering
sciences, one usually deals with the following $2m$-points Fourier
series
\begin{eqnarray*}
f_m(r) = \frac{a_m(0)}{2} + \sum_{l=1}^{m-1}
\left(a_m(l) \cos \frac{\pi rl}{m} +b_m(l) \sin
\frac{\pi rl}{m} \right) + (-1)^r\frac{a_m(m)}{2}
\end{eqnarray*}
with $r=0, 1, 2,\ldots, 2m-1$ and
$m\in\mathbbm{N}$. Thanks to the orthogonality properties of circular
functions, one may determine the coefficients in this expansion:
\begin{eqnarray*}
\everymath{\displaystyle}
\begin{cases}
 a_m(k) = \frac{1}{m} \sum_{r=1}^{2m-1} f_m(r)\cos
\frac{\pi rk}{m} ,
& k=0, 1, 2,\ldots, m \vspace{8pt}\cr
 b_m(k) = \frac{1}{m} \sum_{r=1}^{2m-1} f_m(r)\sin
\frac{\pi rk}{m} ,
& k= 1, 2, 3,\ldots, m-1
\end{cases}
\end{eqnarray*}
as well as derive Parseval's theorem
\begin{eqnarray}
\nonumber
\frac{1}{m} \sum_{r=1}^{2m-1}f^2_m(r)
= \frac{a^2_m(0)}{2} + \sum_{l=1}^{m-1}
\left(a^2_m(l) + b^2_m(l)
\right) + \frac
{a^2_m(m)}{2} ,
\end{eqnarray}
see for more details \cite[Chapter 6]{hamming_01}.

In contrast, in our researches, we encounter the following
$(m-1)$-points finite Fourier series
%
\begin{eqnarray}
\label{lkwcjh3098j} f_m(r) = a_m(0) + \sum
_{l=1}^{m-1} \left(a_m(l) \cos
\frac{2\pi rl}{m} +b_m(l) \sin\frac{2\pi rl}{m} \right)
\end{eqnarray}
$r=1, 2, 3,\ldots, m-1$, $m\in\mathbbm{N}$,
for which inversion formul\ae~and Parseval's theorem are quite
different. Let, first, derive
the inversion formul\ae~for the coefficients of this series.
Multiplying both sides
by $\cos(2\pi rk/m)$, where $k=1, 2, 3,\ldots, m-1$, and summing over
$r\in[1,m-1]$, gives
\be\label{lkjcx20djndmednsd}
\begin{array}{l}
\displaystyle
\displaystyle\sum_{r=1}^{m-1} f_m(r)\cos\dfrac{2\pi rk}{m} \,=\,  
\sum_{r=1}^{m-1} \left[ a_m(0)
+ \sum_{l=1}^{m-1} a_m(l) \cos\dfrac{2\pi rl}{m} 
+ \sum_{l=1}^{m-1} b_m(l) \sin\frac{\,2\pi r l\,}{m} \right]\cos\dfrac{2\pi rk}{m}  \\[6mm]
\quad \displaystyle=\,  a_m(0) 
\underbrace{\sum_{r=1}^{m-1} \cos\dfrac{2\pi rk}{m}}_{-1} 
+ \sum_{l=1}^{m-1} a_m(l) 
\underbrace{\sum_{r=1}^{m-1} \cos\dfrac{2\pi rl}{m} \cdot \cos\dfrac{2\pi rk}{m}}_{\frac{1}{2}m(\delta_{l,k}+\delta_{l,m-k})-1} \\[6mm]
\displaystyle \qquad
+ \sum_{l=1}^{m-1} b_m(l) 
\underbrace{\sum_{r=1}^{m-1} \sin\frac{\,2\pi r l\,}{m}\cdot\cos\dfrac{2\pi rk}{m} }_{0} 
 =\, -a_m(0)  - \sum_{l=1}^{m-1} a_m(l) + \frac{m}{2}\Big\{a_m(k) + a_m(m-k)\!\Big\}
\end{array}
\ee
Similarly, multiplying both sides of {\eqref{lkwcjh3098j}}
by $\sin(2\pi rk/m)$, where $k=1, 2, 3,\ldots, m-1$, and summing over
$r\in[1,m-1]$, yields
\be\label{hiuylckja2kdd}
\begin{array}{l}
\displaystyle
\displaystyle\sum_{r=1}^{m-1} f_m(r)\sin\dfrac{2\pi rk}{m} \,=\,  
\sum_{r=1}^{m-1} \left[ a_m(0)
+ \sum_{l=1}^{m-1} a_m(l) \cos\dfrac{2\pi rl}{m} 
+ \sum_{l=1}^{m-1} b_m(l) \sin\frac{\,2\pi r l\,}{m} \right]\sin\dfrac{2\pi rk}{m}  \\[6mm]
\qquad \displaystyle=\,  a_m(0) 
\underbrace{\sum_{r=1}^{m-1} \sin\dfrac{2\pi rk}{m}}_{0} 
+ \sum_{l=1}^{m-1} a_m(l) 
\underbrace{\sum_{r=1}^{m-1} \cos\dfrac{2\pi rl}{m} \cdot \sin\dfrac{2\pi rk}{m}}_{0}\\[6mm]
\displaystyle \qquad\quad
+ \sum_{l=1}^{m-1} b_m(l) 
\underbrace{\sum_{r=1}^{m-1} \sin\frac{\,2\pi r l\,}{m}\cdot\sin\dfrac{2\pi rk}{m} }_{\frac{1}{2}m(\delta_{l,k}-\delta_{l,m-k})} 
 =\,  \frac{m}{2}\Big\{b_m(k) - b_m(m-k)\Big\}
\end{array}
\ee
Finally, Parseval's equality for the finite series {\eqref{lkwcjh3098j}} reads:
\addtocounter{equation}{-1}
\be\label{hhxqwuybs}
\specialnumber{b}
\begin{array}{l}
\displaystyle\sum_{r=1}^{m-1} f^2_m(r)=\,  
\sum_{r=1}^{m-1} \left[ a_m(0)
+ \sum_{l=1}^{m-1} a_m(l) \cos\dfrac{2\pi rl}{m} 
+ \sum_{l=1}^{m-1} b_m(l) \sin\frac{\,2\pi r l\,}{m} \right]^2 \!\!=\\[6mm]
\quad \displaystyle
=\,  \sum_{r=1}^{m-1}  a^2_m(0) 
+ 2a_m(0) \sum_{l=1}^{m-1} a_m(l) \!
\underbrace{\sum_{r=1}^{m-1} \cos\dfrac{2\pi rl}{m}}_{-1}
+ 2a_m(0) \sum_{l=1}^{m-1} b_m(l) \!
\underbrace{\sum_{r=1}^{m-1} \sin\dfrac{2\pi rl}{m}}_{0} \\[6mm]
\qquad \displaystyle
+ 2\!\sum_{l=1}^{m-1}\sum_{n=1}^{m-1} a_m(l) b_m(n) 
\underbrace{\sum_{r=1}^{m-1} \cos\frac{\,2\pi r l\,}{m}\cdot\sin\dfrac{2\pi rn}{m} }_{0} 
\\[6mm]
\displaystyle 
\qquad 
+ \sum_{l=1}^{m-1} \sum_{n=1}^{m-1}  a_m(l) a_m(n) 
\underbrace{\sum_{r=1}^{m-1} \cos\dfrac{2\pi rl}{m} \cdot \cos\dfrac{2\pi rn}{m}}_{\frac{1}{2}m(\delta_{l,n}+\delta_{l,m-n})-1}\\[6mm]
\displaystyle \qquad
+ \sum_{l=1}^{m-1} \sum_{n=1}^{m-1}  b_m(l) b_m(n) 
\underbrace{\sum_{r=1}^{m-1} \sin\dfrac{2\pi rl}{m} \cdot \sin\dfrac{2\pi rn}{m}}_{\frac{1}{2}m(\delta_{l,n}-\delta_{l,m-n})} 
= \,(m-1)a^2_m(0) \,-\, 2a_m(0)\!\sum_{l=1}^{m-1} a_m(l) \\[6mm]
\displaystyle \qquad\quad
- \left[ \sum_{l=1}^{m-1} a_m(l)\right]^2
\!\! +\, \frac{m}{2}\sum_{l=1}^{m-1} \Big[a^2_m(l) + a_m(l)a_m(m-l) + b^2_m(l) - b_m(l)b_m(m-l) \Big]
\end{array}
\ee

\subsection{Some applications}

The finite Fourier series may be successfully used for the
finite-length summations in a variety of problems and contexts.
Consider, for example, the Gauss' Digamma theorem, which is usually
written in one of three equivalent forms
\be\label{dh38239djws}
\specialnumber{a,b,c}
\left\{
\begin{array}{ll}
\displaystyle
\Psi \biggl(\!\frac{r}{m} \vphantom{\frac{1}{2}} \!\biggr) =\,-\gamma-\ln2m-\frac{\pi}{2}\ctg\frac{\,\pi r\,}{m} 
+ 2\!\!\!\!\!\!\sum_{l=1}^{\lfloor\!\frac{1}{2}(m-1)\!\rfloor} \!\!\!\!\cos\frac{\,2\pi r l\,}{m} \cdot\ln\sin\frac{\pi l}{m}    \\[6mm]
\displaystyle
\Psi \biggl(\!\frac{r}{m} \vphantom{\frac{1}{2}} \!\biggr) =\,-\gamma-\ln2m-\frac{\pi}{2}\ctg\frac{\,\pi r\,}{m} 
+ \sum_{l=1}^{m-1} \cos\frac{\,2\pi r l\,}{m} \cdot\ln\sin\frac{\pi l}{m} \\[6mm]
\displaystyle
\Psi \biggl(\!\frac{r}{m} \vphantom{\frac{1}{2}} \!\biggr) =\,-\gamma-\ln 2\pi m-\frac{\pi}{2}\ctg\frac{\,\pi r\,}{m} - 2\sum_{l=1}^{m-1} 
\cos\frac{\,2\pi r l\,}{m} \cdot\ln\Gamma\!\left(\!\dfrac{l}{m}\!\right)
\end{array}
\right.
\ee
$r=1, 2,\ldots, m-1 $, $m\in\mathbbm{N}_{\geqslant2}$,
first and second of which are due to Gauss\footnote{Strictly speaking,
Gauss wrote them in a slightly different manner,
see \cite[p.~39]{gauss_02}.} \cite[pp.~35--38]{campbell_01}, \cite[vol.~I, p.~19, \S1.7.3]{bateman_01},
while the third one is due to Malmsten \cite[p.~57,
Eq.~(70)]{malmsten_00}, \cite[p.~37, Eq.~(23)]{iaroslav_06}.
Remarking that the cotangent may be represented by {\eqref
{jk90d2dmekwe}}, two latter equations take the
form
%
%
%
\begin{eqnarray}
\everymath{\displaystyle} \label{hebv38ls2w} 
\begin{cases}
\Psi \left( \frac{r}{m}   \right) = -\gamma
-\ln2 m +
\frac{\pi}{m} \sum_{l=1}^{m-1} \sin\frac{2\pi rl}{m} \cdot l
+ \sum_{l=1}^{m-1} \cos\frac{ 2\pi r l }{m} \cdot\ln\sin\frac
{\pi l}{m} \vspace{12pt}\cr
\Psi \left( \frac{r}{m}   \right) = -\gamma
-\ln2\pi m +
\frac{\pi}{m} \sum_{l=1}^{m-1} \sin\frac{2\pi rl}{m} \cdot l
- 2\!\sum_{l=1}^{m-1} \cos\frac{ 2\pi r l }{m} \cdot\ln\Gamma
\left( \frac{l}{m}  \right)
\end{cases}
\end{eqnarray}
$r=1, 2,\ldots, m-1 $, $m\in\mathbbm{N}_{\geqslant2}$,
which represent complete finite Fourier series
of the same type as {\eqref{lkwcjh3098j}}. Hence,
the application of {\eqref{lkjcx20djndmednsd}}--(\ref{hhxqwuybs}b)
straightforwardly yields
the following important summation formul\ae
%
\begin{eqnarray}
\label{cj932230fdncdbds}
\everymath{\displaystyle}
\begin{cases}
\sum_{r=1}^{m-1} \Psi \left( \frac{r}{m}
  \right)
\cdot\cos\frac{2\pi rk}{m} = m\ln \left( 2\sin\frac{ k\pi}{m}
 \right) + \gamma ,
\quad  k=1, 2,\ldots, m-1 \vspace{8pt}\cr
\sum_{r=1}^{m-1} \Psi \left( \frac{r}{m} \right)
\cdot\sin\frac{2\pi rk}{m} = \frac{\pi}{2} (2k-m) ,
\quad  k=1, 2,\ldots, m-1 \vspace
{8pt}\cr
\sum_{r=1}^{m-1} \Psi^2  \left( \frac{r}{m} \right) =
(m-1)\gamma^2 + m(2\gamma+\ln4m)\ln{m} -m(m-1)\ln^2 2 \vspace{2pt}\cr
\phantom{\sum_{r=1}^{m-1} \Psi^2  \left( \frac{r}{m} \right) =}
+\frac{\pi^2 (m^2-3m+2)}{12}
+m \sum_{l=1}^{ m-1 } \ln^2 \sin\frac{\pi l}{m}
\end{cases}
\end{eqnarray}
where the last sum, due to the symmetry of $\ln\sin(\pi l/m)$ about
$l=m/2$, may be also written as
\begin{eqnarray*}
\sum_{l=1}^{ m-1 } \ln^2 \sin
\frac{\pi l}{m} = 2 \!\!\!\!\!\! \sum_{l=1}^{\lfloor\frac{1}{2}(m-1) \rfloor}
\!\!\!\!\!\! \ln^2 \sin\frac
{\pi l}{m}
\end{eqnarray*}
For the purpose of demonstration, we take Malmsten's representation
for the $\Psi$-function.\footnote{The reader
may perform the same procedure with the more usual Gauss'
representation as an exercise.}
Inserting expressions for coefficients $a_m(0)=-\gamma-\ln2\pi m$,
$a_m(l)=-2\ln\Gamma(l/m)$ and $b_m(l)=\pi l/m$ into {\eqref
{lkjcx20djndmednsd}},
yields for the first sum:
\begin{eqnarray*}
\sum_{r=1}^{m-1} \Psi \left(
\frac{r}{m} \right) \cdot\cos\frac{2\pi rk}{m} = \gamma+ \ln2\pi
m - m \underbrace{ \left[\ln\Gamma \left( \frac{k}{m}
\right) + \ln\Gamma \left( 1-\frac{k}{m} \right)
\right]}_{\ln\pi-\ln\sin(\pi k/m)} \\
\qquad\quad + \, 2 \!\sum
_{l=1}^{m-1} \ln\Gamma \left(
\frac{l}{m} \right)
= \gamma+ m\ln \left( 2\sin\frac{ k\pi}{m} \right)
\end{eqnarray*}
where the final simplification is performed with the help of the
reflection formula and
Gauss' multiplication theorem for the logarithm of the $\Gamma$-function {\eqref{kjd203jxndhe8}}.
Analogously, using {\eqref{hiuylckja2kdd}} yields for the second sum:
\begin{eqnarray*}
\sum_{r=1}^{m-1} \Psi \left(
\frac{r}{m} \right) \cdot\sin \frac{2\pi rk}{m} \ =
\frac{m}{2} \left[ \frac{\pi k}{m}- \frac{\pi(m-k)}{m}
\right] = \frac{\pi}{2} (2k-m)
\end{eqnarray*}
By taking advantage of this opportunity, we would like to remark that a
formula of the similar nature appears also in \cite[p.~39]{campbell_01} and
\cite[p.~19, Eq.~(49)]{srivastava_03}. Sadly, the formula given in the
former source contains two errors;
the correct variant of the formula is
\begin{eqnarray*}
\sum_{r=1}^{m} \Psi \left(
\frac{r}{m} \right) \cdot\exp \frac{2\pi rk i}{m} = m\ln
\left( 1-\exp\frac{2\pi k
i}{m} \right) , \quad k\in\mathbbm{Z} , \ m
\in\mathbbm{N} , \ k\neq m.
\end{eqnarray*}

Finally, by formula (\ref{hhxqwuybs}b), we derive Parseval's theorem for
the $\Psi$-function of a discrete argument
\be\label{kjch328hdcnls}
\begin{array}{rl}
\displaystyle\sum_{r=1}^{m-1} \Psi^2\! \biggl(\!\frac{r}{m} \vphantom{\frac{1}{2}} \!\biggr) = \,  
(m-1)(\gamma + \ln 2\pi m)^2 - 4 (\gamma + \ln 2\pi m)
\!\!\!\!\!\! \underbrace{\sum_{l=1}^{m-1} \ln\Gamma\!\left(\!\dfrac{l}{m}\!\right)}_{\frac{1}{2}(m-1)\ln2\pi-\frac{1}{2}\ln{m}}\!\!\!\!\! 
- \,4 \Biggl[\sum_{l=1}^{m-1} \ln\Gamma\!\left(\!\dfrac{l}{m}\!\right)\Biggr]^2
\\[10mm]
\displaystyle  +2m \!\sum_{l=1}^{m-1} \ln\Gamma\!\left(\!\dfrac{l}{m}\!\right) 
\cdot\Biggl[\ln\Gamma\!\left(\!\dfrac{l}{m}\!\right) + \ln\Gamma\!\left(\!1-\dfrac{l}{m}\!\right)\Biggr] 
+ \frac{\pi^2}{m} \cdot\!\sum_{l=1}^{m-1} l^2 
- \, \frac{\pi^2}{2}\cdot\!\sum_{l=1}^{m-1} l =  (m-1)\gamma^2  \\[5mm]
\displaystyle 
+ m(2\gamma+\ln4m)\ln{m} -m(m-1)\ln^2\!2 
+\frac{\pi^2 (m^2-3m+2)}{12} +2m\!\!\!\!\!\!\sum_{l=1}^{\lfloor\!\frac{1}{2}(m-1)\!\rfloor} \!\!\!\!
\ln^2\!\sin\frac{\pi l}{m}
\end{array}
\ee
where the sum from the third line, thanks to the symmetry of
$\ln\sin(\pi l/m)$ about $l=m/2$
and to the fact that $\ln\sin(\pi l/m)=0$ for $l=m/2$, could be
simplified as follows
\be\label{7826gdxhkwx39e4}
\begin{array}{l}
\displaystyle
\sum_{l=1}^{m-1} \ln\Gamma\!\left(\!\dfrac{l}{m}\!\right)
\cdot\Biggl[\ln\Gamma\!\left(\!\dfrac{l}{m}\!\right) + 
\ln\Gamma\!\left(\!1-\dfrac{l}{m}\!\right)\Biggr] = 
\sum_{l=1}^{m-1} \ln\Gamma\!\left(\!\dfrac{l}{m}\!\right)
\cdot\Biggl[\ln\pi-\ln\sin\frac{\pi l}{m}\Biggr] \\[6mm]
\displaystyle\quad
 = \frac{\ln\pi}{2}\Big[(m-1)\ln2\pi-\ln{m}\Big]- \sum_{l=1}^{m-1} \ln\Gamma\!\left(\!\dfrac{l}{m}\!\right) \cdot\ln\sin\frac{\pi l}{m}
= \frac{\ln\pi}{2}\Big[(m-1)\ln2\pi-\ln{m}\Big] \\[6mm]
\displaystyle\qquad
- \!\!\!\!\sum_{l=1}^{\lfloor\!\frac{1}{2}(m-1)\!\rfloor} \!
\Biggl[\ln\pi-\ln\sin\frac{\pi l}{m}\Biggr] \ln\sin\frac{\pi l}{m} = \frac{\ln\pi}{2}\Big[(m-1)\ln2\pi-\ln{m}\Big]
-\ln\pi \!\!\!\!\! \sum_{l=1}^{\lfloor\!\frac{1}{2}(m-1)\!\rfloor} \!\!\!\! \ln\sin\frac{\pi l}{m}\\[6mm]
\displaystyle\qquad
+ \!\!\!\!\!\!\sum_{l=1}^{\lfloor\!\frac{1}{2}(m-1)\!\rfloor} \!\!\!\!\ln^2\!\sin\frac{\pi l}{m} 
= \frac{\ln\pi}{2}\Big[(m-1)\ln4\pi-2\ln{m}\Big] \,+ \!\!\!\!\!\sum_{l=1}^{\lfloor\!\frac{1}{2}(m-1)\!\rfloor} \!\!\!\!\ln^2\!\sin\frac{\pi l}{m} 
\end{array}
\ee
because
\begin{eqnarray*}
\sum_{l=1}^{\lfloor\frac{1}{2}(m-1) \rfloor}\!\!\!\!\!  \ln\sin
\frac{\pi
l}{m} = \ln \!\!\!\!\!\! \prod_{l=1}^{\lfloor\frac{1}{2}(m-1) \rfloor}
 \!\!\!\!\! \sin\frac{\pi
l}{m} = \frac{1-m}{2}\ln2+\frac{1}{2}\ln m
\end{eqnarray*}
and where
%
\begin{eqnarray}
\label{jk0909jjds2eja} \sum_{l=1}^{m-1}
l^2 = \frac{ m(m-1)(2m-1) }{6} \quad\mbox{and}\quad \sum
_{l=1}^{m-1} l = \frac{ m(m-1) }{2}
\end{eqnarray}
respectively, which completes the evaluation of the third formula in
{\eqref{cj932230fdncdbds}}.

In like manner, we may also derive similar summation formul\ae~ for
the Hurwitz $\zeta$-function.
Rewriting Hurwitz's functional equation {\eqref{khjhbcubdbvxhbs}} in the
form analogous to {\eqref{lkwcjh3098j}}
\be\notag
\begin{array}{l}
\displaystyle
\zeta\!\left(a,\frac{r}{m}\right) 
=\, m^{a-1}\zeta(a) + \frac{2 \Gamma(1-a)}{(2\pi m)^{1-a}}\Biggl[ \sin\frac{\pi a}{2}
\sum_{l=1}^{m-1} \cos\frac{2\pi rl}{m} \cdot\zeta\!\left(\!1-a,\,\frac{l}{m}\!\right)  + \\[6mm]
\displaystyle\qquad\qquad\qquad\qquad\qquad\qquad\qquad\qquad\qquad\qquad
+ \cos\frac{\pi a}{2}
\sum_{l=1}^{m-1} \sin \frac{2\pi rl}{m} \cdot\zeta\!\left(\!1-a,\,\frac{l}{m}\!\right) \Biggr] 
\end{array}
\ee
yields
\be\notag
\begin{array}{l}
\displaystyle
\displaystyle\sum_{r=1}^{m-1} \zeta\!\left(a,\frac{r}{m}\right) 
\!\cdot\cos\dfrac{2\pi rk}{m} =\frac{m \Gamma(1-a)}{(2\pi m)^{1-a}}
\sin\frac{\pi a}{2} \cdot \left\{\zeta\!\left(\!1-a,\,\frac{k}{m}\!\right) +
\zeta\!\left(\!1-a,\,1-\frac{k}{m}\!\right) \!\right\}  -\, \zeta(a) \\[6mm]
\displaystyle\sum_{r=1}^{m-1} \zeta\!\left(a,\frac{r}{m}\right) 
\!\cdot\sin\dfrac{2\pi rk}{m} \,=\, \frac{m \Gamma(1-a)}{(2\pi m)^{1-a}}
\cos \frac{\pi a}{2} \cdot \left\{\zeta\!\left(\!1-a,\,\frac{k}{m}\!\right) -
\zeta\!\left(\!1-a,\,1-\frac{k}{m}\!\right) \!\right\}   \\[6mm]
\displaystyle\sum_{r=1}^{m-1} \zeta^2\!\left(a,\frac{r}{m}\right) = \,
\big(m^{2a-1}-1 \big)\zeta^2(a) + \frac{2m\Gamma^2(1-a)}{(2\pi m)^{2-2a}} \times\\
\displaystyle\qquad\qquad\qquad\qquad\qquad
\times\!
\sum_{l=1}^{m-1} \left\{ \zeta\!\left(1-a,\frac{l}{m}\right) - \cos\pi a 
\cdot \zeta\!\left(1-a,1-\frac{l}{m}\right)\!\right\}\cdot \zeta\!\left(1-a,\frac{l}{m}\right)
\end{array}
\ee
which hold for any $r=1, 2,3, \ldots, m-1$ and $k=1, 2,3, \ldots, m-1$,
where $m$ is positive integer.

By the way, there are many other functions which are orthogonal or
semi-orthogonal over some
discrete interval. For instance, by considering another set of circular
functions and their properties
\be\label{olkefj938hjcwdc}
\begin{array}{l}
\displaystyle
\sum_{r=0}^{m-1} \cos\dfrac{(2r+1)k\pi }{m} \,= 
\sum_{r=0}^{m-1} \sin\dfrac{(2r+1)k\pi }{m} \,=  0  \,, \qquad\qquad k=1, 2,\ldots, m-1\\[6mm]
\displaystyle
\sum_{r=0}^{m-1} \cos\dfrac{(2r+1)k\pi }{m} \cdot\sin\dfrac{(2r+1)l\pi }{m} \,=  0  \,,  \qquad\qquad\qquad \quad\!\! k,l=1, 2,\ldots, m-1\\[6mm]
\displaystyle
\sum_{r=0}^{m-1} \cos\dfrac{(2r+1)k\pi }{m} \cdot\cos\dfrac{(2r+1)l\pi }{m} \,=  \frac{n}{2}\big(\delta_{k,l}-\delta_{k,m-l}
- \delta_{k,m+l}+\delta_{k,2m-l}\big)   \\[6mm]
\displaystyle
\sum_{r=0}^{m-1} \sin\dfrac{(2r+1)k\pi }{m} \cdot\sin\dfrac{(2r+1)l\pi }{m} \,=   \frac{n}{2}\big(\delta_{k,l}+\delta_{k,m-l}
- \delta_{k,m+l}-\delta_{k,2m-l}\big)  
\end{array}
\ee
where in last two formul\ae~$k,l=1, 2,\ldots, 2m-1$, as well as
{\eqref{hebv38ls2w}}, one may easily prove that
\be\notag
\begin{cases}
\displaystyle\sum_{r=0}^{m-1} \Psi \biggl(\!\frac{2r+1}{2m} \vphantom{\frac{1}{2}} \!\biggr)
\!\cdot\cos\dfrac{(2r+1)k\pi }{m} \,=\, m\ln\tg\frac{\,\pi k\,}{2m} \,,
\qquad\qquad\qquad & k=1, 2,\ldots, m-1 \\[6mm]
\displaystyle\sum_{r=0}^{m-1} \Psi \biggl(\!\frac{2r+1}{2m} \vphantom{\frac{1}{2}} \!\biggr)
\!\cdot\sin\dfrac{(2r+1)k\pi }{m} \,=\,-\frac{\pi m}{2}    \,,
\qquad\qquad\qquad & k=1, 2,\ldots, m-1   
\end{cases}
\ee

By a similar line of reasoning, we also derive
\be\label{jk920u3jd209rnd}
\begin{array}{ll}
\displaystyle\sum_{r=1}^{2m-1} (-1)^r\!\cdot\Psi \biggl(\!\frac{r}{2m} \vphantom{\frac{1}{2}} \!\biggr) \,=\, 2m\ln2+\gamma\\[7mm]
\displaystyle\sum_{r=0}^{2m-1} (-1)^r\!\cdot\Psi \biggl(\!\frac{2r+1}{4m} \vphantom{\frac{1}{2}} \!\biggr) \,=\,-\pi m\\[7mm]
\displaystyle\sum_{r=1}^{m-1} \ctg\dfrac{\pi r}{m}\cdot \Psi \biggl(\!\frac{r}{m} \vphantom{\frac{1}{2}} \!\biggr)\,=\, -\frac{\pi(m-1)(m-2)}{6} \\[7mm]
\end{array}
\ee

\addtocounter{equation}{-1}
\be
\begin{array}{ll}
\displaystyle\sum_{r=1}^{m-1} \dfrac{r}{m}\cdot \Psi \biggl(\!\frac{r}{m} \vphantom{\frac{1}{2}} \!\biggr)\,=
\,-\frac{\gamma}{2}(m-1)-\frac{m}{2}\ln m -\frac{\pi}{2}\sum_{r=1}^{m-1} \dfrac{r}{m}\cdot\ctg\dfrac{\pi r}{m} \\[7mm]
\displaystyle\sum_{r=1}^{m-1} \cos\dfrac{(2l+1)\pi r}{m}\cdot \Psi \biggl(\!\frac{r}{m} \vphantom{\frac{1}{2}} \!\biggr)\,=\, 
-\frac{\pi}{m}\cdot\!\sum_{r=1}^{m-1} \frac{r \cdot\sin\dfrac{2\pi r}{m}}{\,\cos\dfrac{2\pi r}{m} -\cos\dfrac{(2l+1)\pi }{m} \,}\\[8mm]
\displaystyle\sum_{r=1}^{m-1} \sin\dfrac{(2l+1)\pi r}{m}\cdot \Psi \biggl(\!\frac{r}{m} \vphantom{\frac{1}{2}} \!\biggr)\,=\, 
-(\gamma+\ln2m)\ctg\frac{(2l+1)\pi}{2m} \\[-1mm]
\displaystyle \qquad\qquad\qquad\qquad\qquad\qquad\qquad
+ \sin\dfrac{(2l+1)\pi }{m}\cdot\!\sum_{r=1}^{m-1} \frac{\ln\sin\dfrac{\pi r}{m}}
{\,\cos\dfrac{2\pi r}{m} -\cos\dfrac{(2l+1)\pi }{m} \,}
\end{array}
\ee
where the last two formul\ae~remain valid for any $l\in\mathbbm{Z}$.

\section{An integral formula for the logarithm of the $\Gamma
$-function at rational arguments}\label{894yf3hedbe}

In this part, we evaluate integral {\eqref{kj093dhcg2387bh}} for $p=k/n$
and show that it reduces to the logarithm of the $\Gamma$-function at
rational argument, Euler's constant $\gamma$ and elementary
functions.

From a simple algebraic argument, it follows that
\begin{eqnarray*}
\sum_{r=1}^{n-1} \operatorname{sh} {rx}
\cdot\sin\frac{2\pi r k}{n} = -\frac{1}{2}\cdot \frac{\operatorname{sh}{nx}\cdot\sin\frac{2\pi k}{n}}{
\operatorname{ch}{x}-\cos\frac{2\pi k}{n} } ,
\quad x\in\mathbbm{C} , \ k\in\mathbbm{Z} .
\end{eqnarray*}
Then, for $p=k/n$, where $k$ and $n$ are positive integers such that
$k$ does not exceed $n$,
the denominator of integrand {\eqref{kj093dhcg2387bh}} may be replaced
by the above identity and hence
%
\begin{eqnarray}
\label{ascvwdi209} \int\limits_0^\infty \frac{ e^{-nx} \cdot\ln{x} }{ \operatorname
{ch}{x}-\cos\frac{2\pi k}{n} }
\,dx = -4\csc\frac{2\pi k}{n}\sum_{r=1}^{n-1}
\sin\frac{2\pi r k}{n}\cdot \!\int\limits_0^\infty
\frac{ \operatorname{sh}{rx}\cdot\ln{x} }{
e^{2nx}-1 } \,dx
\end{eqnarray}
The latter integral was already evaluated in our previous work, see
\cite[p.~73, \no25]{iaroslav_06}.
By setting in exercise \no25-a $b=n$, $m=r$, and then by rewriting the
result for $2n$ instead of $n$, we get
%
\begin{eqnarray}
\int\limits_0^\infty \frac{ \operatorname{sh}{rx}\cdot\ln{x} }{
e^{2nx}-1 } \,dx = -
\frac{\pi}{4n}\operatorname{ctg}\frac{r\pi}{2n}\cdot\ln2\pi -
\frac{\gamma+\ln{r}}{2r} +\frac{\pi}{2n} \sum_{l=1}^{2n-1}
\sin\frac{\pi r l}{n}\cdot\ln \Gamma \left( \frac{l}{2n}
\right)\nonumber\\
\end{eqnarray}
By inserting the above formula into {\eqref{ascvwdi209}} and by
recalling that for $k=1, 2, 3,\ldots, n-1$
and $l=1, 2, 3,\ldots, 2n-1$
\begin{eqnarray*}
\everymath{\displaystyle}
\begin{cases}
\sum_{r=1}^{n-1} \sin\frac{2\pi r k}{n} \cdot\operatorname
{ctg}\frac{\pi r}{2n} = n-2k \vspace{6pt}\cr
\sum_{r=1}^{n-1} \sin\frac{2\pi r k}{n} \cdot\sin\frac{\pi rl}{n}
= \frac{n}{2}
 \{\delta_{k,\frac{l}{2}} - \delta_{k,n-\frac{l}{2}}  \}
\end{cases}
\end{eqnarray*}
the expression for integral {\eqref{kj093dhcg2387bh}} at
$p=k/n$ takes its final form
\be\label{kjh298hdnd}
\begin{array}{l}
\displaystyle 
\int\limits_0^\infty \!\!\frac{\,e^{-nx}\!\cdot\ln{x}\,}{\,\ch{x}-\cos\frac{2\pi k}{n}\,}\,dx \,=  \,
2\!\int\limits_0^1\!\!\frac{\,x^n\ln{\ln{\frac{1}{x}}}\,}{\,x^2-2x\cos\frac{2\pi k}{n}+1\,}\,dx\,=\,
2\!\int\limits_1^\infty \!\!\frac{\,\ln{\ln{x}}\,}{\,x^n\!\left(x^2-2x\cos\frac{2\pi k}{n}+1\right)\,}\,dx
\\[8mm]
\displaystyle \quad
=\left\{\frac{\pi(n-2k)\ln2\pi}{n}-2\pi\ln\Gamma 
\biggl(\!\frac{k}{n}\vphantom{\frac{1}{2}} \!\biggr) 
+\pi\ln{\pi}-\pi\ln\sin\frac{\pi k}{n} 
+ 2\!\sum_{r=1}^{n-1}\frac{\gamma+\ln{r}}{r}\cdot\sin\frac{2\pi r k}{n} \right\}\times\\[5mm]
\displaystyle\qquad\;
\times\csc\frac{2\pi k}{n}
\end{array}
\ee
Whence
\begin{eqnarray}
\label{kjh298hdnd2} %
\ln\Gamma \left(
\frac{k}{n} \right) &=& \frac{(n-2k)\ln2\pi}{2n} + \frac{1}{2}\left\{\ln{\pi}-\ln\sin\frac{\pi k}{n}\right\}
+\frac{1}{\pi}\sum_{r=1}^{n-1}\frac{\gamma+\ln{r}}{r}\cdot\sin
\frac{2\pi r k}{n} \nonumber\\[6pt]
&& -\frac{1}{2\pi}\sin\frac{2\pi k}{n}\cdot\!
\int\limits_0^\infty
\frac{ e^{-nx} \cdot\ln{x} }{ \operatorname
{ch}{x}-\cos\frac{2\pi k}{n} } \,dx\,,\qquad k=1, 2, 3,\ldots, n-1\,,\qquad
\end{eqnarray}
$k\neq n/2$.
By the way, {\eqref{kjh298hdnd}}--{\eqref{kjh298hdnd2}} may be proven by other methods
as well. For instance, one may
directly employ {\eqref{908djdsn}}
because $a_0=0$ for $p=k/n$ and all remaining integrals in the
right-hand side are known. Yet, \eqref{kjh298hdnd}--\eqref{kjh298hdnd2} may be also
obtained with the aid of previously derived results in exercises \no
60 and 58 in \cite[Sect.~4]{iaroslav_06}, as well as Malmsten's
representation
for the logarithm of the $\Gamma$-function
%
\begin{eqnarray}
\label{j2039drjm2d} \ln\Gamma(z)= \frac{1}{2}\ln\pi - \frac{1}{2}
\ln\sin\pi z - \frac{2z-1}{2}\ln2\pi - \frac{\sin2\pi z}{2\pi} \int\limits
_0^\infty \frac{ \ln{x} }{ \mathrm{ch} {x}-\cos2\pi z } \,dx
\end{eqnarray}
where $0<\operatorname{Re}z<1$,
see exercises \no2, 29-h, 30 \cite[Sect.~4]{iaroslav_06}.

\vspace{19mm}


\begin{figure}[!h]   
\centering
\includegraphics[width=0.3\textwidth]{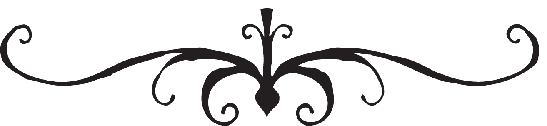}
\end{figure}

\newpage

\bibliographystyle{crelle}

\end{document}